\documentclass[12pt,final,
titlepage]{article}
\usepackage{tikz}
\usepackage{epsfig}
\usepackage{amsfonts}
\usepackage{amssymb}
\usepackage{amsmath}
\usepackage{amsthm}
\usepackage{latexsym}
\usepackage{graphicx}
\usepackage{bm}
\usepackage{tabularx}
\usepackage{booktabs}
\usepackage{enumitem}
\usepackage{multirow}
\usepackage[titletoc]{appendix}
\usepackage[numbers]{natbib}
\usepackage{caption}
\usepackage{subcaption}
\usepackage{mathtools}
\usepackage{multirow}
\usepackage{lineno}
\usepackage{atbegshi}      % hook at shipout (page output)

\newcommand{\field}[1]{\mathbb{#1}}
\def\A{\field{A}}                                % Action space
\def\X{\field{X}}                                % State space

\def\Z{\field{Z}}

\catcode`\@=11

\numberwithin{equation}{section}
\numberwithin{table}{section}
\numberwithin{figure}{section}

\theoremstyle{plain}
\newtheorem{theorem}{Theorem}[section]
\newtheorem{lemma}[theorem]{Lemma}

\newtheorem{proposition}[theorem]{Proposition}

\theoremstyle{definition}
\newtheorem{definition}{Definition}[section]
\newtheorem{assumption}[definition]{Assumption}

\theoremstyle{remark}

\newcommand{\ls}[1]
  {\dimen0=\fontdimen6\the\font \lineskip=#1\dimen0
  \advance\lineskip.5\fontdimen5\the\font \advance\lineskip-\dimen0
  \lineskiplimit=.9\lineskip \baselineskip=\lineskip
  \advance\baselineskip\dimen0 \normallineskip\lineskip
  \normallineskiplimit\lineskiplimit \normalbaselineskip\baselineskip
  \ignorespaces }

\graphicspath{{./}{Figs/}}

\usepackage[colorlinks=true,linkcolor=magenta,citecolor=blue,pagebackref]{hyperref}
\usepackage{authblk}

\title{Control policies for a two-stage queueing system with parallel and single server options}
\author[1]{Shuwen Lu}
\author[2]{Mark E. Lewis}
\author[3] {Jamol Pender}%\\ \\

\affil[1]{\small Department of Systems Engineering, Cornell University, Ithaca, NY 14850, USA {\tt\small sl3243@cornell.edu}}%
\affil[2,3]{\small School of Operations Research and Information Engineering,
        Cornell University, Ithaca, NY 14850, USA
        {\tt\small mark.lewis@cornell.edu, jjp274@cornell.edu}}

\begin{document}
% \pagewiselinenumbers  % resets to 1 each page
\maketitle
\ls{1.25}
\begin{abstract}
{\color{blue} We study a two-stage queueing system attended by two flexible servers. Once a server begins to interact with a job, it remains attached to its job throughout both stages of service. After completing the first-stage service, the server chooses between two second-stage options. Option 1 sends the server-job pair to a parallel-service station where both servers may operate simultaneously; therefore, the server immediately begins the second-stage service on its job. Option 2 routes the pair to a single-service station that admits at most one server-job pair into service at a time. While we consider the complete range of parameters, the usual reason for choosing Option 2 is that jobs may receive higher-quality (faster and/or lower expected cost) service. If the single-service station is occupied, the pair must wait together until the station becomes available.}

The decision-making scenario is modeled via a Markov decision process formulation of a clearing system with holding costs at each station. \textcolor{blue}{Focusing on cases with a large initial queue,} we fully characterize the structural properties of an optimal control policy based on the relationship between the service rates at the downstream stations. A numerical study highlights the significance of \textcolor{blue}{understanding the structure of the optimal policy} by comparing its performance against several natural heuristic policies.

\end{abstract}

\section{Introduction}\label{sec:intro}
\begin{figure}[htbp]
\centering
    \includegraphics[scale=0.5]{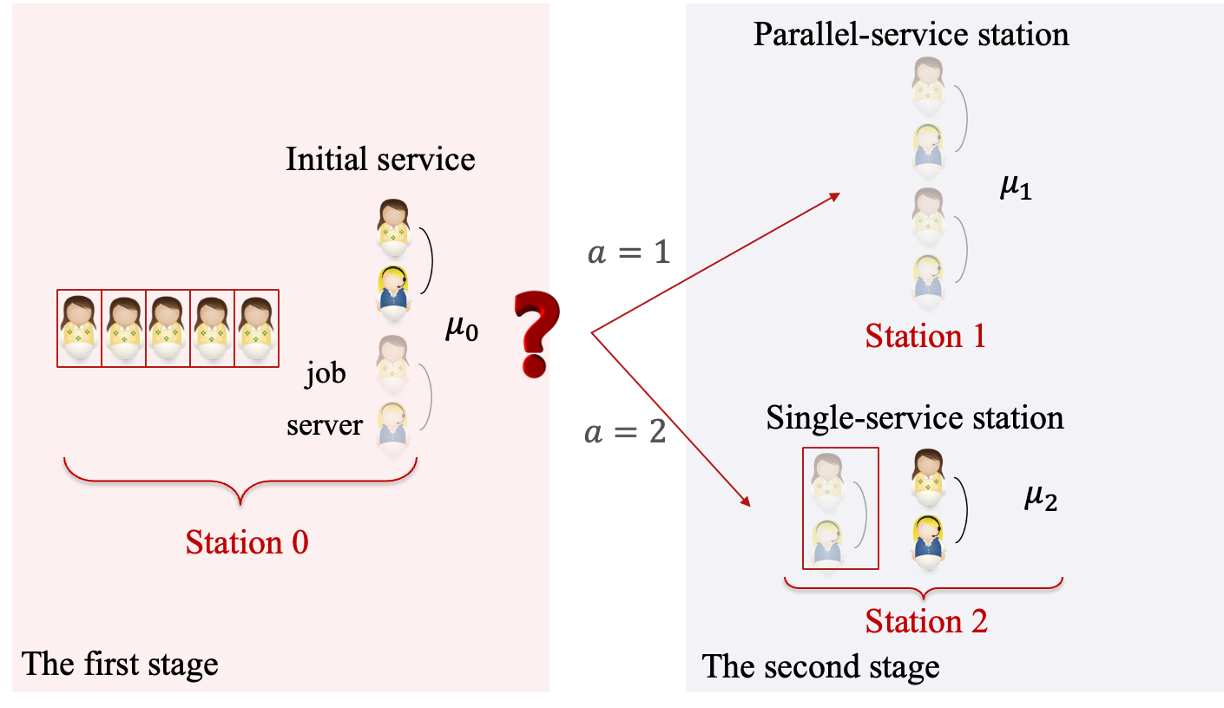}
\caption{\textcolor{blue}{\textbf{Diagram of Two-stage Model.} Queues are represented by rectangles. Server-job pairs are depicted by arcs connecting each job to its job in service. Opaque server-job pairs indicate current work assignments, while semi-transparent pairs denote potential future assignments. The single-service station allows only one server to process its job at a time.}} \label{Fig:two-stage}
\end{figure}

{\color{blue}In this paper, we study a two-stage service system in which each job requires two sequential stages of service. The system is staffed by two fully flexible servers, each capable of performing either stage. 
During the first stage, the servers conduct an initial screening and may process their respective jobs in parallel. Upon completing the first-stage service, the server must choose between two downstream service options and proceed to the second stage along with the job. Under the first option, the server moves with its job to a parallel-service station, allowing concurrent processing so the second-stage service commences immediately. Under the second option, the server-job pair enters a single-service station that may provide higher-quality (faster and/or lower expected cost) service, but allows at most one server to operate at a time. If this station is occupied, sending another pair induces downstream queueing, delaying the server’s return to the first stage. 
We assume that either second-stage option fully completes the job, after which the server becomes available to return to the first stage and begin the initial screening on a new job. The queueing dynamics are illustrated in Figure \ref{Fig:two-stage}.
}

\textcolor{blue}{Because the server makes the routing decision only after completing first-stage service and must process the same job at the selected second-stage station, upstream and downstream operations interact in a nontrivial way. This complicates the decision process. In particular, if routing a server-job pair to the single-service station leads to waiting, the server is blocked while accompanying the job, leaving the first-stage queue unattended and delaying subsequent processing. The impact of these blocking complications depends critically on the relative speeds of upstream and downstream services. When the first-stage service is slow, downstream services are likely to complete before the next decision point. The single-service station is often available at that time. In contrast, when the first-stage service is fast, decisions occur before the single-service station becomes available, thereby compounding the upstream congestion effects. As a result, simple structural properties, such as monotone or threshold-type policies, are unlikely to hold, underscoring the technical challenges.} 

\textcolor{blue}{ 
At a higher level, the system captures a broader class of multi-server models in which servers choose between immediate processing without resource constraints and (potentially) waiting to access a valuable but limited resource.
In our setting, relative to the parallel-service option, the effective access to this resource per server under the alternative option is reduced by half.}

\textcolor{blue}{Such queueing dynamics with control arise in clinical telehealth environments. For example, following an initial diagnosis, frontline providers decide whether to treat low-acuity patients independently or request a teleconference with a general physician for a more comprehensive assessment. Physician involvement can increase patient confidence but may introduce delays due to limited physician availability or constrained telemedicine infrastructure. In the escalated pathway, the frontline provider stays with the patient so that the consultation can begin as soon as a physician becomes available \cite{richmond2017american}, and remains engaged during the session to provide essential assistance and support physician-nurse collaboration \cite{cherba2022essential}. 
Similar settings appear in the manufacturing industry, where engineers, after an initial screening, choose between processing the job immediately using a readily available standard tool and sending it to a specialized machine that offers higher yield or throughput. 
These machines are often scarce or capacity-constrained, thereby limiting the number of concurrent jobs. Because operating them requires dedicated supervision and involves costly setup or changeover, the initiating engineer typically remains assigned to the job and waits with it until the machine becomes available. 
In these cases, limited but higher-quality resources, such as expert personnel or specialized equipment, drive the need for strategic decision-making.}

\textcolor{blue}{Since job-completion requirements and demand bursts frequently arise in practice \cite{argon2008scheduling, minh2011performance}, we focus on how the system empties an existing backlog or returns to its normative state following a surge in arrivals via a clearing system model. Focusing on regimes with a large initial number of queued jobs, we fully characterize the optimal policy in terms of the relative service rates of the two second-stage options. 
In particular, we identify sufficient conditions under which a large initial queue makes one downstream option preferable to the other.}

Additionally, this work provides interpretable and practical insights into tandem queues with a decision point after the initial stage, where one option offers immediate service while the other \textcolor{blue}{offers potentially faster or more cost-efficient service but is capacity-limited and may create blocking. 
In such settings, effective decisions must balance the immediate benefit to the current job against the future congestion imposed on subsequent jobs. Our structural results clarify this trade-off by showing how the optimal routing policy depends on the system state and key parameters. Finally, our numerical study illustrates how the structure of the optimal policy varies with system parameters and evaluates several simple heuristics against the optimum. The results show when simple rules perform well and when more careful state-dependent control is needed.}

The rest of the paper is organized as follows: Section \ref{sec:lit-review} provides a review of relevant literature. The details of the model are covered in Section \ref{sec:description}. A statement of the main results is in Section \ref{sec:main-results}. 
A numerical study in Section \ref{sec:num} shows the complexity of the optimal policy and compares it against several heuristic policies.
Finally, a conclusion is drawn in Section \ref{sec:conc} along with several new directions for future research.
The proofs of the main results showing the structure of \textcolor{blue}{optimal policies} are provided in \textcolor{blue}{Appendix} \ref{sec:appendix}. 

% \begin{figure}[htbp]
% \centering
%     \includegraphics[scale=0.5]{two_stage_queue.png}
% \caption{\textbf{Diagram of Two-stage Model.} Queues are represented by rectangles. Current servers are depicted as arcs connecting to jobs in service. Opaque servers indicate where servers may work in the future.} \label{Fig:two-stage}
% \end{figure}

\section{\textcolor{blue}{Literature review}}\label{sec:lit-review}

There has been extensive research on the optimal control of tandem queues, \textcolor{black}{with much of the classical literature focusing on server allocation problems aimed at minimizing holding and related operational costs.} In a single-server tandem queue, the server is entirely flexible and must complete services in all stages. 
\textcolor{black}{For two-stage systems with Poisson arrivals and general service and switching times, it is known that when the holding cost rate at the second stage exceeds that at the first stage, the optimal policy is greedy and exhaustive at the second stage \cite{iravani1997two}. Moreover, when the server is configured for the first stage, the optimal policy is characterized by a switching curve \cite{duenyas1998control}. 
% In contrast, fully characterizing optimal policies in service systems with three or more stages remains considerably more challenging \cite{duenyas1998control}.
}

\textcolor{black}{Conventional tandem queues often assume multiple flexible servers that can serve jobs at different stages. For example, two-stage tandem systems with two parallel flexible servers have been studied in \cite{ahn1999optimal, schiefermayr2005complete, ahn2002optimal}, and extensions to more servers appear in \cite{pandelis1994optimal}, primarily with the objective of minimizing holding costs.}
The case where two identical servers operate in parallel in the sense that either can process any job is similar to our model in that aspect. Given an initial number of jobs in the two stages without further external arrivals, the necessary and sufficient conditions under which an optimal control policy allocates both servers to upstream or downstream queues are characterized in \cite{ahn1999optimal}, if no preemption is allowed. If instead assuming one job can be preempted once the other is finished so that both servers can serve the same job, a complete characterization of the optimal control under all possible parameters is presented in \cite{schiefermayr2005complete}. With arrivals, the optimal policy is studied in \cite{ahn2002optimal} in two scenarios based on whether the two servers can collaboratively serve one job under the assumption of preemption. 
\textcolor{blue}{In two-stage queues with more parallel flexible servers, general service times and the possibility of departing the system without entering the downstream stage have been considered \cite{pandelis1994optimal}.} 
\textcolor{black}{Within the class of non-preemptive policies, conditions are established under which prioritizing the upstream queue is optimal.}

% There is also extensive research on server allocation policies in series queues with more complex service structures (e.g., the existence of dedicated servers) \cite{farrar1993optimal, wu2006dynamic, pandelis2008optimal} and/or with different objectives such as maximizing the throughput \cite{andradottir2001server, wu2008heuristics}, maximizing rewards with penalty for waiting or abandonment \cite{andradottir2021optimizing, yu2023optimal}, etc.
\textcolor{blue}{When dedicated servers are present, allocation policies for flexible servers in two-stage queues have also been studied, with the objective of minimizing holding costs. In a system with one dedicated server at each stage and one flexible server, the existence of an optimal transition-monotone policy is established in \cite{farrar1993optimal}.} 
The existence of these structured policies is then extended by \cite{wu2006dynamic} with multiple dedicated and re-configurable servers. Further extension is in \cite{pandelis2008optimal}, where the assumption that two stations are identical is relaxed, and all servers incur different operating costs while working at different stations.  
\textcolor{blue}{Other literature considers alternative objectives, including reward maximization \cite{zayas2016dynamic} and throughput maximization studied in \cite{andradottir2001server} without dedicated servers and in \cite{wu2008heuristics} with dedicated servers, among others.} 

\textcolor{blue}{Another stream of research on the control of tandem queues focuses on work prioritization policies involving multiple types of servers, each capable of operating only in specific stages or modes. 
For example, the interaction between attending physicians and residents is modeled as a three-stage tandem queue in \cite{dobson2012queueing}, where throughput-optimal policies are derived under the assumption of at most two residents.}
To study the time allocation of attending physicians between meeting with residents and their own responsibilities, a two-stage service system for such interaction \textcolor{blue}{is considered in \cite{andradottir2021optimizing}}. The analysis allows an arbitrary number of residents with the objective of maximizing the long-term reward, where the throughput-optimal objective is a special case. Related work by \cite{yu2023optimal}, by comparison, considers multiple supervisors (resp. attending physicians in \cite{andradottir2021optimizing}), as opposed to only one in \cite{andradottir2021optimizing}, customer abandonments and a different cost structure.

\textcolor{black}{In systems with parallel service options, the main control is often routing/assignment, dispatching each job to a service type or station \cite{hordijk1992assignment, down2006dynamic}. 
In settings where multiple job classes share a single server (or a pooled capacity), the control is typically scheduling, i.e., selecting which class of job receives service next \cite{nain1989interchange, huang2022dynamically}.}
In a G/M/1 queue, it is proved that it is optimal to schedule according to the $c\mu$-rule \cite{nain1989interchange}. 
In the presence of server deterioration, sufficient conditions of the optimality of a state-dependent $c\mu$-rule \textcolor{blue}{are identified in \cite{huang2022dynamically}}. 
\textcolor{blue}{The N-network is a queueing system with two job classes and two server pools, where the dedicated pool serves only class-1 jobs while the flexible pool can serve both class-1 and class-2 jobs \cite{ahn2004optimal, bell2001dynamic}. The service structure can be represented as a bipartite graph in the shape of an “N,” hence the name. Control typically focuses on scheduling the flexible pool, determining when to prioritize class-2 jobs versus when to assist the dedicated pool with class-1 jobs.}

\textcolor{blue}{Despite its tandem structure and a routing decision, our model differs substantially from the existing queueing-control literature. From a decision-making perspective, the server-job pair chooses between two distinct service options: a parallel-service station that guarantees immediate processing and a single-service station that offers higher service quality \textcolor{blue}{(lower expected cost or faster processing)} but may involve waiting. This trade-off, together with the requirement that the server must complete the selected service, distinguishes our work from classical routing control. When the decision is made upon initial job interaction, the problem reduces to the control of a single-stage queue, as studied in \cite{lu2026balancingindependentcollaborativeservice}, where an optimal policy is shown to be monotone.
From a queueing perspective, we assume that once a server begins processing a job in the first stage, it remains attached to that job through both service stages before handling a new one. This eliminates the flexibility of server assignment across stages and contrasts with the classic Jackson network, where each server remains fixed to a particular queue \cite{massey1984operator, massey1986family, massey1987calculating, pender2017approximating}.} 

\section{Model description and main results}\label{sec:description}
We formulate the problem using a stochastic clearing system model, where there are no external arrivals, but a fixed number of jobs in the system initially. All jobs begin service on a first-come-first-served basis, and the goal is to find a control policy that minimizes the cost of emptying the system (formalized below).

\textcolor{blue}{Jobs that begin in the first stage await initial service by one of two servers. Upon completion of the first stage, a decision determines whether the server-job pair proceeds to a parallel-service station or a single-service station for second-stage processing. For each job, this is the only time a decision is made throughout the process.
The server is committed to the job until both stages are completed and cannot serve another job during this period. 
We refer to the initial service station, the parallel-service station, and the single-service station as Stations 0, 1, and 2, respectively.
The service times at Stations 0, 1, and 2} are assumed to be exponentially distributed with rates \textcolor{blue}{$\mu_0$, $\mu_1$, and $\mu_2$, respectively.} 

Next, we define a stochastic process describing the evolution of the system and the decision-making scenario. Define the state space 
\begin{align*}
    \X & = \Big\{(i,j,k,\ell)\in \Z_+^4\Big{|} i = 0, j+k+\ell = 1 \text{ or } i \geq 0, j+k+\ell = 2\Big\}, 
\end{align*}
where
\begin{itemize}
    \item $i$ is the number of jobs in the Station 0 queue (not including those in service),
    \item $j$ is the number in service at Station 0,
    \item $k$ is the number of jobs currently receiving Station 1 service, and
    \item $\ell$ is the number of jobs at Station 2 (including those
        in service).
\end{itemize}
We make the following assumption, which is important to define the state and action spaces that follow.
\begin{assumption}\label{assm:non-preemptive}
    All services are non-preemptive. 
\end{assumption}
When a server completes a Station 0 service, the control action is to decide what kind of service should be provided. Since there are only two servers, and the last three components of the state space capture where service is currently occurring, the set of states where a decision is made is 
\textcolor{blue}{
\begin{align*}
    \X_{D} := \Big\{x \in \X | x\in  \{(0,1,0,0)\} \cup \{(i,2,0,0), (i,1,0,1), (i,1,1,0)\text{ for }i \geq 0\}\Big\}. 
\end{align*} 
}
In that sense, the action set is
\begin{align*}
    \A(i,j,k,\ell) & = \begin{cases}
        \{1,2\} & \text{if $x \in \X_D$,}\\
        \{0\} & \text{otherwise.}
    \end{cases}
\end{align*}
When $x \in \X_D$, each allowable action $a \in \A(i,j,k,\ell)$ is a binary variable, where $a = 1$ \textcolor{blue}{routes the server-job pair to Station 1 and $a = 2$ sends the pair to Station 2}. \textcolor{blue}{When $x \notin \X_D$, the action set reduces to the singleton $\{0\}$, where $a = 0$ is a placeholder (dummy) action indicating that no decision is taken.}
Note that if the server decides to complete the \textcolor{black}{second-stage} service in the parallel-service \textcolor{black}{station}, there is no queueing and the service begins immediately.

Let $\{Q_{0,0}(t), t\geq 0\}$ and $\{Q_{0,1}(t), t\geq 0\}$ denote the stochastic processes describing the number of jobs in queue and in service at Station 0, respectively. Let $\{Q_1(t), t\geq 0\}$ and $\{Q_2(t), t\geq 0\}$ denote the number of jobs at Stations 1 and 2, respectively. \textcolor{blue}{Note that at any time $t$, $\big(Q_{0,0}(t), Q_{0,1}(t), Q_1(t), Q_2(t)\big) \in \X$.}
As a measure of service preference, we assume the system is charged $h_0$ per job, per unit of time either waiting or in service at Station 0. Similarly, the system is charged $h_1$ and $h_2$ per job, per unit time at Stations 1 and 2, respectively. As is common in the queueing literature, throughout the paper we term these costs \textit{holding costs}. 
% Unless otherwise stated, Assumption \ref{assm:cost-to-serve} below holds in our main results.
% \begin{assumption}\label{assm:cost-to-serve}
%     \textcolor{blue}{$\frac{h_1}{\mu_1} \geq \frac{h_2}{\mu_2}$}. 
% \end{assumption}
% This indicates that ignoring any other considerations, Station 2 is the preferred choice for a single service. 
For any control policy $\pi$, the total cost $T^{\pi}$ given initial state $s=(i,j,k,\ell) \in \X$ is
\begin{align*}
  T^{\pi}(s) &  = \int_{0}^{\infty}\Big(h_0\big(Q^{\pi}_{0,0}(t) + Q^{\pi}_{0,1}(t)\big) +
    h_1 Q^{\pi}_1(t) + h_2 Q^{\pi}_2(t)\Big) dt,
\end{align*}
where $\big(Q^{\pi}_{0,0}(0), Q^{\pi}_{0,1}(0), Q^{\pi}_1(0), Q^{\pi}_2(0)\big) =(i,j,k,\ell) = s$, and the dependence on the control policy has been added to each stochastic process.
Note that if service completion generates revenue, it remains constant across all policies in a clearing system. Therefore, it suffices to consider only the total cost incurred.

Starting with any finite number of jobs initially, under any non-idling policy $T^{\pi} < \infty$ almost surely. The goal is to find a policy that minimizes the expectation of $T^{\pi}$ for any (fixed) initial state $s=(i,j,k,\ell) \in \X$, i.e.,
\begin{align*}
  \mathbb{E}\big[T^{\pi}(s)\big] &  = \mathbb{E}_{s} \Bigg[\int_{0}^{\infty}\Big(h_0\big(Q^{\pi}_{0,0}(t) + Q^{\pi}_{0,1}(t)\big) +
    h_1 Q^{\pi}_1(t) + h_2 Q^{\pi}_2(t)\Big) dt \Bigg].
\end{align*}

Define $v(i,j,k,\ell)$, often referred to as the \emph{value function}, to be the optimal expected total cost incurred starting at state $(i,j,k,\ell)\in\X$ until the system clears all the jobs (so that $v(0,0,0,0) = 0$). Define $d(j,k,\ell) = j\mu_0 + k\mu_1 + \min\{\ell,1\}\mu_2$. By the \emph{Principle of Optimality} (see
Section 4.3 in \cite{puterman2014markov}), any non-zero state $v(i,j,k,\ell)$ satisfies the optimality equations \textcolor{blue}{(also known as the \textit{Bellman equations})} as follows (c.f. Theorem 7.3.3 in \cite{puterman2014markov} for the optimality equations of discrete-time MDP).
\begin{enumerate}
    \item If $i = 0$,
    \begin{align}\label{eq:opt-zero}
        \begin{split}
            v(0,j,k,\ell) & = \frac{j h_0+ k h_1 + \ell h_2 }{d(j,k,\ell)}
                + \Big[ \frac{k \mu_1}{d(j,k,\ell)} v(0, j, k-1, \ell)  \\
                & \quad + \frac{\min\{\ell, 1\}\mu_2}{d(j,k,\ell)} v(0, j, k, \ell-1) \\
                & \quad + \frac{j \mu_0}{d(j,k,\ell)} \min\{v(0, j-1, k+1, \ell), v(0, j-1, k, \ell+1)\} \Big].
        \end{split}
    \end{align}
    \item If $i \geq 1$,
    \begin{align}\label{eq:opt1}
        \begin{split}
            v(i,j,k,\ell) & = \frac{(i+j) h_0+ k h_1 + \ell h_2 }{d(j,k,\ell)}
                + \Big[ \frac{k \mu_1}{d(j,k,\ell)} v(i-1, j+1, k-1, \ell) \\
                & \quad + \frac{\min\{\ell, 1\}\mu_2}{d(j,k,\ell)} v(i-1, j+1, k, \ell-1) \\
                & \quad + \frac{j \mu_0}{d(j,k,\ell)} \min\{v(i, j-1, k+1, \ell), v(i, j-1, k, \ell+1)\} \Big].
        \end{split}
    \end{align}
\end{enumerate}

% As is simple to see, the decision of whether to work at Station 1 or 2 hinges on the difference in the total cost starting in different states (as captured by the value functions).
{\color{blue}
These equations admit a natural interpretation. Starting from state $s=(i,j,k,\ell) \in \X$, the next service completion occurs after an expected time of $\frac{1}{d(j,k,\ell)}$. Over this interval, the system incurs the expected holding cost equal to $\frac{(i+j) h_0+ k h_1 + \ell h_2}{d(j,k,\ell)}$, after which it moves to a successor state $s'$. 
When the next completion is an initial screening, which occurs with probability $\frac{j \mu_0}{d(j,k,\ell)}$, the optimal decision-maker chooses the downstream transition that minimizes the cost-to-go. In particular, if $v(i, j-1, k+1, \ell) < v(i, j-1, k, \ell+1)$, the server-job pair is sent to Station 1 and $s' = (i, j-1, k+1, \ell)$; otherwise, it is sent to Station 2, yielding $s' = (i, j-1, k, \ell+1)$. 
If some other service completion occurs first, the state changes accordingly, and no decision is made. The Bellman equations therefore express the value $v(s)$ as the expected one-step cost plus the optimal expected continuation cost after the next event.}

\textcolor{black}{From the optimality equations \eqref{eq:opt-zero}--\eqref{eq:opt1}, it follows that, upon completion of first-stage service in state $(i,j,k,\ell) \in \X_D$, the routing decision is determined by the sign of the following difference:}
\begin{align}
    D(i,j,k,\ell) := v(i, j-1, k+1, \ell) - v(i, j-1, k, \ell+1). \label{def:D}
\end{align}
\textcolor{black}{A negative (resp. positive) value favors sending the job to Station 1 (resp. Station 2).}

\subsection{Main results}\label{sec:main-results}
Consider any decision state $(i,j,k,\ell) \in \X_D$. Immediately after a service completion at Station 0, \textcolor{black}{the decision-maker compares the value functions corresponding to sending the server-job pair to} Station 1 or Station 2, namely $v(i,j-1,k+1,\ell)$ and $v(i,j-1,k,\ell+1)$, respectively. 
{\color{blue}For example, an optimal policy chooses Station 1 (Station 2) in state $(0,1,0,0)$ if $D(0,1,0,0) = v(0,0,1,0) - v(0,0,0,1) \leq (\geq) 0$. A straightforward calculation shows that
\begin{align}
    D(0,1,0,0)  & = \textcolor{blue}{\frac{h_1}{\mu_1} - \frac{h_2}{\mu_2}}. \label{eq:opt_single}
\end{align}
Consequently, the sign of $\frac{h_1}{\mu_1} - \frac{h_2}{\mu_2}$ fully determines an optimal decision at $(0,1,0,0)$. 
In particular, if $\frac{h_1}{\mu_1} > \frac{h_2}{\mu_2}$, Station 2 is preferred from the perspective of minimizing the expected cost of serving a single job at the second stage. 
More generally, for $n = 0,1,2$, the ratio $\frac{h_n}{\mu_n}$ represents the expected holding cost incurred during service at Station $n$. It therefore serves as a natural measure of in-service cost efficiency: a smaller ratio indicates lower expected cost while a job is being served, and hence greater cost effectiveness. In this sense, although Station 1 offers the advantage of immediate parallel service, Station 2 may still be attractive when its service is sufficiently fast and/or sufficiently cost-efficient.}

It remains to discuss the general decision states $(i,1,1,0)$, $(i,1,0,1)$ and $(i,2,0,0)$ in $\X_D$. For ease of exposition, let $\X_{\widehat{D}} := \X_D \setminus (0,1,0,0)$. In these decision states, both servers are busy, so the total number of jobs following (and including) the initial service is exactly two, i.e., $j+k+\ell = 2$. 
\textcolor{blue}{We begin by presenting our main results in Theorem \ref{thm:large_queue}, which characterizes the structure of optimal policies in the heavy-initial-queue regime; see also Tables \ref{table:policy_empty_Station2}--\ref{table:policy_queue_in_Station2}.} 
% We present our main results in Theorems \ref{thm:large_queue}--\ref{thm:threshold}; \textcolor{blue}{see also Tables \ref{table:policy_empty_Station2}--\ref{table:policy_queue_in_Station2} for results related to Theorem \ref{thm:large_queue}}. The proofs are provided in Appendix \ref{sec:appendix}.

\begin{theorem} \label{thm:large_queue}
    The following results hold for decision states $(i,j,k,\ell) \in \X_{\widehat{D}}$:
    \begin{enumerate} 
        \item Consider $\ell = 0$, that is, decision states of the form $(i,1,1,0)$ or $(i,2,0,0)$.
        \begin{enumerate} [ref=\theenumi(\alph*)]
            \item \label{state:empty_Station2_slower} 
            If \textcolor{black}{$\mu_1 > \mu_2$}, there exists a finite $i'_{j,k,\ell}$ such that an optimal policy chooses Station 1 for all $i \geq i'_{j,k,\ell}$. If, in addition, $\frac{h_1}{\mu_1} \leq \frac{h_2}{\mu_2}$ holds, then there exists an optimal policy that chooses Station 1 for all $i$.
            \item \label{state:empty_Station2_equal} 
            If $\mu_1 = \mu_2$, \textcolor{blue}{there exists an optimal policy that chooses Station 1 for all $i$ when $\frac{h_1}{\mu_1} \leq \frac{h_2}{\mu_2}$}, and chooses Station 2 for all $i$ when $\frac{h_1}{\mu_1} \geq \frac{h_2}{\mu_2}$.
            \item \label{state:empty_Station2_faster} 
            If $\mu_1 < \mu_2$, \textcolor{blue}{there exists a finite $i'_{j,k,\ell}$ such that an optimal policy chooses Station 2 for all $i \geq i'_{j,k,\ell}$.} If, in addition, $\frac{h_1}{\mu_1} \geq \frac{h_2}{\mu_2}$ holds, there exists an optimal policy that chooses Station 2 for all $i$.
        \end{enumerate}
        \item Consider $\ell = 1$, that is, decision states of the form $(i,1,0,1)$.
        \begin{enumerate} [ref=\theenumi(\alph*)]
            \item \label{state:queue_in_Station2_slower} 
            If \textcolor{black}{$\mu_2 \leq 2\mu_1$}, there exists a finite $i'_{1,0,1}$ such that an optimal policy chooses Station 1 for all $i \geq i'_{1,0,1}$. 
            \textcolor{blue}{If, in addition, $\frac{h_1}{\mu_1} \leq \frac{2h_2}{\mu_2}$ holds, there exists an optimal policy that chooses Station 1 for all $i$.} 
            \item \label{state:queue_in_Station2_faster}
            Suppose $\mu_2 > 2\mu_1$.
            \begin{enumerate} [ref=\theenumii\roman*]
                \item  \label{state:queue_in_Station2_faster_slow_initial} 
                There exists a sufficiently small $\mu_0'$ such that, for any fixed $\mu_0 < \mu_0'$, there is a finite $i'_{1,0,1}$ for which an optimal policy chooses Station 2 for all $i \geq i'_{1,0,1}$.
                \item \label{state:queue_in_Station2_faster_fast_initial} 
                \textcolor{blue}{There exists a sufficiently large $\mu_0'$ such that, for any fixed $\mu_0 > \mu_0'$, there is a finite $i'_{1,0,1}$ for which an optimal policy chooses Station 1 for all $i \geq i'_{1,0,1}$.}
            \end{enumerate}
        \end{enumerate}
    \end{enumerate}
\end{theorem}

\begin{table}[htbp]
\captionsetup{labelfont={color=blue}, textfont={color=blue}}
\caption{Optimal policies for decision states of the form $(i,1,1,0)$ and $(i,2,0,0)$. All statements refer to Theorem \ref{thm:large_queue}.}
\label{table:policy_empty_Station2}
\centering
{\color{blue}
\renewcommand{\arraystretch}{1.8}
\begin{tabular}{|c|c|c|}
\hline
& $\dfrac{h_1}{\mu_1}\leq \dfrac{h_2}{\mu_2}$ & $\dfrac{h_1}{\mu_1}\geq \dfrac{h_2}{\mu_2}$ \\[0.3em]
\hline
$\mu_1>\mu_2$ & $a=1$ (Statement \ref{state:empty_Station2_slower})  & $a=1$ for large $i$ (Statement \ref{state:empty_Station2_slower}) \\ [0.3em]
\hline
$\mu_1=\mu_2$ & $a=1$ (Statement \ref{state:empty_Station2_equal}) & $a=2$ (Statement \ref{state:empty_Station2_equal}) \\ [0.3em]
\hline
$\mu_1<\mu_2$ & $a=2$ for large $i$ (Statement \ref{state:empty_Station2_faster}) & $a=2$ (Statement \ref{state:empty_Station2_faster}) \\[0.3em]
\hline
\end{tabular}
}
\end{table}

\begin{table}[htbp]
\captionsetup{labelfont={color=blue}, textfont={color=blue}}
\caption{Optimal policies for decision states of the form $(i,1,0,1)$. All statements refer to Theorem \ref{thm:large_queue}.}
\label{table:policy_queue_in_Station2}
\centering
{\color{blue}
\renewcommand{\arraystretch}{1.8}
\begin{tabular}{|c|c|c|}
\hline
 & $\dfrac{h_1}{\mu_1}\leq \dfrac{2h_2}{\mu_2}$ & $\dfrac{h_1}{\mu_1}> \dfrac{2h_2}{\mu_2}$ \\ [0.3em]
\hline
$\mu_2 \leq 2\mu_1$ & $a=1$ (Statement \ref{state:queue_in_Station2_slower}) & $a=1$ for large $i$ (Statement \ref{state:queue_in_Station2_slower}) \\[0.3em]
\hline
$\mu_2 > 2\mu_1$ & \multicolumn{2}{c|}{\begin{tabular}{c}
For any $\mu_0$ small enough, $a=2$ for large $i$ (Statement \ref{state:queue_in_Station2_faster_slow_initial}) \\[-1em]
For any $\mu_0$ large enough, $a=1$ for large $i$ (Statement \ref{state:queue_in_Station2_faster_fast_initial})
\end{tabular}
}  \\
\hline
\end{tabular}}
\end{table}

Next, we provide intuitive explanations to highlight the implications of these results. 
{\color{blue}
Consider the decision states $(i,j,k,\ell) \in \X_{\widehat{D}}$ with $\ell = 0$. That is, states of the form $(i,1,1,0)$ and $(i,2,0,0)$. In this case, sending the current server-job pair to Station 2 results in immediate service. Consequently, choosing Station 1 (resp. Station 2) yields expected immediate service time and cost equal to $\frac{1}{\mu_1}$ and $\frac{h_1}{\mu_1}$ (resp. $\frac{1}{\mu_2}$ and $\frac{h_2}{\mu_2}$), respectively. The analysis therefore reduces to the sign of the following time and cost differences:}
\begin{align}
    &\textcolor{blue}{\frac{1}{\mu_1} - \frac{1}{\mu_2}} \quad & \textcolor{blue}{\frac{h_1}{\mu_1} - \frac{h_2}{\mu_2} \label{eq:empty-Station2-time-cost}}
\end{align}
\begin{itemize}
    \item {\color{black} When $\mu_1 > \mu_2$, i.e., $\frac{1}{\mu_1} < \frac{1}{\mu_2}$, the server completes its current job faster in Station 1 than in Station 2 (See Statement \ref{state:empty_Station2_slower} of Theorem \ref{thm:large_queue}).} This time savings allows the server to return to upstream service more quickly, reducing the cost incurred at \textcolor{black}{the first stage} due to waiting jobs at Station 0. When the initial queue is large, this upstream cost reduction outweighs the (potentially) higher service cost of this job at Station 1 compared to Station 2, making Station 1 the preferred option. 
    \begin{itemize}
        \item If, in addition, \textcolor{black}{$\frac{h_1}{\mu_1} \leq \frac{h_2}{\mu_2}$}, then Station 1 is not only faster but also incurs lower or equal cost, further reinforcing its preference. In this case, there exists an optimal policy that always chooses Station 1, regardless of the size of the initial queue. 
    \end{itemize}
    {\color{blue}
    \item By contrast, when \textcolor{blue}{$\mu_1 < \mu_2$}, the result is symmetric since Station 2 is free of queueing effects in this case, as captured in Statement \ref{state:empty_Station2_faster} of Theorem \ref{thm:large_queue}.

   \item When $\mu_1 = \mu_2$, the two options offer identical service speeds, so the optimal policy is determined by the relative service costs: choose Station 1 if \(\frac{h_1}{\mu_1} \leq \frac{h_2}{\mu_2}\), and Station 2 otherwise. This behavior is characterized in Statement \ref{state:empty_Station2_equal} of Theorem \ref{thm:large_queue}.
   }
\end{itemize}
{\color{blue}
The optimal policy becomes more subtle at decision states $(i,1,0,1)$, where sending the server-job pair to Station 2 results in downstream waiting before starting service.}
Given that the previous job is being processed at Station 2, the expected time (cost) to complete the current job at Station 2 equals \textcolor{black}{$\frac{2}{\mu_2}$ ($\frac{2h_2}{\mu_2}$)}. 
{\color{blue}
As in the previous case, we first examine the signs of the corresponding differences in immediate expected time and cost:
\begin{align}
    &\frac{1}{\mu_1} - \frac{2}{\mu_2} \quad & \frac{h_1}{\mu_1} - \frac{2h_2}{\mu_2} \label{eq:queue-in-Station2-time-cost}
\end{align}

The downstream queueing complicates the policy structure beyond the comparison of these immediate time and cost differences, as it also has a system-wide impact. In particular, from the optimality equations \eqref{eq:opt-zero}--\eqref{eq:opt1}, proceeding to Station 2 yields a smaller total transition rate ($d(0,0,2) = \mu_2$) compared to Station 1 ($d(0,1,1) = \mu_1 + \mu_2$). Operationally, this means that the server is \textit{blocked} with its paired job and cannot initiate any task until service at Station 2 begins. 
We therefore refer to states of the form $(i,1,0,1)$ as ``blocked decision states''. By contrast, states $(i,1,1,0)$ and $(i,2,0,0)$, where routing to Station 2 initiates service immediately, are called ``unblocked decision states''.

Moreover, the requirements that each server remains engaged with its assigned job through both service stages and that decisions occur only upon an upstream service completion introduce a nontrivial interplay between the relative speeds of upstream and downstream processing. 
When $\mu_0$ is small, downstream work typically finishes before the next upstream completion and hence before the next decision epoch. Consequently, Station 2 is often available at subsequent decision points, and blocking effects of routing to Station 2 are largely limited to the current job. In this case, the optimal policy is driven primarily by the immediate time and cost trade-off in \eqref{eq:queue-in-Station2-time-cost}. 
By contrast, when $\mu_0$ is large, upstream service tends to complete before downstream service. Decisions are then often made while Station 2 is still occupied, so blocking propagates to subsequent jobs, making Station 1 more attractive.
}

\noindent We now consider the following cases:
\begin{itemize}
    \item When {\color{blue}$\mu_2 \leq 2\mu_1$, i.e., $\frac{1}{\mu_1} \leq \frac{2}{\mu_2}$}, serving the job at Station 1 allows the server to return to upstream service at Station 0 more quickly, reducing the cost associated with waiting jobs at Station 0. Since Station 1 involves no blocking considerations, this upstream cost reduction dominates the (potentially) higher downstream cost at Station 1 (if \textcolor{black}{$\frac{h_1}{\mu_1} > \frac{2h_2}{\mu_2}$}) when the initial queue is large, leading to the decision to choose Station 1 service. This is formalized in Statement \ref{state:queue_in_Station2_slower} of Theorem \ref{thm:large_queue}. 
    \begin{itemize}
        \item Statement \ref{state:queue_in_Station2_slower} of Theorem \ref{thm:large_queue} further identifies cases where the expected cost at Station 1 is also lower (\textcolor{black}{$\frac{h_1}{\mu_1} \leq \frac{2h_2}{\mu_2}$}), resulting in an optimal policy that always selects Station 1, regardless of the initial queue length.
    \end{itemize}
    \item {\color{blue}
    If $\mu_2 > 2\mu_1$, i.e., $\frac{1}{\mu_1} > \frac{2}{\mu_2}$, proceeding to Station~2 yields a smaller immediate expected completion time for the current job. However, unlike in previous cases, the optimal policy for a large initial queue is not solely determined by the sign of this time difference, but also depends on the interplay between upstream and downstream service rates.
    }
    \begin{itemize}
        \item 
        When $\mu_0$ is small, blocking effects are minimal beyond the current job. As a result, the immediate reduction in upstream cost, driven by the faster return of the server from Station 2, becomes the dominant factor when the initial queue is large, making Station 2 the preferred choice. This result aligns with Statement \ref{state:queue_in_Station2_faster_slow_initial} of Theorem \ref{thm:large_queue}. 
        \item 
        {\color{blue} 
        When $\mu_0$ is large, the accumulation of blocking effects beyond the current job makes Station 1 preferable for large initial queues. This corresponds to Statement \ref{state:queue_in_Station2_faster_fast_initial} of Theorem \ref{thm:large_queue}.}
    \end{itemize}
\end{itemize}

{\color{blue}
Taken together, for decision states $(i,j,k,\ell) \in \X_{\widehat{D}}$, Theorem \ref{thm:large_queue} characterizes optimal policies in the regime of heavy initial queues at Station 0. In particular, it identifies conditions under which the optimal policy stabilizes beyond a threshold in the initial queue length, denoted by $i'_{j,k,\ell}$, as well as the corresponding stabilizing decision. To accommodate different parameter regimes, this threshold may be zero or a finite number, corresponding to policies that are constant for all $i$ or eventually stabilize, respectively. Consequently, for sufficiently large $i$, the optimal policy becomes independent of the exact value of $i$.}

{\color{blue}The proof of Theorem \ref{thm:large_queue} is detailed in Appendix \ref{sec:appendix}; here we provide an outline. 
All arguments rely on characterizing the sign of $D(i,j,k,\ell)$, defined in \eqref{def:D}, for the decision state of interest $(i,j,k,\ell) \in \X_{\widehat{D}}$, leveraging the optimality equations \eqref{eq:opt-zero} and \eqref{eq:opt1}.
For example, Statement \ref{state:queue_in_Station2_slower} is established by showing that, when $\frac{1}{\mu_1} \leq \frac{2}{\mu_2}$, we have $D(i,1,0,1) \leq 0$ for all $i \geq i'_{1,0,1}$ for some finite $i'_{1,0,1}$. In contrast, Statement \ref{state:queue_in_Station2_faster_slow_initial} is established by showing that when $\frac{1}{\mu_1} > \frac{2}{\mu_2}$, there exists $\mu_0$ sufficiently small, such that $D(i,1,0,1) \geq 0$ for all $i \geq i'_{1,0,1}$ for some (possibly different) finite $i'_{1,0,1}$.}

\textcolor{blue}{Theorem \ref{thm:threshold} below establishes a stronger monotone characterization under the assumption $\mu_1 = \mu_2$, leading to a threshold-type optimal policy. }

\begin{theorem} \label{thm:threshold}
    Suppose \textcolor{black}{$\mu_1 = \mu_2$} and consider $(i,j,k,\ell) \in \X_{\widehat{D}}$.
    Recall that the number of jobs following (and including) the initial service is a constant two. There exists an optimal policy that is monotone in the following sense.
    \begin{enumerate} [label=\textcolor{blue}{\arabic*.}, ref=\arabic*]
        \item \label{state:threshold_i}
        For fixed $j,k, \ell$, it is monotone decreasing (e.g., from $a = 2$ to $a=1$) in $i$. That is, there exists a control threshold such that it is optimal to work at Station 2 in $(i,j,k,\ell) \in \X_{\widehat{D}}$ but not in $(i+m,j,k,\ell) \in \X_{\widehat{D}}$ for all $m > 0$.
        \item Consider fixed $i$.
        \begin{enumerate} [label=\textcolor{blue}{(\alph*)}, ref=\theenumi(\alph*)]
            \item \label{state:threshold12}
            There exists a control threshold such that it is optimal to work at Station 2 in $(i, j, k, \ell) \in \X_{\widehat{D}}$ and not in $(i, j, k-m, \ell+m) \in \X_{\widehat{D}}$ for all $m > 0$.
            \item \label{state:threshold02}
            There exists a control threshold such that it is optimal to work at Station 2 in $(i, j, k, \ell) \in \X_{\widehat{D}}$ and not in $(i, j-m, k, \ell+m) \in \X_{\widehat{D}}$ for all $m > 0$.
            \item \label{state:threshold01}
            There exists a control threshold such that it is optimal to work at Station 2 in $(i, j, k, \ell) \in \X_{\widehat{D}}$ and not in $(i, j+m, k-m, \ell) \in \X_{\widehat{D}}$ for all $m > 0$.
        \end{enumerate}
    \end{enumerate}
\end{theorem}

{\color{blue}
% Building on the insights from Theorem \ref{thm:large_queue}, which focuses on large initial queues, we now turn to Theorem \ref{thm:threshold}. This result establishes a stronger structural characterization under the assumption $\mu_1 = \mu_2$, leading to a threshold-type optimal policy. 
% We now turn to Theorem \ref{thm:threshold}, which establishes a stronger monotone characterization under the assumption $\mu_1 = \mu_2$, leading to a threshold-type optimal policy. 
Theorem \ref{thm:threshold} shows monotonicity with respect to both the initial queue length and the distribution of jobs across stations.
Statement \ref{state:threshold_i} indicates that, as the number of initially queued jobs increases, it becomes more desirable to avoid blocking and obtain immediate service, making the single-service station less attractive.
Statement \ref{state:threshold12} (respectively, Statement \ref{state:threshold02}) shows that when more jobs are concentrated at Station 2 and fewer at Station 1 (respectively, Station 0), there exists an optimal policy that tends to favor proceeding to Station 1.
Finally, Statement \ref{state:threshold01} indicates that when more jobs are present at Station 0 and fewer at Station 1, an optimal policy again prefers routing to Station 1.}

\textcolor{blue}{
The proof of Theorem \ref{thm:threshold}, which is based on the monotonicity properties of $D$ with respect to the state variables, is also provided in Appendix \ref{sec:appendix}. Specifically, for all $(i,j,k,\ell) \in \X_{\widehat{D}}$, we show:}
\begin{align}
    & \textcolor{blue}{D(i,j,k,\ell) - D(i+1,j,k,\ell) \geq 0,} \label{eq:D_mono_i} \\
    & \textcolor{blue}{D(i,j,k,\ell) - D(i,j,k-1,\ell+1) \geq 0,  \quad \text{for }(i,j,k-1,\ell+1) \in \X_{\widehat{D}},} \nonumber \\
    & \textcolor{blue}{D(i,j,k,\ell) - D(i,j-1,k,\ell+1) \geq 0, \quad \text{for }(i,j-1,k,\ell+1) \in \X_{\widehat{D}},} \nonumber \\
    & \textcolor{blue}{D(i,j,k,\ell) - D(i,j+1,k-1,\ell) \geq 0, \quad \text{for }(i,j+1,k-1,\ell) \in \X_{\widehat{D}}. \nonumber}
\end{align}

{\color{blue}
We conclude this section with a discussion of technical challenges. 
A natural question is whether the structural insights for large initial queues in Theorem \ref{thm:large_queue} can be strengthened to a monotone or threshold-type structure, whereby the optimal policy selects one action if $i <N_{j,k,\ell}$ and the other otherwise, for some threshold $N_{j,k,\ell}$ (see, e.g., Statement \ref{state:threshold_i} in Theorem \ref{thm:threshold}). However, such a characterization does not hold in general.
In particular, when $\mu_1 < \mu_2$, the interaction between upstream and downstream service rates may cause the optimal policy at states $(i,1,0,1)$ to switch more than once as $i$ increases under certain parameter regimes, leading to a non-monotonic structure. 
This non-threshold behavior is further illustrated through numerical examples (see \ref{example:non-monotone_slow_station2} and \ref{example:non-monotone_fast_station2}) presented in Section \ref{sec:num}.

Furthermore, even in parameter regimes where an optimal policy appears monotone, establishing such a structure analytically is challenging. For example, consider decision states $(i,j,k,\ell) \in \X_{\widehat{D}}$ such that $\ell = 0$ under $\mu_1 > \mu_2$ with $\frac{h_1}{\mu_1} > \frac{h_2}{\mu_2}$ (note that if $\frac{h_1}{\mu_1} \leq \frac{h_2}{\mu_2}$, then $a=1$ for all $i$ by Statement \ref{state:empty_Station2_slower} in Theorem \ref{thm:large_queue}). Standard approaches rely on showing that $D(i,j,k,\ell)$ is monotone non-increasing in $i$ (see \eqref{eq:D_mono_i}), which would imply a threshold structure with the threshold defined as $N_{j,k,\ell} := \min\{i \in \Z_+|D(i,j,k,\ell) \leq 0\}$. However, this property fails in general. For instance, under $\mu_0 = 5, \mu_1 = 3.1, \mu_2 = 3, h_0 = 0.1, h_1 = 22, h_2 = 10$, we have $D(1,2,0,0) <D(2,2,0,0)$. This non-monotonicity prevents the use of standard arguments and obstructs a general proof of threshold structure.
}

{\color{blue}
The insights of our main theorems extend naturally to controlled queues with multiple servers. In this setting, the time and cost differences in \eqref{eq:queue-in-Station2-time-cost} generalize to
\begin{align*}
    &\frac{1}{\mu_1} - \frac{\ell+1}{\mu_2}, \quad & \frac{h_1}{\mu_1} - \frac{(\ell+1)h_2}{\mu_2},
\end{align*}reflecting increased downstream congestion.
Despite this structural similarity, extending the analysis to multiple servers is technically challenging, and key arguments in the two-server case no longer carry over directly. In particular, the right-hand side of inequality \eqref{eq:decision_bd2} in Statement \ref{state4:general_rates} of Lemma \ref{lemma:general_rates} does not, in general, upper-bound the difference $D(i,j,k,\ell)$ on the left-hand side. For example, under $\mu_0 = 0.1, \mu_1 = 10, \mu_2 = 15, h_0 = 0.1, h_1 = 1, h_2 = 1.12$ with three servers, \eqref{eq:decision_bd2} is violated at state $(2,1,0,2)$: the left-hand side equals $D(2,1,0,2) = -0.1268$, while the right-hand side equals $-0.1273$. This breakdown highlights a fundamental obstacle to extending our proof techniques beyond the two-server setting, even though the qualitative insights are expected to remain valid.
}
\section{Numerical study}\label{sec:num}
{\color{black}
In this section, we further illustrate the structural properties of optimal policies through a numerical study. In particular, we provide practical insights by comparing several natural heuristic policies against the optimal policy.}
Since the optimal policy remains unchanged when all holding costs are scaled by a positive constant (with service rates fixed), we set $h_1=1$ without loss of generality. Similarly, we fix $\mu_1 = 10$.
{\color{blue}
To align with the partition of the parameter space in Theorem \ref{thm:large_queue} and Tables \ref{table:policy_empty_Station2}--\ref{table:policy_queue_in_Station2}, we select three representative pairs $(\mu_2,h_2)$ for each regime; see Table \ref{table:mu2h2}. For each pair, we test all configurations of $\mu_0$ and $h_0$ listed in Table \ref{table:params}.
}

% In this section, we compare the values, i.e., the total expected costs of some natural heuristic policies with that of the optimal policy. In particular, the values beginning in states $(20,j,k,\ell) \in \X_{\widehat{D}}$ are calculated under these policies with various parameter configurations. The parameters are listed in Table \ref{table:params}. 
% % All combinations are tested under the natural constraint that $\frac{h_1}{\mu_1} > \frac{h_2}{\mu_2}$, \textcolor{blue}{which implies that the single-service station is individually preferable from the perspective of a single job}.
% Note that since the optimal policy remains the same when all the holding costs are multiplied by a positive constant (while service rates remain unchanged), without loss of generality, we fix $h_1=1$. Analogously, we let $\mu_1 = 10$. 

\begin{table}[bpth]
    \captionsetup{labelfont={color=blue}, textfont={color=blue}}
    \caption{Representative $(\mu_2,h_2)$ pairs used in the numerical study, organized by downstream rate and cost regimes.} \label{table:mu2h2}
    \centering
    {\color{blue}
    \begin{tabular}{|c|c|c|c|}
        \hline
        \rule{0pt}{4.2ex}
        & \(\dfrac{h_1}{\mu_1} \le \dfrac{h_2}{\mu_2}\)
        & \(\dfrac{h_2}{\mu_2} < \dfrac{h_1}{\mu_1} \le \dfrac{2h_2}{\mu_2}\)
        & \(\dfrac{h_1}{\mu_1} > \dfrac{2h_2}{\mu_2}\) \\ [0.8em]
        \hline
        \(\mu_1 \ge \mu_2\)
        &
        \begin{tabular}{c}
        \((8,1.0)\)\\[-1mm]
        \((8,1.5)\)\\[-1mm]
        \((10,1.0)\)
        \end{tabular}
        &
        \begin{tabular}{c}
        \((8,0.5)\)\\[-1mm]
        \((8,0.7)\)\\[-1mm]
        \((10,0.5)\)
        \end{tabular}
        &
        \begin{tabular}{c}
        \((8,0.3)\)\\[-1mm]
        \((8,0.2)\)\\[-1mm]
        \((10,0.3)\)
        \end{tabular}
        \\
        \hline
        \(\mu_1 < \mu_2 \le 2\mu_1\)
        &
        \begin{tabular}{c}
        \((15,2.0)\)\\[-1mm]
        \((15,1.5)\)\\[-1mm]
        \((18,2.0)\)
        \end{tabular}
        &
        \begin{tabular}{c}
        \((15,1.0)\)\\[-1mm]
        \((15,1.4)\)\\[-1mm]
        \((18,1.0)\)
        \end{tabular}
        &
        \begin{tabular}{c}
        \((15,0.5)\)\\[-1mm]
        \((15,0.7)\)\\[-1mm]
        \((12,0.5)\)
        \end{tabular}
        \\
        \hline
        \(\mu_2 > 2\mu_1\)
        &
        \begin{tabular}{c}
        \((25,3.0)\)\\[-1mm]
        \((25,2.6)\)\\[-1mm]
        \((30,3.0)\)
        \end{tabular}
        &
        \begin{tabular}{c}
        \((25,2.0)\)\\[-1mm]
        \((25,1.3)\)\\[-1mm]
        \((30,2.0)\)
        \end{tabular}
        &
        \begin{tabular}{c}
        \((25,1.0)\)\\[-1mm]
        \((25,0.2)\)\\[-1mm]
        \((30,1.0)\)
        \end{tabular}
        \\
        \hline
    \end{tabular}
}
    \label{tab:placeholder}
\end{table}

\begin{table}[htbp]
\caption{Parameters used in the numerical study.}\label{table:params}
\centering
\begin{tabular}{cc}
\hline
\textbf{Parameters} & \textbf{Values}         \\ \hline
\textbf{$\mu_0$}        & 2, 5, 10, \textcolor{blue}{15,} 20, 30, 40     \\ 
\textbf{$h_0$}         & 0.01, 0.05, 0.1, 0.5, 1 \\ 
\textbf{$\mu_1$}        & 10                \\ 
\textbf{$h_1$}         & 1                       \\ 
% \textbf{$\mu_2$}        & 6, 8, \textcolor{blue}{10,} 12, 15, 25             \\ 
% \textbf{$h_2$}         & 0.2, 0.5, 1, 1.5, 2, \textcolor{blue}{3}     \\ 
\hline
\end{tabular}
\end{table}

\noindent In comparison, the following heuristic policies are evaluated. For $(i,j,k,\ell) \in \X_{\widehat{D}}$,
\begin{enumerate}
[leftmargin=2.5cm, label=Policy \arabic*.]
    \item Always go to \textbf{Station 1}: $a \equiv 1$.
    \item Go to \textbf{Station 1} if the number of jobs waiting at Station 0 is at least the threshold $10$: $a=\mathbb{I}\{i < 10\} + 1$.

    \item Go to \textbf{Station 1} if the number of jobs waiting at Station 0 is at least the threshold $15$: $a=\mathbb{I}\{i < 15\} + 1$.
    
    \item Always go to \textbf{Station 2}: $a \equiv 2$.
    \item Go to \textbf{Station 2} if the single server is available at Station 2: $a=\mathbb{I}\{\ell = 0\} + 1$.
    \item \textcolor{blue}{Go to \textbf{Station 2} if the number of jobs waiting at Station 0 is at least the threshold $10$: $a=\mathbb{I}\{i \geq 10\} + 1$.
    }
    \textcolor{blue}{
    \item Go to \textbf{Station 2} if the number of jobs waiting at Station 0 is at least the threshold $15$: $a=\mathbb{I}\{i \geq 15\} + 1$.
    }
\end{enumerate}
{\color{blue}
For ease of exposition, we say that a policy is of \textbf{Station-1 (resp. Station-2) threshold type with threshold $N$, denoted as $(1, N)$ (resp. $(2,N)$)}, if it routes the server-job pair to Station 1 (resp. Station 2) when the initial queue length is below $N$, and to the other station otherwise. 
Thus, the label indicates the station chosen for relatively small initial queues. 
Here, the dependence of $N$ on $(j,k,\ell)$ is suppressed when clear from context.  
Under this terminology, Policies 1--4 are of \textbf{Station-2 threshold type} with thresholds $N_1 = 0$ , $N_2 = 10$, $N_3 = 15$, and $N_4 = \infty$, respectively, independent of the state. Policy 5 is also of this type, but with a threshold that depends on the state variable $\ell$: 
\begin{align}
    N_5 = \begin{cases}
        \infty & \text{if $\ell = 0$,}\\
        0 & \text{if $\ell = 1$.}
    \end{cases}\label{eq:N5}
\end{align} 
By contrast, Policies 6 and 7 are of \textbf{Station-1 threshold type} with thresholds $N_6 = 10$ and $N_7 = 15$, respectively.
Note that any policy that always selects the same action (for fixed $j,k,\ell$) can be represented as either threshold type: as the type corresponding to that action with threshold $\infty$, or as the opposite type with threshold $0$.
For instance, Policy 1 may equivalently be viewed as a Station-1 threshold-type policy with threshold $\infty$. The same observation applies to Policies 4 and 5.
}

{\color{blue}
We compare the value functions, namely, the total expected holding costs incurred until all jobs are cleared, starting from initial states $(20,j_0,k_0,\ell_0) \in \X_{\widehat{D}}$ under these heuristic policies against an optimal.} 
Performance is evaluated using the relative error, defined as follows:
\begin{align*}
    \textrm{Relative Error} = \frac{v^{\pi}(20,j_0,k_0,\ell_0) - v^*(20,j_0,k_0,\ell_0)}{v^*(20,j_0,k_0,\ell_0)}, 
\end{align*}where $v^{\pi}$ is the value function under each heuristic policy $\pi$ and $v^*$ denotes the optimal value function.
{\color{blue}
The numerical results are summarized in Tables \ref{table:slow-Station2}, \ref{table:mid-Station2}, and \ref{table:fast-Station2}, corresponding to the three service-rate regimes $\mu_1\geq\mu_2$, $\mu_1 < \mu_2 \leq 2\mu_1$, and $\mu_2 > 2\mu_1$, repsectively, as characterized in Theorem \ref{thm:large_queue}. Within each table, results are further organized by the downstream cost regimes. For each heuristic policy and each regime cell, we first compute the relative error defined above for every tested parameter instance in that cell and for each initial state $(20,j_0,k_0,\ell_0) \in \X_{\widehat{D}}$, namely the decision states $(20,1,1,0)$, $(20,1,0,1)$, and $(20,2,0,0)$. The reported maximum, mean, and standard deviation for each policy are then taken over this collection of relative errors. 
Underlined entries correspond to heuristics which align exactly with an optimal in the corresponding regime, while bold entries indicate the heuristic with the smallest mean relative error in that regime.
}

\begin{table}[htbp]
\centering
\renewcommand{\arraystretch}{0.9}
\setlength{\tabcolsep}{5pt}
\small
\captionsetup{labelfont={color=blue}, textfont={color=blue}}
\caption{Relative errors ($\%$) of heuristic policies for decision states $(20,j_0,k_0,\ell_0) \in \X_{\widehat{D}}$ when $\mu_1 \geq \mu_2$, where $N_5$ is given in \eqref{eq:N5}.} \label{table:slow-Station2}
{\color{blue}
\begin{tabular}{|c|c|r|r|r|r|r|r|r|}
\hline
Cost regime & Error &
\begin{tabular}{@{}c@{}}Policy 1\\ $(2,0)$\end{tabular} &
\begin{tabular}{@{}c@{}}Policy 2\\ $(2,10)$\end{tabular} &
\begin{tabular}{@{}c@{}}Policy 3\\ $(2,15)$\end{tabular} &
\begin{tabular}{@{}c@{}}Policy 4\\ $(2,\infty)$\end{tabular} &
\begin{tabular}{@{}c@{}}Policy 5\\ $(2,N_5)$\end{tabular} &
\begin{tabular}{@{}c@{}}Policy 6\\ $(1,10)$\end{tabular} &
\begin{tabular}{@{}c@{}}Policy 7\\ $(1,15)$\end{tabular} \\
\hline
\multirow{3}{*}{$\dfrac{h_1}{\mu_1}\le \dfrac{h_2}{\mu_2}$}
& Max  & \textbf{\underline{0.00}} & 97.77 & 151.44 & 228.41 & 54.65 & 122.36 & 68.61 \\
& Mean & \textbf{\underline{0.00}} & 28.25 & 47.19 & 75.23 & 11.66 & 44.60 & 25.22 \\
& Std  & \textbf{\underline{0.00}} & 22.74 & 34.65 & 50.87 & 12.99 & 27.90 & 15.72 \\
\hline
\multirow{3}{*}{$\dfrac{h_2}{\mu_2}< \dfrac{h_1}{\mu_1}\le \dfrac{2h_2}{\mu_2}$}
& Max  & 46.72 & 29.46 & 44.65 & 92.40 & \textbf{7.09} & 71.07 & 45.37 \\
& Mean & 9.32 & 13.79 & 20.99 & 34.92 & \textbf{1.15} & 28.73 & 21.07 \\
& Std  & 11.06 & 6.96 & 10.97 & 22.70 & \textbf{1.94} & 14.62 & 9.38 \\
\hline
\multirow{3}{*}{$\dfrac{h_1}{\mu_1}> \dfrac{2h_2}{\mu_2}$}
& Max  & 114.64 & 56.75 & 34.50 & 81.53 & \textbf{26.00} & 65.83 & 84.25 \\
& Mean & 27.62 & 15.45 & 14.76 & 19.54 & \textbf{3.68} & 30.26 & 30.45 \\
& Std  & 28.95 & 12.18 & 7.56 & 20.70 & \textbf{5.82} & 13.92 & 17.22 \\
\hline

\end{tabular}
}
\end{table}

\begin{table}[htbp]
\centering
\renewcommand{\arraystretch}{0.95}
\setlength{\tabcolsep}{5pt}
\small
\captionsetup{labelfont={color=blue}, textfont={color=blue}}
\caption{Relative errors ($\%$) of heuristic policies for decision states $(20,j,k,\ell) \in \X_{\widehat{D}}$ when $\mu_1 < \mu_2 \leq 2\mu_1$, where $N_5$ is given in \eqref{eq:N5}.} \label{table:mid-Station2}
{\color{blue}
\begin{tabular}{|c|c|c|c|c|c|c|c|c|}
\hline
Cost regime & Error & 
\begin{tabular}{@{}c@{}}Policy 1\\ $(2,0)$\end{tabular} &
\begin{tabular}{@{}c@{}}Policy 2\\ $(2,10)$\end{tabular} &
\begin{tabular}{@{}c@{}}Policy 3\\ $(2,15)$\end{tabular} &
\begin{tabular}{@{}c@{}}Policy 4\\ $(2,\infty)$\end{tabular} &
\begin{tabular}{@{}c@{}}Policy 5\\ $(2,N_5)$\end{tabular} &
\begin{tabular}{@{}c@{}}Policy 6\\ $(1,10)$\end{tabular} &
\begin{tabular}{@{}c@{}}Policy 7\\ $(1,15)$\end{tabular} \\
\hline

\multirow{3}{*}{$\dfrac{h_1}{\mu_1}\le \dfrac{h_2}{\mu_2}$}
& Max  & 22.08 & 49.74 & 77.77 & 119.10 & \textbf{21.64} & 62.81 & 34.73 \\
& Mean & 5.26 & 18.53 & 25.54 & 35.00 & \textbf{2.96} & 19.91 & 12.68 \\
& Std  & 5.42 & 9.87 & 16.50 & 25.91 & \textbf{5.32} & 13.00 & 6.88 \\
\hline
\multirow{3}{*}{$\dfrac{h_2}{\mu_2}< \dfrac{h_1}{\mu_1}\le \dfrac{2h_2}{\mu_2}$}
& Max  & 49.07 & 30.92 & 41.64 & 63.14 & \textbf{\underline{0.00}} & 35.52 & 36.62 \\
& Mean & 18.95 & 18.18 & 18.16 & 19.12 & \textbf{\underline{0.00}} & 18.53 & 18.30 \\
& Std  & 10.96 & 6.97 & 9.01 & 14.62 & \textbf{\underline{0.00}} & 8.96 & 8.33 \\
\hline
\multirow{3}{*}{$\dfrac{h_1}{\mu_1}> \dfrac{2h_2}{\mu_2}$}
& Max  & 104.46 & 53.82 & 30.15 & 43.30 & \textbf{10.27} & 51.71 & 75.79 \\
& Mean & 32.78 & 19.72 & 14.53 & 10.12 & \textbf{0.75} & 22.07 & 26.94 \\
& Std  & 20.45 & 9.24 & 5.73 & 9.81 & \textbf{1.78} & 9.95 & 14.02 \\
\hline

\end{tabular}
}
\end{table}

\begin{table}[htbp]
\centering
\renewcommand{\arraystretch}{0.95}
\setlength{\tabcolsep}{5pt}
\small
\captionsetup{labelfont={color=blue}, textfont={color=blue}}
\caption{Relative errors ($\%$) of heuristic policies for decision states $(20,j,k,\ell) \in \X_{\widehat{D}}$ when $\mu_2 > 2\mu_1$, where $N_5$ is given in \eqref{eq:N5}.} \label{table:fast-Station2}
{\color{blue}
\begin{tabular}{|c|c|c|c|c|c|c|c|c|}
\hline
Cost regime & Error & 
\begin{tabular}{@{}c@{}}Policy 1\\ $(2,0)$\end{tabular} &
\begin{tabular}{@{}c@{}}Policy 2\\ $(2,10)$\end{tabular} &
\begin{tabular}{@{}c@{}}Policy 3\\ $(2,15)$\end{tabular} &
\begin{tabular}{@{}c@{}}Policy 4\\ $(2,\infty)$\end{tabular} &
\begin{tabular}{@{}c@{}}Policy 5\\ $(2,N_5)$\end{tabular} &
\begin{tabular}{@{}c@{}}Policy 6\\ $(1,10)$\end{tabular} &
\begin{tabular}{@{}c@{}}Policy 7\\ $(1,15)$\end{tabular} \\
\hline

\multirow{3}{*}{$\dfrac{h_1}{\mu_1}\le \dfrac{h_2}{\mu_2}$}
& Max  & 47.87 & 42.57 & 53.73 & 83.51 & \textbf{12.53} & 43.34 & 26.09 \\
& Mean & 12.97 & 19.54 & 20.77 & 21.04 & \textbf{1.24} & 13.05 & 11.75 \\
& Std  & 11.97 & 8.81 & 11.05 & 18.71 & \textbf{2.89} & 8.55 & 5.83 \\
\hline
\multirow{3}{*}{$\dfrac{h_2}{\mu_2}< \dfrac{h_1}{\mu_1}\le \dfrac{2h_2}{\mu_2}$}
& Max  & 60.48 & 45.09 & 33.90 & 45.65 & \textbf{1.37} & 34.71 & 43.97 \\
& Mean & 32.42 & 24.73 & 18.62 & 9.79 & \textbf{0.06} & 16.58 & 22.62 \\
& Std  & 12.07 & 8.28 & 7.59 & 10.34 & \textbf{0.19} & 8.25 & 8.69 \\
\hline
\multirow{3}{*}{$\dfrac{h_1}{\mu_1}> \dfrac{2h_2}{\mu_2}$}
& Max  & 529.04 & 282.70 & 162.09 & \textbf{9.26} & 134.91 & 254.20 & 376.48 \\
& Mean & 95.27 & 53.92 & 30.68 & \textbf{0.80} & 11.54 & 42.20 & 65.34 \\
& Std  & 103.70 & 52.71 & 28.69 & \textbf{2.02} & 25.81 & 51.88 & 75.96 \\
\hline

\end{tabular}
}
\end{table}

{\color{blue}
An immediate observation from Table \ref{table:slow-Station2} is that, when $\mu_1 \geq \mu_2$ and $\frac{h_1}{\mu_1} \leq \frac{h_2}{\mu_2}$, Policy 1 attains zero error. This is consistent with Statements \ref{state:empty_Station2_slower}, \ref{state:empty_Station2_equal}, and \ref{state:queue_in_Station2_slower} of Theorem \ref{thm:large_queue}, which imply that the optimal policy always routes to Station 1 in this regime. Moreover, Policy 1 continues to perform well in adjacent regimes, namely when either the cost condition shifts to $\frac{h_2}{\mu_2} < \frac{h_1}{\mu_1} \leq \frac{2h_2}{\mu_2}$ while $\mu_1 \geq \mu_2$ (Table \ref{table:slow-Station2}), or the rate condition shifts to $\mu_1 < \mu_2 \leq 2\mu_1$ while $\frac{h_1}{\mu_1} \leq \frac{h_2}{\mu_2}$ (Table \ref{table:mid-Station2}).
Likewise, Table \ref{table:mid-Station2} shows that Policy 5 attains zero error when $\mu_1 < \mu_2 \leq 2\mu_1$ and $\frac{h_2}{\mu_2} < \frac{h_1}{\mu_1} \leq \frac{2h_2}{\mu_2}$, in agreement with Statements \ref{state:empty_Station2_faster} and \ref{state:queue_in_Station2_slower}. Policy 5 also remains highly accurate in adjacent regimes, both within the same rate regime when the cost condition changes, and across neighboring rate regimes when the intermediate cost condition is preserved.
}

{\color{blue}
Consider now the parameter region in which $\mu_1 \geq \mu_2$ and $\frac{h_1}{\mu_1} > \frac{h_2}{\mu_2}$ unless stated otherwise. Recall that $i'_{j,k,\ell}$ in Theorem \ref{thm:large_queue} denotes an initial queue length beyond which the optimal action stabilizes, for fixed $j,k,\ell$ satisfying $j+k+\ell = 2$. \textcolor{blue}{When the inequality is strict, i.e., $\mu_1 > \mu_2$},} 
Statement \ref{state:empty_Station2_slower} and \ref{state:queue_in_Station2_slower} in Theorem \ref{thm:large_queue} establish the existence of a finite $i'_{j,k,\ell}$ such that an optimal decision-maker goes to Station 1 ($a=1$) if $i\geq i'_{j,k,\ell}$. While the result does not guarantee a monotone optimal policy in general (recall \eqref{eq:D_mono_i} need not hold), our numerical experiments strongly suggest that an optimal policy is in fact of \textcolor{blue}{\textbf{Station-2 threshold type with a finite threshold $N_{j,k,\ell}$}. This threshold may be zero; for example, $N_{1,0,1} = 0$ when $\frac{h_2}{\mu_2} < \frac{h_1}{\mu_1} \leq \frac{2h_2}{\mu_2}$, or it may take other finite values.} 
When $\mu_1 = \mu_2$, a similar result holds at states $(i,1,0,1)$; for other decision states $(i,1,1,0)$ and $(i,2,0,0)$, an optimal policy also follows a \textcolor{blue}{\textbf{Station-2 threshold type}}, but with an infinite threshold (see Statements \ref{state:empty_Station2_equal}). Consequently, when $\mu_1 \geq \mu_2$ and $\frac{h_1}{\mu_1} > \frac{h_2}{\mu_2}$, an \textcolor{blue}{optimal policy is numerically well described by a \textbf{Station-2 threshold structure}.}

{\color{blue}
Despite this seemingly simple threshold behavior, the cases differ substantially depending on whether or not Station 2 is available at the decision epoch.
In unblocked decision states $(i,1,1,0)$ and $(i,2,0,0)$, the numerically observed Station-2 thresholds $N_{j,k,0}$ vary systematically with the parameter configuration. In particular, the threshold increases with the relative downstream cost advantage of Station 2 and decreases with its time disadvantage (recall the time and cost differences in \eqref{eq:empty-Station2-time-cost}), while larger values of $h_0$ tend to reduce the threshold by making delayed return from Station 2 to upstream service more costly. For example, under $\mu_0 = 10, \mu_1 = 10, \mu_2 = 8, h_0 = 0.2, h_1 = 1, h_2 = 0.5$, we obatin $N_{1,1,0} = 4$. Varying one factor at a time illustrates this pattern: 
\begin{itemize} [leftmargin=2\parindent]
    \item [Cost:] Decreasing $h_2$ to $0.1$ yields $N_{1,1,0} = 8$, since Station 2 becomes less costly.
    \item [Time:] Decreasing $\mu_2$ to $4$ and $h_2$ to $0.25$ so that the ratio $\frac{h_2}{\mu_2}$ remains unchanged, yields $N_{1,1,0} = 2$, reflecting the larger service-time disadvantage of Station 2.
    \item [$h_0$:] Decreasing $h_0$ to $0.1$ yields $N_{1,1,0} = 31$, since a lower upstream holding cost makes delayed return to Station 0 less expensive.
\end{itemize}}
By contrast, in blocked states $(i,1,0,1)$, Station 2 is already occupied, so choosing Station 2 blocks the current server-job pair until service can begin, preventing the server from returning to provide upstream service. 
Accordingly, the \textcolor{blue}{Station-2 threshold $N_{1,0,1}$} of an optimal policy should be smaller than $N_{1,1,0}$ and $N_{2,0,0}$. 
{\color{blue}Moreover, because faster first-stage service compounds the downstream blocking effects at Station 2 (see the analysis in Section \ref{sec:main-results}), the threshold $N_{1,0,1}$ is additionally decreasing in $\mu_0$, beyond its dependence on the downstream cost and time differences and the upstream holding cost $h_0$.}

\textcolor{black}{Policy 5 explicitly accounts for this distinction by setting $N_{1,0,1} = 0$ and $N_{1,1,0}=N_{2,0,0}=\infty$, thereby avoiding blocking at states $(i,1,0,1)$ while fully favoring Station 2 at the other decision states. Although this routing rule is extreme, it performs best in the parameter regime where $\mu_1 \geq \mu_2$ and $\frac{h_1}{\mu_1} > \frac{h_2}{\mu_2}$, yielding} the smallest maximum error, mean error and standard deviation over decision states $(20,j,k,\ell) \in \X_{\widehat{D}}$; see Table \ref{table:slow-Station2}. 
{\color{blue}
The worst cases for Policy 5 arise when an optimal policy goes to Station 2 for all $(i,j,k,\ell)\in\X_D$ where $i \leq 20$, whereas Policy 5 instead routes to Station 1 for decision states $(i,1,0,1)$. For example, the maximum error $26.00\%$ occurs at state $(20,1,0,1)$ when $\mu_0 = 15,\mu_1 = 10, \mu_2 = 8, h_0 = 0.01, h_1 = 1, h_2 = 0.2$.}

By contrast, Policies 1--4 impose the same \textcolor{blue}{threshold} across all $(j,k,\ell)$ and therefore ignore the effect of blocking. Their performance depends strongly on how closely their fixed Station-2 threshold values match the effective thresholds of the optimal policy under a given parameter configuration. Since the routing rule becomes progressively more favorable to Station 2 from Policy 1 to Policy 4, \textcolor{blue}{Policies 1 and 2} perform poorly when the optimal policy routes to Station 2 over most of the relevant state space, whereas Policies 3 and 4 deteriorate when the optimal policy routes more frequently to Station 1. Thus, the worst cases for Policies 1–4 arise precisely when their fixed thresholds are misaligned with the state-dependent thresholds suggested by the optimal policy.

% The performance of Policies 1--4 with constant thresholds depends largely on the thresholds of optimal policies based on different input parameters. For instance, as \textcolor{blue}{Policy} 1 always decides to go to Station 1 and Policy 2 goes to Station 1 most of the time (except for states $(i,j,k,\ell)\in\X_D$ with $0\leq i < 10$), the worst performance cases of Policies 1 and 2 occur when the optimal policy goes to Station 2 for all $(i,j,k,\ell)\in\X_D$ where $i \leq 20$. In particular, the maximum error occurs for both when $\mu_0 = 5,\mu_1 = 10, \mu_2 = 8, h_0 = 0.01, h_1 = 1, h_2 = 0.2$ at state $(20,2,0,0)$. Under this circumstance, Policies 3 and 4 achieve better performance where the former has the error 31.2\% and the latter has no error since it coincides with the optimal \textcolor{blue}{policy} for all states $(i,j,k,\ell)\in\X_D$. Similarly, the worst cases for Policies 3 and 4 happen when the optimal policy goes to Station 1 more often. 
% For example, the maximum error for Policy 4 occurs at state $(20,2,0,0)$ when $\mu_0 = 40,\mu_1 = 10, \mu_2 = 8, h_0 = 1, h_1 = 1, h_2 = 0.3$, where the optimal policy always routes to Station 1 except $(0,1,0,0), (0,1,1,0), (1,1,1,0)$ and $(0,2,0,0)$ . Similarly, the worst-case error for Policy 3 occurs at the same state when $h_0 = 0.01, h_1 = 1, h_2  = 0.7$ for the same rates, where an optimal policy always routes to Station 1 for states $(i,1,0,1)$ for $i\leq 20$.
{\color{blue} 
Consider now $\mu_2 > \mu_1$, unless stated otherwise. We first note that the structural behavior at the unblocked decision states is symmetric to the previously discussed case $\mu_1 \geq \mu_2$. In place of the Station-2 threshold structure observed there, the optimal policy is now numerically well described by a Station-1 threshold structure, with $N_{j,k,0} = 0$ when $\frac{h_1}{\mu_1} \geq \frac{h_2}{\mu_2}$, and with other finite threshold values otherwise. This is captured in Statement \ref{state:empty_Station2_faster} in Theorem \ref{thm:large_queue}.
By contrast, at blocked decision states $(i,1,0,1)$, the relative speeds of upstream and downstream service create a more intricate operational interaction, substantially complicating the structure of the optimal policy. 
In particular, a large first-stage service rate $\mu_0$ compounds downstream blocking effects, whereas a small $\mu_0$ tends to limit blocking to the current server-job pair. Consequently, unlike the case $\mu_1 \geq \mu_2$, where the optimal policy is numerically described as a Station-2 threshold-type structure, such a characterization does not generally hold. 
To illustrate this phenomenon, we examine how the optimal policy at $(i,1,0,1)$ varies under different parameter configurations corresponding to the signs of the time and cost differences in \eqref{eq:queue-in-Station2-time-cost}. Throughout, we fix $\mu_1 = 10$ and $h_1 = 1$ and vary the remaining parameters.
\begin{enumerate} [label= \textbf{Case} \arabic*:, leftmargin=2.5\parindent, ref= \arabic*]
    \item 
    Suppose $\mu_1 < \mu_2 \leq 2 \mu_1$, i.e., $\frac{1}{\mu_2}<\frac{1}{\mu_1} \leq \frac{2}{\mu_2}$.
    \begin{enumerate} [label= \textbf{Subcase} (\alph*):, leftmargin=2\parindent, ref= \textbf{Subcase} \theenumi (\alph*)]
        \item \label{example:non-monotone_slow_station2}
        If $\frac{h_1}{\mu_1} \leq \frac{2h_2}{\mu_2}$, Statement \ref{state:queue_in_Station2_slower} of Theorem \ref{thm:large_queue} establishes that $a=1$ for all $i$.
        \item If $\frac{h_1}{\mu_1} > \frac{2h_2}{\mu_2}$, non-monotonic behavior can arise. For example, under $\mu_0 = 0.1, \mu_2 = 18, h_0 = 0.1, h_2 = 0.8824$, the optimal policy is: 
    \begin{align*}
        a & =\begin{cases} 
        2 &\text{ if } i=0 \text{ or } 2 \leq i < 4,\\
        1 &\text{ if } i=1 \text{ or } i \geq 4.
       \end{cases}
    \end{align*}
    Numerically, this pattern appears only under rather extreme parameter configurations and is typically confined to very small initial queue lengths, after which the policy stabilizes quickly.
    Increasing the initial service rate to $\mu_0 = 10$ restores monotonicity:
     \begin{align*}
        a & =\begin{cases} 
        2 &\text{ if } i=0,\\
        1 &\text{ if } i \geq 1.
       \end{cases}
    \end{align*}
    In both cases, the policy eventually stabilizes to $a=1$ as $i$ increases, consistent with Statement \ref{state:queue_in_Station2_slower} of Theorem \ref{thm:large_queue}. 
    \end{enumerate}
    \item Now consider $\mu_2 > 2\mu_1$, i.e., $\frac{1}{\mu_1} > \frac{2}{\mu_2}$. 
    Statment \ref{state:queue_in_Station2_faster_slow_initial} in Theorem \ref{thm:large_queue} confirms that if $\mu_0$ is sufficiently small, then $a=2$ for all large $i$, whereas Statment \ref{state:queue_in_Station2_faster_fast_initial} shows that if $\mu_0$ is sufficiently large, then $a=1$ for all large $i$.
    \begin{enumerate} [label= \textbf{Subcase} (\alph*):, leftmargin=2\parindent, ref= \textbf{Subcase} \theenumi (\alph*)]
        \item Suppose $\frac{h_1}{\mu_1} \leq \frac{2h_2}{\mu_2}$. For example, with parameters $\mu_0 = 1, \mu_2 = 28, h_0 = 0.1, h_2 = 1.8667$. The optimal policy is indeed monotone:
        \begin{align*}
            a & =\begin{cases} 
            1 &\text{ if } i < 26,\\
            2 &\text{ if } i \geq 26.
           \end{cases}
        \end{align*}
        Increasing $\mu_0$ to $\mu_0 = 100$ collapses the policy to $a = 1$ for all $i$.
        \item \label{example:non-monotone_fast_station2}
        Now consider $\frac{h_1}{\mu_1} > \frac{2h_2}{\mu_2}$. With $\mu_0 = 0.2, \mu_2 = 28, h_0 = 0.1, h_2 = 1.12$, the optimal policy is $a = 2$ for all $i$. Increasing $\mu_0$ again shifts the optimal policy to Station 1, and intermediate values of $\mu_0$ may induce non-monotonic behavior. For instance, when $\mu_0 = 20$, the policy becomes
        \begin{align*}
            a & =\begin{cases} 
            2 &\text{ if } i \leq 2 \text{ or } i > 15,\\
            1 &\text{ if } 2< i \leq 15.
           \end{cases}
        \end{align*}
        Further increasing $\mu_0$ to $\mu_0 = 200$ yields a threshold-type policy again:
        \begin{align*}
            a & =\begin{cases} 
            2 &\text{ if } i=0,\\
            1 &\text{ if } i \geq 1.
           \end{cases}
        \end{align*}
    \end{enumerate}
\end{enumerate}
These examples demonstrate that, under general service rates, optimal policies need not exhibit monotone or threshold behavior, even though they eventually stabilize for large $i$. A key takeaway is that higher values of $\mu_0$ strengthen the preference for Station 1, as faster upstream service amplifies downstream blocking effects associated with Station 2.
}

{\color{blue}
These structural properties help explain the heuristic performance reported in Tables \ref{table:mid-Station2} and \ref{table:fast-Station2}. Consider first the regime $\mu_1 < \mu_2 \leq 2\mu_1$ in Table \ref{table:mid-Station2}. When, in addition, $\frac{h_2}{\mu_2} < \frac{h_1}{\mu_1} \leq \frac{2h_2}{\mu_2}$, then Policy 5 is optimal. When $\frac{h_1}{\mu_1} \leq \frac{h_2}{\mu_2}$, Policies 1 and 5 perform best, as both agree with the optimal action at states $(i,1,0,1)$. Their relative performance is then largely determined by how the test initial queue length 20 compares with the effective Station-1 thresholds at decision states $(i,j,k,0) \in \X_{\widehat{D}}$. 
In particular, if $N_{j,k,0} \geq 20$, the optimal policy coincides with Policy 1 throughout the tested range $i \leq 20$; otherwise, Policy 5 may perform better, as it deviates from the optimal policy only at smaller initial queue lengths.
Finally, when $\frac{h_1}{\mu_1} > \frac{2h_2}{\mu_2}$, Policies 4 and 5 align with the optimal policy at $(i,j,k,0) \in \X_{\widehat{D}}$, and their relative performance is then governed by how the stabilizing threshold $i'_{1,0,1}$ compares with the test initial queue length 20. 

When $\mu_2 > 2\mu_1$, as in Table \ref{table:fast-Station2}, the policy structure at states $(i,1,0,1)$ becomes even more intricate because the immediate time and cost advantages of Station 2 must be weighed against blocking effects that are amplified by the initial service rate $\mu_0$. 
Consequently, heuristic performance becomes especially sensitive to the balance among these competing effects. Policy 5 often remains competitive because it accounts for blocking, but its advantage diminishes when the optimal policy favors Station 2 even at $(i,1,0,1)$. 
Likewise, the fixed-threshold Policies 1--4 become increasingly sensitive to threshold misalignment, since a single threshold cannot capture the more complex dependence on $\mu_0$ and the downstream cost and time tradeoffs.

In summary, the numerical study shows that heuristic performance is highly sensitive to the parameter regime, which in turn highlights the value of exploiting structural insights from the optimal policy when designing practical routing rules. From an operational perspective, the routing decision should balance the immediate service-time and cost differences between the two downstream options against the blocking effects associated with Station 2. In general, a station becomes more attractive when its service is faster or less costly. In addition, increasing the upstream service rate $\mu_0$ strengthens the preference for Station 1 by amplifying the blocking consequences of routing to Station 2, while increasing the upstream holding cost $h_0$ tends to reduce the queue-length thresholds at which the decision switches. 
}

\section{Conclusions} \label{sec:conc}
A holding-cost clearing system model of a two-stage tandem queue with two servers was considered. Upon completing an initial screening service at the upstream stage, each server decides whether to \textcolor{black}{proceed with the job to a parallel-service station or to a single-service station} at the downstream stage. \textcolor{blue}{We formulated this sequential decision problem as a Markov decision process and established structural properties of the optimal policy based on the relative downstream service speeds, focusing on large initial queues. Numerical experiments examined how an optimal policy changes with system parameters and evaluated the performance of several natural heuristics, highlighting the importance of understanding the policy structures.} This paper motivates several directions for future study:
\begin{enumerate}
    \item \textbf{Increasing the number of servers}. The current study is confined to only two servers to make the analysis more tractable. A system with more flexible servers that move with the jobs and more dedicated servers for special requests is of interest. For example in healthcare systems, there are often multiple flexible servers (nurse practitioners) that accompany the patients throughout their visit. There is a dedicated space for doctors where collaborative service (with the nurses) can be provided.
    \item \textbf{Designing a better heuristic policy}. As demonstrated in the numerical study, the optimal structure differs significantly depending on whether \textcolor{black}{$\mu_1 \geq \mu_2$ or $\mu_1 < \mu_2$}. 
    In particular, when \textcolor{black}{$\mu_1 < \mu_2$}, the \textcolor{black}{optimal policy does not admit a threshold structure}.
    Furthermore, even though \textcolor{black}{an optimal policy is numerically well described by a threshold structure in $i$ when $\mu_1 \geq \mu_2$}, the threshold itself should depend on the state variables $j,k,\ell$. For software implementation (e.g., in a manufacturing setting), heuristics are primarily useful when computational constraints arise, such as with a large initial queue. In any system with human servers, simple yet well-designed heuristics that effectively leverage structural properties become even more critical to decision-makers and should be explored.
\end{enumerate}

\appendix
\section{Appendix} \label{sec:appendix}
To establish the structure of an optimal policy, we first present the following preliminary results, which provide bounds on the cost differences between various initial states. 
% Notably, unlike the main results, these preliminaries do not rely on Assumption \ref{assm:cost-to-serve}.
\subsection{Preliminaries for Theorem \ref{thm:large_queue}}
\begin{lemma} \label{lemma:general_rates}
The difference in value functions when evaluated at various starting states yields several bounds. 
    \begin{enumerate}
        \item \label{state1:general_rates} Consider $(0,j,k, \ell) \in \X$ where $j+k+\ell = 1$. 
            \begin{equation} \label{eq:empty-queue-bound1}
                v(0,j,k+1,\ell) - v(0,j,k,\ell) \geq \frac{h_1}{\mu_1},
            \end{equation}
            \begin{equation} \label{eq:empty-queue-bound2}
                v(0,j,k,\ell+1) - v(0,j,k,\ell) \geq \frac{h_2}{\mu_2}, 
            \end{equation}
            and 
            \begin{equation} \label{eq:mono_bdry}
                v(0,j+1,k,\ell) - v(0,j,k,\ell) \geq \frac{h_0}{\mu_0} + \min\left\{\frac{h_1}{\mu_1},\frac{h_2}{\mu_2}\right\}.
            \end{equation}
        \item \label{state2:general_rates}
        Consider $s = (i,j,k,\ell) \in \X$ such that $i \geq 0$ and $j+k+\ell = 2$. 
            \begin{equation}  \label{eq:mono}
                v(i+1,j,k,\ell) - v(i,j,k,\ell) \geq \frac{h_0}{\mu_0} + \min\left\{\frac{h_1}{\mu_1},\frac{h_2}{\mu_2}\right\}.
            \end{equation}
            In addition, if $j\geq 1$,
            \begin{align}
                v(i+1,j-1,k+1,\ell) - v(i,j,k,\ell) & \geq  \min\left\{\frac{h_1}{\mu_1},\frac{h_2}{\mu_2}\right\}. \label{eq:onemore1}
            \end{align}
            \begin{align}
                v(i+1,j-1,k,\ell+1) - v(i,j,k,\ell) \geq  \min\left\{\frac{h_1}{\mu_1},\frac{h_2}{\mu_2}\right\}. \label{eq:onemore2}
            \end{align}
        \item \label{state3:general_rates}
        Suppose $\mu_1 \geq \mu_2$ and consider $s = (i,j,k,\ell) \in \X_{\widehat{D}}$. 
            \begin{equation}\label{eq:decision_bd1}
                D(i,j,k,\ell) = v(i,j-1,k+1,\ell) - v(i,j-1,k,\ell+1) \leq \frac{h_1}{\mu_1} - \frac{h_2}{\mu_2}.
            \end{equation}
            \item \label{state4:general_rates}
            \textcolor{blue}{Suppose $\mu_1 \leq \mu_2$} and consider $s = (i,j,k,\ell) \in \X_{\widehat{D}}$.
            \begin{align}
                \lefteqn{D(i,j,k,\ell) = v(i,j-1,k+1,\ell) - v(i,j-1,k,\ell+1)}& \nonumber \\
                &\quad \leq \frac{h_1}{\mu_1} - \frac{(\ell+1)h_2}{\mu_2} + \frac{(i+j-1)h_0}{2}\left(\frac{1}{\mu_1} - \frac{\ell+1}{\mu_2}\right). \label{eq:decision_bd2}
        \end{align}
        As an intermediate result in proving \eqref{eq:decision_bd2}, we also have
        \begin{equation}
            \textcolor{blue}{\min\{v(i+1,0,2,0), v(i+1,0,1,1)\}}-v(i,2,0,0) \leq \frac{h_1}{\mu_1} + \frac{(i+2)h_0}{2}\left(\frac{1}{\mu_1} - \frac{1}{\mu_0}\right). \label{eq:decision_bd2_inte}
        \end{equation}
    \end{enumerate}
\end{lemma}
\begin{proof} [Proof of Statement \ref{state1:general_rates}]
    A little algebra in \textcolor{blue}{the optimality equations \eqref{eq:opt-zero}} yield the values of the following boundary cases:
    \begin{align*}
        \begin{split}
            v(0,0,1,0) &= \frac{h_1}{\mu_1}\\
            v(0,0,0,1) &= \frac{h_2}{\mu_2}\\
            v(0,1,0,0) 
            &=\frac{h_0}{\mu_0}+\min\left\{\frac{h_1}{\mu_1}, \frac{h_2}{\mu_2}\right\}\\
        \end{split}
        \begin{split}
            v(0,0,2,0) 
            &=\frac{2h_1}{\mu_1}\\
            v(0,0,0,2) 
            &=\frac{3h_2}{\mu_2}\\
            v(0,0,1,1) 
            &=\frac{h_1}{\mu_1} + \frac{h_2}{\mu_2}\\
        \end{split}
    \end{align*}
    \begin{align*}
         v(0,1,1,0) 
            &=\frac{h_0}{\mu_0} + \frac{h_1}{\mu_1} + \min\left\{\frac{h_1}{\mu_1}, \frac{h_2}{\mu_2}\right\}\\
            v(0,1,0,1) 
            &=\frac{h_0}{\mu_0} + \frac{h_2}{\mu_2}+\frac{\mu_0}{{\mu_0+\mu_2}}\min\{\frac{h_1}{\mu_1},\frac{2h_2}{\mu_2}\}+\frac{\mu_2}{{\mu_0+\mu_2}}\min\left\{\frac{h_1}{\mu_1}, \frac{h_2}{\mu_2}\right\}\\
            &\geq \frac{h_0}{\mu_0} + \frac{h_2}{\mu_2}+\frac{\mu_0}{{\mu_0+\mu_2}}\min\left\{\frac{h_1}{\mu_1},\frac{h_2}{\mu_2}\right\}+\frac{\mu_2}{{\mu_0+\mu_2}}\min\left\{\frac{h_1}{\mu_1}, \frac{h_2}{\mu_2}\right\}\\
        &= \frac{h_0}{\mu_0} + \frac{h_2}{\mu_2}+\min\left\{\frac{h_1}{\mu_1},\frac{h_2}{\mu_2}\right\}\\
        v(0,2,0,0)
        &=\frac{2h_0}{2\mu_0} + \min\{v(0,1,1,0), v(0,1,0,1)\}\\
        &\geq \frac{h_0}{\mu_0} + \min\left\{\frac{h_0}{\mu_0} + \frac{h_1}{\mu_1} + \min\left\{\frac{h_1}{\mu_1}, \frac{h_2}{\mu_2}\right\},\frac{h_0}{\mu_0} + \frac{h_2}{\mu_2} + \min\left\{\frac{h_1}{\mu_1}, \frac{h_2}{\mu_2}\right\}\right\}\\
        &=\frac{2h_0}{\mu_0} + 2\min\left\{\frac{h_1}{\mu_1}, \frac{h_2}{\mu_2}\right\}
    \end{align*}
    It suffices to show the following three cases:
    \begin{enumerate}[label= \textbf{Case} \arabic*:, leftmargin=3\parindent]
    \item If $(j,k,\ell) = (0,1,0)$, the three inequalities respectively become
    \begin{align*}
        v(0,0,2,0) - v(0,0,1,0) &= \frac{h_1}{\mu_1},\\
        v(0,0,1,1) - v(0,0,1,0) &= \frac{h_2}{\mu_2},\\
        v(0,1,1,0) - v(0,0,1,0) &= \frac{h_0}{\mu_0}+\min\left\{\frac{h_1}{\mu_1}, \frac{h_2}{\mu_2}\right\}.
    \end{align*}
    \item If $(j,k,\ell) = (0,0,1)$, the three inequalities respectively become
    \begin{align*}
        v(0,0,1,1) - v(0,0,0,1) &= \frac{h_1}{\mu_1},\\
        v(0,0,0,2) - v(0,0,0,1) &= \frac{2h_2}{\mu_2} \geq \frac{h_2}{\mu_2},\\
        v(0,1,0,1) - v(0,0,0,1)
        &\geq \frac{h_0}{\mu_0}+\min\left\{\frac{h_1}{\mu_1}, \frac{h_2}{\mu_2}\right\}.
    \end{align*}
    \item If $(j,k,\ell) = (0,0,1)$, the three inequalities respectively become
    \begin{align*}
        v(0,1,1,0) - v(0,1,0,0) &= \frac{h_1}{\mu_1},\\
        v(0,1,0,1) - v(0,1,0,0) &\geq \frac{h_2}{\mu_2},\\
        v(0,2,0,0) - v(0,1,0,0)
        &\geq \frac{h_0}{\mu_0}+\min\left\{\frac{h_1}{\mu_1}, \frac{h_2}{\mu_2}\right\}.
    \end{align*}
    \end{enumerate}
\end{proof}
\begin{proof} [Proof of Statement \ref{state2:general_rates}]
    We prove \eqref{eq:mono} first, which is equivalent to showing the following
    \begin{align}
        v(i+1,0,2,0) - v(i,0,2,0) &\geq \frac{h_0}{\mu_0}+\min\left\{\frac{h_1}{\mu_1}, \frac{h_2}{\mu_2}\right\}, \label{eq:mono_020}\\
        v(i+1,0,1,1) - v(i,0,1,1) &\geq \frac{h_0}{\mu_0}+\min\left\{\frac{h_1}{\mu_1}, \frac{h_2}{\mu_2}\right\},\label{eq:mono_011}\\
        v(i+1,0,0,2) - v(i,0,0,2) &\geq \frac{h_0}{\mu_0}+\min\left\{\frac{h_1}{\mu_1}, \frac{h_2}{\mu_2}\right\},\label{eq:mono_002}\\
        v(i+1,1,1,0) - v(i,1,1,0) &\geq \frac{h_0}{\mu_0}+\min\left\{\frac{h_1}{\mu_1}, \frac{h_2}{\mu_2}\right\},\label{eq:mono_110}\\
        v(i+1,1,0,1) - v(i,1,0,1) &\geq \frac{h_0}{\mu_0}+\min\left\{\frac{h_1}{\mu_1}, \frac{h_2}{\mu_2}\right\},\label{eq:mono_101}\\
        v(i+1,2,0,0) - v(i,2,0,0) &\geq \frac{h_0}{\mu_0}+\min\left\{\frac{h_1}{\mu_1}, \frac{h_2}{\mu_2}\right\}\label{eq:mono_200},
    \end{align}for all $i\geq0$.
    
    Proof by induction on $i$. When $i = 0$, \textcolor{blue}{applying the optimality equations \eqref{eq:opt1} and \eqref{eq:opt-zero} at states $(1,0,2,0)$ and $(0,0,2,0)$, respectively,} the first inequality \eqref{eq:mono_020} becomes
    \begin{align*}
        \lefteqn{v(1,0,2,0) - v(0,0,2,0)}&\\
        &\quad = \textcolor{blue}{\left(\frac{h_0+2h_1}{2\mu_1} + v(0,1,1,0)\right) - \left(\frac{2h_1}{2\mu_1} + v(0,0,1,0)\right)}\\
        &\quad = \frac{h_0}{2\mu_1} + v(0,1,1,0) - v(0,0,1,0)\\
        &\quad \geq \frac{h_0}{\mu_0}+\min\left\{\frac{h_1}{\mu_1}, \frac{h_2}{\mu_2}\right\},
    \end{align*}where the inequality follows from the non-negativity of the term $\frac{h_0}{2\mu_1}$ and \eqref{eq:mono_bdry} in Statement \ref{state1:general_rates}. The proof for \eqref{eq:mono_011} and \eqref{eq:mono_002} follows similarly by using \eqref{eq:mono_bdry} and is omitted for brevity. For \eqref{eq:mono_110} at $i=0$,
    using the non-negativity of the term $\frac{h_0}{\mu_0+\mu_1}$ and \eqref{eq:mono_bdry} yields
    \begin{align*}
        \lefteqn{v(1,1,1,0) - v(0,1,1,0)}&\\
        &\quad = \frac{h_0}{\mu_0+\mu_1} + \frac{\mu_1}{\mu_0+\mu_1} \Big[v(0,2,0,0) - v(0,1,0,0)\Big]\\
        &\qquad + \frac{\mu_0}{\mu_0+\mu_1} \Big[\min\{v(1,0,2,0), v(1,0,1,1)\} - \min\{v(0,0,2,0), v(0,0,1,1)\}\Big]
        \\
        &\quad \geq \frac{\mu_1}{\mu_0+\mu_1} \left(\frac{h_0}{\mu_0}+\min\left\{\frac{h_1}{\mu_1}, \frac{h_2}{\mu_2}\right\}\right)\\
        &\qquad + \frac{\mu_0}{\mu_0+\mu_1} \Big[\min\{v(1,0,2,0), v(1,0,1,1)\} - \min\{v(0,0,2,0), v(0,0,1,1)\}\Big].
    \end{align*}

There are two cases to consider based on the optimal decision at state $(1,1,1,0)$.
\begin{enumerate}[label= \textbf{Case} \arabic*:, leftmargin=3\parindent]
\item If $v(1,0,2,0) \leq v(1,0,1,1)$, substituting the second minimum with an upper bound $v(0,0,2,0)$ yields
\begin{align*}
    \lefteqn{v(1,1,1,0) - v(0,1,1,0)}&\\
    &\quad \geq \frac{\mu_1}{\mu_0+\mu_1} \left(\frac{h_0}{\mu_0}+\min\left\{\frac{h_1}{\mu_1}, \frac{h_2}{\mu_2}\right\}\right) + \frac{\mu_0}{\mu_0+\mu_1} \Big[v(1,0,2,0) -v(0,0,2,0)\Big]\\
    &\quad \geq \frac{\mu_1}{\mu_0+\mu_1} \left(\frac{h_0}{\mu_0}+\min\left\{\frac{h_1}{\mu_1}, \frac{h_2}{\mu_2}\right\}\right) + \frac{\mu_0}{\mu_0+\mu_1} \left(\frac{h_0}{\mu_0}+\min\left\{\frac{h_1}{\mu_1}, \frac{h_2}{\mu_2}\right\}\right)\\
    &\quad = \frac{h_0}{\mu_0}+\min\left\{\frac{h_1}{\mu_1}, \frac{h_2}{\mu_2}\right\},
\end{align*}where the last inequality holds by \eqref{eq:mono_020} at $i=0$.
\item If $v(1,0,2,0) > v(1,0,1,1)$, substituting the second minimum with an upper bound $v(0,0,1,1)$ yields
\begin{align*}
    \lefteqn{v(1,1,1,0) - v(0,1,1,0)}&\\
    &\quad \geq \frac{\mu_1}{\mu_0+\mu_1} \left(\frac{h_0}{\mu_0}+\min\left\{\frac{h_1}{\mu_1}, \frac{h_2}{\mu_2}\right\}\right) + \frac{\mu_0}{\mu_0+\mu_1} \Big[v(1,0,1,1) -v(0,0,1,1)\Big]\\
    &\quad \geq \frac{\mu_1}{\mu_0+\mu_1} \left(\frac{h_0}{\mu_0}+\min\left\{\frac{h_1}{\mu_1}, \frac{h_2}{\mu_2}\right\}\right) + \frac{\mu_0}{\mu_0+\mu_1} \left(\frac{h_0}{\mu_0}+\min\left\{\frac{h_1}{\mu_1}, \frac{h_2}{\mu_2}\right\}\right)\\
    &\quad = \frac{h_0}{\mu_0}+\min\left\{\frac{h_1}{\mu_1}, \frac{h_2}{\mu_2}\right\},
\end{align*}where the last inequality holds by \eqref{eq:mono_011} at $i=0$.
\end{enumerate}
The proof for \eqref{eq:mono_101} at $i=0$ follows analogously by using \eqref{eq:mono_bdry}, and the established results of \eqref{eq:mono_011} and \eqref{eq:mono_002} to discuss the two cases. It remains to check \eqref{eq:mono_200} at $i=0$. By noticing the non-negativity of $\frac{h_0}{2\mu_0}$,
\begin{align}
    \lefteqn{v(1,2,0,0) - v(0,2,0,0)}& \nonumber \\
    &\quad = \frac{h_0}{2\mu_0} + \min\{v(1,1,1,0), v(1,1,0,1)\} - \min\{v(0,1,1,0), v(0,1,0,1)\} \nonumber \\
    &\quad \geq \min\{v(1,1,1,0), v(1,1,0,1)\} - \min\{v(0,1,1,0), v(0,1,0,1)\} \label{eq:mono_200_inte1}
\end{align}
Similarly, as it was in proving \eqref{eq:mono_110}, there are also two cases to consider in \eqref{eq:mono_200_inte1} \textcolor{blue}{depending} on whether $v(1,1,1,0) \leq v(1,1,0,1)$. A comparable argument applies to find the desired lower bound of the difference between the two minima, by letting the state $(0,2,0,0)$ follow the same (potentially sub-optimal) action as an optimal decision of the state $(1,2,0,0)$. Mathematically, there is a slightly more concise way to present this process.
\begin{small}
\begin{align}
    \lefteqn{\min\{v(1,1,1,0), v(1,1,0,1)\} - \min\{v(0,1,1,0), v(0,1,0,1)\}}& \nonumber \\
    &\quad = \min\Big\{v(1,1,1,0)- \min\{v(0,1,1,0), v(0,1,0,1)\}, v(1,1,0,1)- \min\{v(0,1,1,0), v(0,1,0,1)\}\Big\} \nonumber \\
    &\quad \geq \min\{v(1,1,1,0)- v(0,1,1,0),v(1,1,0,1) - v(0,1,0,1)\} \nonumber \\
    &\quad \geq \min\left\{\frac{h_0}{\mu_0}+\min\left\{\frac{h_1}{\mu_1}, \frac{h_2}{\mu_2}\right\},\frac{h_0}{\mu_0}+\min\left\{\frac{h_1}{\mu_1}, \frac{h_2}{\mu_2}\right\}\right\} \nonumber \\
    &\quad = \frac{h_0}{\mu_0}+\min\left\{\frac{h_1}{\mu_1}, \frac{h_2}{\mu_2}\right\} \label{eq:mono_200_inte2},
\end{align}
\end{small}where the first inequality holds by replacing $\min\{v(0,1,1,0), v(0,1,0,1)\}$ with $v(0,1,1,0)$ and $v(0,1,0,1)$, respectively, for the first and second position of the outside-most minimum. The second inequality applies the established results \eqref{eq:mono_110} and \eqref{eq:mono_101} at $i=0$. Substituting the lower bound given by \eqref{eq:mono_200_inte2} in \eqref{eq:mono_200_inte1} yields the desired inequality.
% \begin{enumerate}[label= \textbf{Case} \arabic*:, leftmargin=3\parindent]
% \item If $v(1,1,1,0) \leq v(1,1,0,1)$, select an upper bound $(0,1,1,0)$ in the second minimum.
% \begin{align*}
%     \lefteqn{v(1,2,0,0) - v(0,2,0,0)}&\\
%     &\quad \geq v(1,1,1,0) - v(0,1,1,0)\\
%     &\quad \geq \frac{h_0}{\mu_0}+\min\left\{\frac{h_1}{\mu_1}, \frac{h_2}{\mu_2}\right\}
% \end{align*}by \eqref{eq:mono_110} at $i=0$.
% \item If $v(1,1,1,0) > v(1,1,0,1)$, select an upper bound $(0,1,0,1)$ in the second minimum.
% \begin{align*}
%     \lefteqn{v(1,2,0,0) - v(0,2,0,0)}&\\
%     &\quad \geq v(1,1,0,1) - v(0,1,0,1)\\
%     &\quad \geq \frac{h_0}{\mu_0}+\min\left\{\frac{h_1}{\mu_1}, \frac{h_2}{\mu_2}\right\}
% \end{align*}by \eqref{eq:mono_101} at $i=0$.
% \end{enumerate}

Now assume all the results hold at $i-1$, consider them at $i \geq 1$. By the non-negativity of the term $\frac{h_0}{2\mu_1}$ and the inductive hypothesis of \eqref{eq:mono_110}, the first inequality \eqref{eq:mono_020} becomes
\begin{align*}
    \lefteqn{v(i+1,0,2,0) - v(i,0,2,0)}&\\
    &\quad = \frac{h_0}{2\mu_1} + v(i,1,1,0) - v(i-1,1,1,0)\\
        &\quad \geq \frac{h_0}{\mu_0}+\min\left\{\frac{h_1}{\mu_1}, \frac{h_2}{\mu_2}\right\},
\end{align*}which verifies \eqref{eq:mono_020}. A comparable reasoning holds for \eqref{eq:mono_011} using the inductive hypothesis of \eqref{eq:mono_110} and \eqref{eq:mono_101}, and \eqref{eq:mono_002} using that of \eqref{eq:mono_101}.

For \eqref{eq:mono_110}, considering the non-negativity of the term $\frac{h_0}{\mu_0+\mu_1}$ and the inductive hypothesis of \eqref{eq:mono_200} yields
    \begin{align*}
        \lefteqn{v(i+1,1,1,0) - v(i,1,1,0)}&\\
        &\quad = \frac{h_0}{\mu_0+\mu_1} + \frac{\mu_1}{\mu_0+\mu_1} \Big[v(i,2,0,0) - v(i-1,2,0,0)\Big]\\
        &\qquad + \frac{\mu_0}{\mu_0+\mu_1} \Big[\min\{v(i+1,0,2,0), v(i+1,0,1,1)\} - \min\{v(i,0,2,0), v(i,0,1,1)\}\Big]
        \\
        &\quad \geq \frac{\mu_1}{\mu_0+\mu_1} \left(\frac{h_0}{\mu_0}+\min\left\{\frac{h_1}{\mu_1}, \frac{h_2}{\mu_2}\right\}\right)\\
        &\qquad + \frac{\mu_0}{\mu_0+\mu_1} \Big[\min\{v(i+1,0,2,0), v(i+1,0,1,1)\} - \min\{v(i,0,2,0), v(i,0,1,1)\}\Big].
    \end{align*}
The remaining discussion is analogous to that of proving \eqref{eq:mono_110} at $i=0$ by utilizing the established results of \eqref{eq:mono_020} and \eqref{eq:mono_011} for the analysis of two cases, which is omitted here for brevity. A comparable argument applies to prove \eqref{eq:mono_101} by using the inductive hypothesis of \eqref{eq:mono_200} at $i-1$ and the verified results of \eqref{eq:mono_011} and \eqref{eq:mono_002}. Finally, the proof of \eqref{eq:mono_200} for the generic $i$ is similar to that at $i=0$ using \eqref{eq:mono_110} and \eqref{eq:mono_101} that have been proved.

It remains to show \eqref{eq:onemore1} and \eqref{eq:onemore2}, where $(j,k,\ell)$ takes values from $(1,1,0)$, $(1,0,1)$ and $(2,0,0)$. Now consider \eqref{eq:onemore1}, which is implied by the following three inequalities
\begin{align}
    v(i+1,0,2,0) - v(i,1,1,0) &\geq \min\left\{\frac{h_1}{\mu_1}, \frac{h_2}{\mu_2}\right\}, \label{eq:onemore1_110}\\
    v(i+1,0,1,1) - v(i,1,0,1) &\geq \min\left\{\frac{h_1}{\mu_1}, \frac{h_2}{\mu_2}\right\},\label{eq:onemore1_101}\\
    v(i+1,1,1,0) - v(i,2,0,0) &\geq \min\left\{\frac{h_1}{\mu_1}, \frac{h_2}{\mu_2}\right\}\label{eq:onemore1_200}.
\end{align}

Proof by induction on $i$. Consider first \eqref{eq:onemore1_110} at $i=0$. 
{\color{blue}
The optimality equation \eqref{eq:opt1} evaluated at $(1,0,2,0)$ can be rearranged as
\begin{align}
    & v(1,0,2,0) = \frac{h_0 + 2 h_1}{\mu_0+2\mu_1} + \frac{2\mu_1}{\mu_0+2\mu_1} v(0,1,1,0) + \frac{\mu_0}{\mu_0+2\mu_1} v(1,0,2,0). \label{eq:thinning_i020}
\end{align}
This representation can also be obtained through a standard thinning argument applied at state $(1,0,2,0)$ under an optimal policy. Since the effective total transition rate out of $(1,0,2,0)$ is $d(0,2,0) = 2\mu_1$, thinning at a higher rate $\mu_0 + 2\mu_1$ introduces a self-transition (dummy transition) with rate $\mu_0$, under which the state remains unchanged. Therefore, equation \eqref{eq:thinning_i020} admits the same interpretation as the original optimality equation: the value equals the sum of the expected one-step holding cost and the expected continuation value after the next event, except that the next event may now be a dummy self-transition.

Similarly, straightforward algebra applied to the optimality equation \eqref{eq:opt-zero} at state $(0,1,1,0)$ (or equivalently, thinning at a higher rate $\mu_0 + 2\mu_1$) yields
\begin{align}
    & v(0,1,1,0) = \frac{h_0 + h_1}{\mu_0+2\mu_1} + \frac{\mu_1}{\mu_0+2\mu_1} v(0,1,0,0) \nonumber \\
    &\qquad + \frac{\mu_0}{\mu_0+2\mu_1}\min\{v(0,0,2,0), v(0,0,1,1)\}+ \frac{\mu_1}{\mu_0+2\mu_1}v(0,1,1,0). \label{eq:thinning_i110}
\end{align}
Substracting \eqref{eq:thinning_i110} from \eqref{eq:thinning_i020}, the left-hand side of \eqref{eq:onemore1_110} at $i=0$ becomes}
\begin{align*}
    \lefteqn{v(1,0,2,0) - v(0,1,1,0)}&\\
    &\quad = \frac{h_1}{\mu_0+2\mu_1} + \frac{\mu_1}{\mu_0+2\mu_1} \Big[v(0,1,1,0) - v(0,1,0,0)\Big]\\
    &\qquad + \frac{\mu_0}{\mu_0+2\mu_1}\Big[v(1,0,2,0) - \min\{v(0,0,2,0), v(0,0,1,1)\}\Big]\\
    &\qquad + \frac{\mu_1}{\mu_0+2\mu_1}\Big[v(0,1,1,0) - v(0,1,1,0)\Big]\\
    &\quad \geq \frac{h_1}{\mu_0+2\mu_1} + \frac{\mu_1}{\mu_0+2\mu_1} \Big[v(0,1,1,0) - v(0,1,0,0)\Big]\\
    &\qquad + \frac{\mu_0}{\mu_0+2\mu_1}\Big[v(1,0,2,0) - v(0,0,2,0)\Big],
\end{align*}where the inequality holds by selecting an upper bound $v(0,0,2,0)$ in the minimum, which corresponds to going to Station 1. Applying \eqref{eq:empty-queue-bound1} and \eqref{eq:mono} yields
\begin{align*}
    \lefteqn{v(1,0,2,0) - v(0,1,1,0)}&\\
    &\quad \geq \frac{h_1}{\mu_0+2\mu_1} + \frac{\mu_1}{\mu_0+2\mu_1} \frac{h_1}{\mu_1} + \frac{\mu_0}{\mu_0+2\mu_1}\left(\frac{h_0}{\mu_0}+\min\left\{\frac{h_1}{\mu_1}, \frac{h_2}{\mu_2}\right\}\right)\\
    &\quad \geq \textcolor{blue}{\frac{\mu_1}{\mu_0+2\mu_1} \frac{h_1}{\mu_1} + \frac{\mu_1}{\mu_0+2\mu_1} \frac{h_1}{\mu_1} + \frac{\mu_0}{\mu_0+2\mu_1} \min\left\{\frac{h_1}{\mu_1}, \frac{h_2}{\mu_2}\right\}}\\
    &\quad \geq \textcolor{blue}{\frac{\mu_1}{\mu_0+2\mu_1} \min\left\{\frac{h_1}{\mu_1}, \frac{h_2}{\mu_2}\right\} + \frac{\mu_1}{\mu_0+2\mu_1} \min\left\{\frac{h_1}{\mu_1}, \frac{h_2}{\mu_2}\right\} + \frac{\mu_0}{\mu_0+2\mu_1} \min\left\{\frac{h_1}{\mu_1}, \frac{h_2}{\mu_2}\right\}}\\
    &\quad = \min\left\{\frac{h_1}{\mu_1}, \frac{h_2}{\mu_2}\right\}.
\end{align*}
It follows similarly to prove \eqref{eq:onemore1_101} when $i=0$ by thinning the MDPs with a higher rate $\mu_0+\mu_1+\mu_2$. Finally for \eqref{eq:onemore1_200}, thinning the MDPs by $2\mu_0+\mu_1$ yields
\begin{align}
    \lefteqn{v(1,1,1,0) - v(0,2,0,0)}& \nonumber \\
    &\quad = \frac{h_1}{2\mu_0+\mu_1} + \frac{\mu_0}{2\mu_0+\mu_1}\Big[\min\{v(1,0,2,0), v(1,0,1,1)\} - \min\{v(0,1,1,0), v(0,1,0,1)\}\Big] \nonumber \\
    &\qquad + \frac{\mu_0}{2\mu_0+\mu_1} \Big[v(1,1,1,0) - \min\{v(0,1,1,0), v(0,1,0,1)\}\Big] \nonumber \\
    &\qquad + \frac{\mu_1}{2\mu_0+\mu_1} \Big[v(0,2,0,0) - v(0,2,0,0)\Big] \nonumber \\
    &\quad \geq \frac{h_1}{2\mu_0+\mu_1} + \frac{\mu_0}{2\mu_0+\mu_1}\Big[\min\{v(1,0,2,0), v(1,0,1,1)\} - \min\{v(0,1,1,0), v(0,1,0,1)\}\Big] \nonumber \\
    &\qquad + \frac{\mu_0}{2\mu_0+\mu_1} \Big[v(1,1,1,0) - v(0,1,1,0)\Big] \label{eq:onemore1_inte1},
\end{align}where the inequality replaces the third minimum with an upper bound $(0,1,1,0)$. For the difference $\min\{v(1,0,2,0), v(1,0,1,1)\} - \min\{v(0,1,1,0), v(0,1,0,1)\}$ appearing in the second term, it remains to distinguish two cases according to an optimal decision at state $(1,1,1,0)$, i.e., whether $v(1,0,2,0) \leq v(1,0,1,1)$. \textcolor{blue}{Let $a^*(1,1,1,0)$ denote an optimal decision at this state. We then upper bound the second minimum by selecting the value function corresponding to choosing action $a^*(1,1,1,0)$ at state $(0,2,0,0)$.} Algebraically, this is saying
\begin{align*}
    \lefteqn{\min\{v(1,0,2,0), v(1,0,1,1)\} - \min\{v(0,1,1,0), v(0,1,0,1)\}}&\\
    &\quad \geq \min\{v(1,0,2,0) - v(0,1,1,0), v(1,0,1,1) - v(0,1,0,1)\}\\
    &\quad \geq \min\left\{\min\left\{\frac{h_1}{\mu_1}, \frac{h_2}{\mu_2}\right\}, \min\left\{\frac{h_1}{\mu_1}, \frac{h_2}{\mu_2}\right\}\right\}\\
    &\quad = \min\left\{\frac{h_1}{\mu_1}, \frac{h_2}{\mu_2}\right\},
\end{align*}where the last inequality applies the proved results \eqref{eq:onemore1_110} and \eqref{eq:onemore1_101} at $i=0$. Using this lower bound and \eqref{eq:mono} in \eqref{eq:onemore1_inte1} yields
\begin{align*}
    \lefteqn{v(1,1,1,0) - v(0,2,0,0)}&\\
    &\quad \geq \frac{h_1}{2\mu_0+\mu_1} + \frac{\mu_0}{2\mu_0+\mu_1}\min\left\{\frac{h_1}{\mu_1}, \frac{h_2}{\mu_2}\right\} + \frac{\mu_0}{2\mu_0+\mu_1} \left(\frac{h_0}{\mu_0}+\min\left\{\frac{h_1}{\mu_1}, \frac{h_2}{\mu_2}\right\}\right)\\
     &\quad \geq \textcolor{blue}{\frac{\mu_1}{2\mu_0+\mu_1}\min\left\{\frac{h_1}{\mu_1}, \frac{h_2}{\mu_2}\right\} + \frac{\mu_0}{2\mu_0+\mu_1}\min\left\{\frac{h_1}{\mu_1}, \frac{h_2}{\mu_2}\right\} + \frac{\mu_0}{2\mu_0+\mu_1} \min\left\{\frac{h_1}{\mu_1}, \frac{h_2}{\mu_2}\right\}}\\
    &\quad = \min\left\{\frac{h_1}{\mu_1}, \frac{h_2}{\mu_2}\right\}.
\end{align*}
Assume all the inequalities \eqref{eq:onemore1_110} - \eqref{eq:onemore1_200} hold at $i-1$ and consider them at $i$. Thinning the MDPs by a higher rate $\mu_0+2\mu_1$ in \eqref{eq:onemore1_110} yields
\begin{align*}
    \lefteqn{v(i+1,0,2,0) - v(i,1,1,0)}&\\
    &\quad = \frac{h_1}{\mu_0+2\mu_1} + \frac{\mu_1}{\mu_0+2\mu_1} \Big[v(i,1,1,0) - v(i-1,2,0,0)\Big]\\
    &\qquad + \frac{\mu_0}{\mu_0+2\mu_1}\Big[v(i+1,0,2,0) - \min\{v(i,0,2,0), v(i,0,1,1)\}\Big]\\
    &\qquad + \frac{\mu_1}{\mu_0+2\mu_1}\Big[v(i,1,1,0) - v(i,1,1,0)\Big]\\
    &\quad \geq \frac{h_1}{\mu_0+2\mu_1} + \frac{\mu_1}{\mu_0+2\mu_1} \Big[v(i,1,1,0) - v(i-1,2,0,0)\Big]\\
    &\qquad + \frac{\mu_0}{\mu_0+2\mu_1}\Big[v(i+1,0,2,0) - v(i,0,2,0)\Big],
\end{align*}where the inequality holds by selecting an upper bound $v(i,0,2,0)$ in the minimum, which corresponds to going to Station 1. Applying the inductive hypothesis of \eqref{eq:onemore1_200} and the result in \eqref{eq:mono} yields
\begin{align*}
    \lefteqn{v(i+1,0,2,0) - v(i,1,1,0)}&\\
    &\quad \geq \frac{h_1}{\mu_0+2\mu_1} + \frac{\mu_1}{\mu_0+2\mu_1} \min\left\{\frac{h_1}{\mu_1}, \frac{h_2}{\mu_2}\right\} + \frac{\mu_0}{\mu_0+2\mu_1}\left(\frac{h_0}{\mu_0}+\min\left\{\frac{h_1}{\mu_1}, \frac{h_2}{\mu_2}\right\}\right)\\
    &\quad \geq \textcolor{blue}{\frac{\mu_1}{\mu_0+2\mu_1} \min\left\{\frac{h_1}{\mu_1}, \frac{h_2}{\mu_2}\right\} + \frac{\mu_1}{\mu_0+2\mu_1} \min\left\{\frac{h_1}{\mu_1}, \frac{h_2}{\mu_2}\right\} + \frac{\mu_0}{\mu_0+2\mu_1} \min\left\{\frac{h_1}{\mu_1}, \frac{h_2}{\mu_2}\right\}}\\
    &\quad = \min\left\{\frac{h_1}{\mu_1}, \frac{h_2}{\mu_2}\right\}.
\end{align*}
Again, one could notice that it follows analogously to prove \eqref{eq:onemore1_101} by thinning with a higher rate $\mu_0+\mu_1+\mu_2$. By thinning the MDPs with $2\mu_0+\mu_1$, it is also true that the proof for \eqref{eq:onemore1_200} is similar to that at $i=0$, using the established results \eqref{eq:onemore1_110} and \eqref{eq:onemore1_101} and the result \eqref{eq:mono}.

Finally, a comparable argument by induction on $i$ also works for \eqref{eq:onemore2}. The slight difference is choosing the upper bound corresponding to the decision of going to Station 2 instead in the minimum that occurs during the additional transition of initial service completion of the base state $(i,j,k,\ell)$. In addition, when proving for the base case at $i=0$, the result \eqref{eq:empty-queue-bound2} needs to be resorted to rather than \eqref{eq:empty-queue-bound1}. The detailed induction is also presented here for the reader's convenience. 

It is equivalent to showing the following three inequalities for \eqref{eq:onemore2}:
\begin{align}
    v(i+1,0,1,1) - v(i,1,1,0) &\geq \min\left\{\frac{h_1}{\mu_1}, \frac{h_2}{\mu_2}\right\}, \label{eq:onemore2_110}\\
    v(i+1,0,0,2) - v(i,1,0,1) &\geq \min\left\{\frac{h_1}{\mu_1}, \frac{h_2}{\mu_2}\right\},\label{eq:onemore2_101}\\
    v(i+1,1,0,1) - v(i,2,0,0) &\geq \min\left\{\frac{h_1}{\mu_1}, \frac{h_2}{\mu_2}\right\}\label{eq:onemore2_200}.
\end{align}
When $i=0$,
thinning the MDPs \textcolor{blue}{in \eqref{eq:onemore2_110}} by a higher rate $\mu_0+\mu_1+\mu_2$ yields
\begin{align*}
    \lefteqn{v(1,0,1,1) - v(0,1,1,0)}&\\
    &\quad = \frac{h_2}{\mu_0+\mu_1+\mu_2} + \frac{\mu_1}{\mu_0+\mu_1+\mu_2} \Big[v(0,1,0,1) - v(0,1,0,0)\Big]\\
    &\qquad + \frac{\mu_0}{\mu_0+\mu_1+\mu_2}\Big[v(1,0,1,1) - \min\{v(0,0,2,0), v(0,0,1,1)\}\Big]\\
    &\qquad + \frac{\mu_2}{\mu_0+\mu_1+\mu_2}\Big[v(0,1,1,0) - v(0,1,1,0)\Big]\\
    &\quad \geq \frac{h_2}{\mu_0+\mu_1+\mu_2} + \frac{\mu_1}{\mu_0+\mu_1+\mu_2} \Big[v(0,1,0,1) - v(0,1,0,0)\Big]\\
    &\qquad + \frac{\mu_0}{\mu_0+\mu_1+\mu_2}\Big[v(1,0,1,1) - v(0,0,1,1)\Big],
\end{align*}where the inequality holds by choosing the upper bound $v(0,0,1,1)$ corresponding to the decision of going to Station 2 at state $(0,1,1,0)$. The application of \eqref{eq:empty-queue-bound2} and \eqref{eq:mono} results in
\begin{align*}
    \lefteqn{v(1,0,1,1) - v(0,1,1,0)}&\\
    &\quad \geq \frac{h_2}{\mu_0+\mu_1+\mu_2} + \frac{\mu_1}{\mu_0+\mu_1+\mu_2} \frac{h_2}{\mu_2} + \frac{\mu_0}{\mu_0+\mu_1+\mu_2}\left(\frac{h_0}{\mu_0}+\min\left\{\frac{h_1}{\mu_1}, \frac{h_2}{\mu_2}\right\}\right)\\
    &\quad \geq \textcolor{blue}{\frac{\mu_2}{\mu_0+\mu_1+\mu_2}\min\left\{\frac{h_1}{\mu_1}, \frac{h_2}{\mu_2}\right\} + \frac{\mu_1}{\mu_0+\mu_1+\mu_2} \min\left\{\frac{h_1}{\mu_1}, \frac{h_2}{\mu_2}\right\} + \frac{\mu_0}{\mu_0+\mu_1+\mu_2}\min\left\{\frac{h_1}{\mu_1}, \frac{h_2}{\mu_2}\right\}}\\
    &\quad = \min\left\{\frac{h_1}{\mu_1}, \frac{h_2}{\mu_2}\right\}.
\end{align*}
The proof \textcolor{blue}{of \eqref{eq:onemore2_101}} at $i=0$ follows similarly by thinning the MDPs with rate $\mu_0+\mu_2$. Thinning the MDPs by $2\mu_0+\mu_2$ \textcolor{blue}{in \eqref{eq:onemore2_200}} results in
\begin{align}
    \lefteqn{v(1,1,0,1) - v(0,2,0,0)}& \nonumber \\
    &\quad = \frac{h_2}{2\mu_0+\mu_2} + \frac{\mu_0}{2\mu_0+\mu_2}\Big[\min\{v(1,0,1,1), v(1,0,0,2)\} - \min\{v(0,1,1,0), v(0,1,0,1)\}\Big] \nonumber \\
    &\qquad + \frac{\mu_0}{2\mu_0+\mu_2} \Big[v(1,1,0,1) - \min\{v(0,1,1,0), v(0,1,0,1)\}\Big] \nonumber \\
    &\qquad + \frac{\mu_2}{2\mu_0+\mu_2} \Big[v(0,2,0,0) - v(0,2,0,0)\Big] \nonumber \\
    &\quad \geq \frac{h_2}{2\mu_0+\mu_2} + \frac{\mu_0}{2\mu_0+\mu_2}\Big[\min\{v(1,0,1,1), v(1,0,0,2)\} - \min\{v(0,1,1,0), v(0,1,0,1)\}\Big] \nonumber \\
    &\qquad + \frac{\mu_0}{2\mu_0+\mu_2} \Big[v(1,1,0,1) - v(0,1,0,1)\Big]\label{eq:onemore2_inte1},
\end{align}where the inequality replaces the third minimum with an upper bound $(0,1,0,1)$, i.e., choosing to go to Station 2 at state $(0,2,0,0)$. Consider now the difference between the two minima $\min\{v(1,0,1,1), v(1,0,0,2)\} - \min\{v(0,1,1,0), v(0,1,0,1)\}$ in the second term in \eqref{eq:onemore2_inte1}. Letting the state $(0,2,0,0)$ follow the optimal action taken at state $(1,1,0,1)$ yields
\begin{align*}
    \lefteqn{\min\{v(1,0,1,1), v(1,0,0,2)\} - \min\{v(0,1,1,0), v(0,1,0,1)\}}&\\
    &\quad \geq \min\{v(1,0,1,1) - v(0,1,1,0), v(1,0,0,2) - v(0,1,0,1)\}\\
    &\quad \geq \min\left\{\frac{h_1}{\mu_1}, \frac{h_2}{\mu_2}\right\} 
\end{align*}by \eqref{eq:onemore2_110} and \eqref{eq:onemore2_101} at $i=0$. Applying this bound and \eqref{eq:mono} in \eqref{eq:onemore2_inte1} yields
\begin{align*}
    \lefteqn{v(1,1,0,1) - v(0,2,0,0)}&\\
    &\quad \geq \frac{h_2}{2\mu_0+\mu_2} + \frac{\mu_0}{2\mu_0+\mu_2}\min\left\{\frac{h_1}{\mu_1}, \frac{h_2}{\mu_2}\right\} + \frac{\mu_0}{2\mu_0+\mu_2} \left(\frac{h_0}{\mu_0}+\min\left\{\frac{h_1}{\mu_1}, \frac{h_2}{\mu_2}\right\}\right)\\
    &\quad \geq \textcolor{blue}{\frac{\mu_2}{2\mu_0+\mu_2}\min\left\{\frac{h_1}{\mu_1}, \frac{h_2}{\mu_2}\right\} + \frac{\mu_0}{2\mu_0+\mu_2}\min\left\{\frac{h_1}{\mu_1}, \frac{h_2}{\mu_2}\right\} + \frac{\mu_0}{2\mu_0+\mu_2} \min\left\{\frac{h_1}{\mu_1}, \frac{h_2}{\mu_2}\right\}}\\
    &\quad = \min\left\{\frac{h_1}{\mu_1}, \frac{h_2}{\mu_2}\right\}.
\end{align*}
Assume the results \eqref{eq:onemore2_110} - \eqref{eq:onemore2_200} hold at $i-1$ to prove them at $i$. After thinning the MDPs by the rate $\mu_0+\mu_1+\mu_2$, the difference of the value functions \textcolor{blue}{in \eqref{eq:onemore2_110}} becomes
\begin{align*}
    \lefteqn{v(i+1,0,1,1) - v(i,1,1,0)}&\\
    &\quad = \frac{h_2}{\mu_0+\mu_1+\mu_2} + \frac{\mu_1}{\mu_0+\mu_1+\mu_2} \Big[v(i,1,0,1) - v(i-1,2,0,0)\Big]\\
    &\qquad + \frac{\mu_0}{\mu_0+\mu_1+\mu_2}\Big[v(i+1,0,1,1) - \min\{v(i,0,2,0), v(i,0,1,1)\}\Big]\\
    &\qquad + \frac{\mu_2}{\mu_0+\mu_1+\mu_2}\Big[v(i,1,1,0) - v(i,1,1,0)\Big]\\
    &\quad \geq \frac{h_2}{\mu_0+\mu_1+\mu_2} + \frac{\mu_1}{\mu_0+\mu_1+\mu_2} \Big[v(i,1,0,1) - v(i-1,2,0,0)\Big]\\
    &\qquad + \frac{\mu_0}{\mu_0+\mu_1+\mu_2}\Big[v(i+1,0,1,1) - v(i,0,1,1)\Big],
\end{align*}where the inequality is by replacing the minimum with the upper bound $v(i,0,1,1)$, indicating going to Station 2. Applying the inductive hypothesis of \eqref{eq:onemore2_200} and the result in \eqref{eq:mono} gives
\begin{align*}
    \lefteqn{v(i+1,0,1,1) - v(i,1,1,0)}&\\
    &\quad \geq \frac{h_2}{\mu_0+\mu_1+\mu_2} + \frac{\mu_1}{\mu_0+\mu_1+\mu_2} \min\left\{\frac{h_1}{\mu_1}, \frac{h_2}{\mu_2}\right\} + \frac{\mu_0}{\mu_0+\mu_1+\mu_2}\left(\frac{h_0}{\mu_0}+\min\left\{\frac{h_1}{\mu_1}, \frac{h_2}{\mu_2}\right\}\right)\\
    &\quad \geq \textcolor{blue}{\frac{\mu_2}{\mu_0+\mu_1+\mu_2}\min\left\{\frac{h_1}{\mu_1}, \frac{h_2}{\mu_2}\right\} + \frac{\mu_1}{\mu_0+\mu_1+\mu_2} \min\left\{\frac{h_1}{\mu_1}, \frac{h_2}{\mu_2}\right\} + \frac{\mu_0}{\mu_0+\mu_1+\mu_2}\min\left\{\frac{h_1}{\mu_1}, \frac{h_2}{\mu_2}\right\}}\\
    &\quad = \min\left\{\frac{h_1}{\mu_1}, \frac{h_2}{\mu_2}\right\}.
\end{align*}
The proof of \eqref{eq:onemore2_101} is again comparable by thinning with the rate $\mu_0+\mu_2$. In addition, after thinning the MDPs by $2\mu_0+\mu_2$, the proof for \eqref{eq:onemore2_200} follows similarly to that at $i=0$, by applying the established result \eqref{eq:onemore2_110} and \eqref{eq:onemore2_101}, and the result \eqref{eq:mono}. This finally completes the proof of Statement \ref{state2:general_rates}.
\end{proof}
\begin{proof} [Proof of Statement \ref{state3:general_rates}] In order to prove \eqref{eq:decision_bd1}, an intermediate result is also required:
\begin{align}
    v(i+1,j-1,k,\ell+1) - v(i,j,k,\ell) \geq \frac{h_2}{\mu_2}.\label{eq:decision_bd1_inte}
\end{align}
Notice that if $\frac{h_1}{\mu_1} \geq \frac{h_2}{\mu_2}$ holds, the intermediate result follows immediately by \eqref{eq:onemore2}, and hence without loss of generality we assume $\frac{h_1}{\mu_1} < \frac{h_2}{\mu_2}$ in proving \eqref{eq:decision_bd1_inte} only, in addition to $\mu_1\geq\mu_2$. Under this circumstance, \textcolor{blue}{the result \eqref{eq:mono} becomes
\begin{align}
    v(i+1,j,k,\ell) - v(i,j,k,\ell) \geq \frac{h_1}{\mu_1}, \quad \forall (i,j,k,\ell) \in \X.\label{eq:decision_bd1_inte1}
\end{align}
}
For each set of values that $(j,k,\ell) \in \X_{\widehat{D}}$ can take, we need to show the following for \eqref{eq:decision_bd1} and \eqref{eq:decision_bd1_inte}:
\begin{enumerate}[label= \textbf{Case} \arabic*:, leftmargin=3\parindent]
    \item For $(j,k,\ell) = (1,1,0)$,
    \begin{align}
        D(i,1,1,0) = v(i,0,2,0) - v(i,0,1,1) &\leq \frac{h_1}{\mu_1} - \frac{h_2}{\mu_2}, \label{eq:decision_bd1_110}\\
        v(i+1,0,1,1) - v(i,1,1,0) &\geq \frac{h_2}{\mu_2}. \label{eq:decision_bd1_inte_110}
    \end{align}
    \item For $(j,k,\ell) = (1,0,1)$,
    \begin{align}
        D(i,1,0,1) = v(i,0,1,1) - v(i,0,0,2) &\leq \frac{h_1}{\mu_1} - \frac{h_2}{\mu_2}, \label{eq:decision_bd1_101}\\
        v(i+1,0,0,2) - v(i,1,0,1) &\geq  \frac{h_2}{\mu_2}. \label{eq:decision_bd1_inte_101}
    \end{align}
    \item
    For $(j,k,\ell) = (2,0,0)$,
    \begin{align}
        D(i,2,0,0) = v(i,1,1,0) - v(i,1,0,1) &\leq \frac{h_1}{\mu_1} - \frac{h_2}{\mu_2}, \label{eq:decision_bd1_200}\\
        v(i+1,1,0,1) - v(i,2,0,0) &\geq \frac{h_2}{\mu_2}.\label{eq:decision_bd1_inte_200}
    \end{align}
\end{enumerate}
We prove the four inequalities together by induction on $i$ under the assumption that $\mu_1 \geq \mu_2$. When $i=0$, consider first $(j,k,\ell) = (1,1,0)$. Thinning the MDPs by a higher rate $2\mu_1$ in \eqref{eq:decision_bd1_110} results in
\begin{align*}
    \lefteqn{D(0,1,1,0) = v(0,0,2,0) - v(0,0,1,1)}&\\
    &\quad =\frac{h_1-h_2}{2\mu_1} + \frac{\mu_1}{2\mu_1} \Big[v(0,0,1,0) - v(0,0,0,1)\Big]\\
    &\qquad + \frac{\mu_2}{2\mu_1} \Big[v(0,0,1,0) - v(0,0,1,0)\Big] + \frac{\mu_1-\mu_2}{2\mu_1} \Big[v(0,0,1,0) - v(0,0,1,1)\Big]\\
    &\quad \leq\frac{h_1-h_2}{2\mu_1} + \frac{1}{2} \left(\frac{h_1}{\mu_1} - \frac{h_2}{\mu_2}\right) + \frac{\mu_1-\mu_2}{2\mu_1} \left(-\frac{h_2}{\mu_2}\right)\\
    &\quad = \frac{h_1}{\mu_1} - \frac{h_2}{\mu_2},
\end{align*}where the inequality is by recalling \eqref{eq:opt_single} and \eqref{eq:empty-queue-bound2}.
As for the intermediate result \eqref{eq:decision_bd1_inte} at $i=0$, recall that thinning the MDPs in \eqref{eq:decision_bd1_inte_110} by a higher rate $\mu_0+\mu_1+\mu_2$ yields
\begin{align}
    \lefteqn{v(1,0,1,1) - v(0,1,1,0)}& \nonumber \\
    &\quad = \frac{h_2}{\mu_0+\mu_1+\mu_2} + \frac{\mu_1}{\mu_0+\mu_1+\mu_2} \Big[v(0,1,0,1) - v(0,1,0,0)\Big] \nonumber \\
    &\qquad + \frac{\mu_0}{\mu_0+\mu_1+\mu_2}\Big[v(1,0,1,1) - \min\{v(0,0,2,0), v(0,0,1,1)\}\Big] \nonumber \\
    &\qquad + \frac{\mu_2}{\mu_0+\mu_1+\mu_2}\Big[v(0,1,1,0) - v(0,1,1,0)\Big]. \label{eq:decision_bd1_inte2}
\end{align}
Consider the difference $v(1,0,1,1) - \min\{v(0,0,2,0), v(0,0,1,1)\}$ in the third term of \eqref{eq:decision_bd1_inte2}. Choosing the upper bound $v(0,0,2,0)$ in the minimum (see the proof of yields \eqref{eq:onemore2_110} for comparison) 
\begin{align}
    \lefteqn{v(1,0,1,1) - \min\{v(0,0,2,0), v(0,0,1,1)\}}& \nonumber \\
    &\quad \geq v(1,0,1,1) - v(0,0,2,0) \nonumber \\
    &\quad = \Big[v(1,0,1,1) - v(0,0,1,1)\Big] - \Big[v(0,0,2,0) - v(0,0,1,1)\Big] \nonumber \\
    &\quad \geq \frac{h_1}{\mu_1} - \left(\frac{h_1}{\mu_1} - \frac{h_2}{\mu_2}\right) \nonumber \\
    &\quad = \frac{h_2}{\mu_2} \label{eq:decision_bd1_inte3}
\end{align}by \eqref{eq:decision_bd1_inte1} and \eqref{eq:decision_bd1_110} at $i=0$.
Using the bound in \eqref{eq:decision_bd1_inte3} in \eqref{eq:decision_bd1_inte2} and applying \eqref{eq:empty-queue-bound2} yield
\begin{align*}
    \lefteqn{v(1,0,1,1) - v(0,1,1,0)}& \nonumber \\
    &\quad \geq \frac{h_2}{\mu_0+\mu_1+\mu_2} + \frac{\mu_1}{\mu_0+\mu_1+\mu_2} \frac{h_2}{\mu_2} + \frac{\mu_0}{\mu_0+\mu_1+\mu_2}\frac{h_2}{\mu_2}\\
    &\quad = \frac{h_2}{\mu_2}.
\end{align*}
The case that $(j,k,\ell) = (1,0,1)$ at $i=0$ for both \eqref{eq:decision_bd1_101} and \eqref{eq:decision_bd1_inte_101} follows analagously by thinning with rates $\mu_1+\mu_2$ and $\mu_0+\mu_2$, respectively. 

It remains to examine when $(j,k,\ell) = (2,0,0)$. Thinning by the higher rate $\mu_0+\mu_1$ in \eqref{eq:decision_bd1_200} gives
\begin{align*}
    \lefteqn{D(0,2,0,0) = v(0,1,1,0) - v(0,1,0,1)}&\\
    &\quad =\frac{h_1-h_2}{\mu_0+\mu_1} + \frac{\mu_0}{\mu_0+\mu_1} \Big[\min\{v(0,0,2,0), v(0,0,1,1)\} - \min\{v(0,0,1,1), v(0,0,0,2)\}\Big]\\
    &\qquad + \frac{\mu_2}{\mu_0+\mu_1} \Big[v(0,1,0,0) - v(0,1,0,0)\Big] + \frac{\mu_1-\mu_2}{\mu_0+\mu_1} \Big[v(0,1,0,0) - v(0,1,0,1)\Big].
\end{align*}
There are two cases to discuss based on whether $v(0,0,1,1) \leq v(0,0,0,2)$, and we use the same trick as before that chooses the potentially sub-optimal action at state $(0,1,1,0)$, which is identical to the \textcolor{blue}{decision of an optimal policy at} state $(0,1,0,1)$. See details below.
\begin{enumerate}[label= \textbf{Case} \arabic*:, leftmargin=3\parindent]
    \item Suppose $v(0,0,1,1) \leq v(0,0,0,2)$, replacing the first minimum with an upper bound $v(0,0,2,0)$ yields
    \begin{align*}
        \lefteqn{v(0,1,1,0) - v(0,1,0,1)}&\\
        &\quad \leq\frac{h_1-h_2}{\mu_0+\mu_1} + \frac{\mu_0}{\mu_0+\mu_1} \Big[v(0,0,2,0) - v(0,0,1,1)\Big]\\
        &\qquad + \frac{\mu_1-\mu_2}{\mu_0+\mu_1} \Big[v(0,1,0,0) - v(0,1,0,1)\Big]\\
        &\quad \leq\frac{h_1-h_2}{\mu_0+\mu_1} + \frac{\mu_0}{\mu_0+\mu_1} \left(\frac{h_1}{\mu_1} - \frac{h_2}{\mu_2}\right) + \frac{\mu_1-\mu_2}{\mu_0+\mu_1}\left(-\frac{h_2}{\mu_2}\right)\\
        &\quad = \frac{h_1}{\mu_1} - \frac{h_2}{\mu_2},
    \end{align*}where the last inequality is by \eqref{eq:decision_bd1_110} at $i=0$ and by \eqref{eq:empty-queue-bound2}.
    \item Suppose $v(0,0,1,1) > v(0,0,0,2)$, replacing the first minimum with an upper bound $v(0,0,1,1)$ yields
    \begin{align*}
        \lefteqn{v(0,1,1,0) - v(0,1,0,1)}&\\
        &\quad \leq\frac{h_1-h_2}{\mu_0+\mu_1} + \frac{\mu_0}{\mu_0+\mu_1} \Big[v(0,0,1,1) - v(0,0,0,2)\Big] \\
        & \qquad + \frac{\mu_1-\mu_2}{\mu_0+\mu_1} \Big[v(0,1,0,0) - v(0,1,0,1)\Big]\\
        &\quad \leq\frac{h_1-h_2}{\mu_0+\mu_1} + \frac{\mu_0}{\mu_0+\mu_1} \left(\frac{h_1}{\mu_1} - \frac{h_2}{\mu_2}\right) + \frac{\mu_1-\mu_2}{\mu_0+\mu_1}\left(-\frac{h_2}{\mu_2}\right)\\
        &\quad = \frac{h_1}{\mu_1} - \frac{h_2}{\mu_2},
    \end{align*}where the last inequality is by \eqref{eq:decision_bd1_101} at $i=0$ and by \eqref{eq:empty-queue-bound2}.   
\end{enumerate}

For \eqref{eq:decision_bd1_inte_200} at $i=0$, thinning the MDPs by the rate $2\mu_0+\mu_2$ results in  
    \begin{align}
    \lefteqn{v(1,1,0,1) - v(0,2,0,0)}& \nonumber \\
    &\quad = \frac{h_2}{2\mu_0+\mu_2} + \frac{\mu_0}{2\mu_0+\mu_2}\Big[\min\{v(1,0,1,1), v(1,0,0,2)\} - \min\{v(0,1,1,0), v(0,1,0,1)\}\Big] \nonumber \\
    &\qquad + \frac{\mu_0}{2\mu_0+\mu_2} \Big[v(1,1,0,1) - \min\{v(0,1,1,0), v(0,1,0,1)\}\Big] \nonumber \\
    &\qquad + \frac{\mu_2}{2\mu_0+\mu_2} \Big[v(0,2,0,0) - v(0,2,0,0)\Big]\label{eq:decision_bd1_inte4}
\end{align}
Consider the difference $\min\{v(1,0,1,1), v(1,0,0,2)\} - \min\{v(0,1,1,0), v(0,1,0,1)\}$ in the second term in \eqref{eq:decision_bd1_inte4}. We use a similar approach that makes the state $(0,2,0,0)$ choose the same but potentially sub-optimal action as an optimal decision at state $(1,1,0,1)$ to find
\begin{align}
    \lefteqn{\min\{v(1,0,1,1), v(1,0,0,2)\} - \min\{v(0,1,1,0), v(0,1,0,1)\}}& \nonumber \\
    &\quad \geq \min\{v(1,0,1,1)- v(0,1,1,0),v(1,0,0,2) - v(0,1,0,1)\} \nonumber \\
    &\quad \geq \min\left\{\frac{h_2}{\mu_2},\frac{h_2}{\mu_2}\right\} = \frac{h_2}{\mu_2} \label{eq:decision_bd1_inte5}
\end{align}
by what has been proved for \eqref{eq:decision_bd1_inte_110} and \textcolor{blue}{\eqref{eq:decision_bd1_inte_101} at} $i=0$. 
Similarly as \eqref{eq:decision_bd1_inte3}, one can show that the difference $v(1,1,0,1) - \min\{v(0,1,1,0), v(0,1,0,1)\}$ in the third term in \eqref{eq:decision_bd1_inte4} is lower bounded by $\frac{h_2}{\mu_2}$ as well, by using \eqref{eq:decision_bd1_inte1} and \eqref{eq:decision_bd1_200} at $i=0$. Applying this bound and that of \eqref{eq:decision_bd1_inte5} back to \eqref{eq:decision_bd1_inte4} yields
\begin{align*}
    v(1,1,0,1) - v(0,2,0,0)\geq \frac{h_2}{2\mu_0+\mu_2} + \frac{\mu_0}{2\mu_0+\mu_2}\frac{h_2}{\mu_2} + \frac{\mu_0}{2\mu_0+\mu_2} \frac{h_2}{\mu_2}= \frac{h_2}{\mu_2}.
\end{align*}
Assume all the results hold at $i-1$ and consider them at $i\geq1$. When $(j,k,\ell) = (1,1,0)$, again thinning the MDPs by the higher rate $2\mu_1$ in \eqref{eq:decision_bd1_110} results in
\begin{align*}
    \lefteqn{D(i,1,1,0) = v(i,0,2,0) - v(i,0,1,1)}&\\
    &\quad =\frac{h_1-h_2}{2\mu_1} + \frac{\mu_1}{2\mu_1} \Big[v(i-1,1,1,0) - v(i-1,1,0,1)\Big]\\
    &\qquad + \frac{\mu_2}{2\mu_1} \Big[v(i-1,1,1,0) - v(i-1,1,1,0)\Big] + \frac{\mu_1-\mu_2}{2\mu_1} \Big[v(i-1,1,1,0) - v(i,0,1,1)\Big]\\
    &\quad \leq\frac{h_1-h_2}{2\mu_1} + \frac{1}{2} \left(\frac{h_1}{\mu_1} - \frac{h_2}{\mu_2}\right) + \frac{\mu_1-\mu_2}{2\mu_1} \left(-\frac{h_2}{\mu_2}\right)\\
    &\quad = \frac{h_1}{\mu_1} - \frac{h_2}{\mu_2},
\end{align*}where the inequality is by the inductive hypotheses of \eqref{eq:decision_bd1_200} and \eqref{eq:decision_bd1_inte_110}. Thinning the MDPs in \eqref{eq:decision_bd1_inte_110} by a higher rate $\mu_0+\mu_1+\mu_2$ yields
\begin{align*}
    \lefteqn{v(i+1,0,1,1) - v(i,1,1,0)}&\\
    &\quad = \frac{h_2}{\mu_0+\mu_1+\mu_2} + \frac{\mu_1}{\mu_0+\mu_1+\mu_2} \Big[v(i,1,0,1) - v(i-1,2,0,0)\Big]\\
    &\qquad + \frac{\mu_0}{\mu_0+\mu_1+\mu_2}\Big[v(i+1,0,1,1) - \min\{v(i,0,2,0), v(i,0,1,1)\}\Big]\\
    &\qquad + \frac{\mu_2}{\mu_0+\mu_1+\mu_2}\Big[v(i,1,1,0) - v(i,1,1,0)\Big].
\end{align*}
Similarly \textcolor{blue}{to} proving \eqref{eq:decision_bd1_inte_110} at $i=0$, the difference $v(i+1,0,1,1) - \min\{v(i,0,2,0), v(i,0,1,1)\}$ in the third term is lower bounded by $\frac{h_2}{\mu_2}$ by using \eqref{eq:decision_bd1_inte1} and \eqref{eq:decision_bd1_110}. Using this and the inductive hypothesis of \eqref{eq:decision_bd1_inte_200} yields the desired result.
A comparable argument applies for both \eqref{eq:decision_bd1_101} and \eqref{eq:decision_bd1_inte_101} when $(j,k,\ell) = (1,0,1)$, by thinning with rates $\mu_1+\mu_2$ and $\mu_0+\mu_2$, respectively. 

Finally consider the case when $(j,k,\ell) = (2,0,0)$. Thinning the MDPs in \eqref{eq:decision_bd1_200} with rate $\mu_0+\mu_1$ yields
\begin{align*}
    \lefteqn{D(i,2,0,0) = v(i,1,1,0) - v(i,1,0,1)}&\\
    &\quad =\frac{h_1-h_2}{\mu_0+\mu_1} + \frac{\mu_0}{\mu_0+\mu_1} \Big[\min\{v(i,0,2,0), v(i,0,1,1)\} - \min\{v(i,0,1,1), v(i,0,0,2)\}\Big]\\
    &\qquad + \frac{\mu_2}{\mu_0+\mu_1} \Big[v(i-1,2,0,0) - v(i-1,2,0,0)\Big] + \frac{\mu_1-\mu_2}{\mu_0+\mu_1} \Big[v(i-1,2,0,0) - v(i,1,0,1)\Big].
\end{align*}
It can be verified that $\min\{v(i,0,2,0), v(i,0,1,1)\} - \min\{v(i,0,1,1), v(i,0,0,2)\} \leq \frac{h_1}{\mu_1} - \frac{h_2}{\mu_2}$ for the difference in the second term by the same reasoning as we have done for $i=0$ using \eqref{eq:decision_bd1_110} and \eqref{eq:decision_bd1_101}. Using this bound and the inductive hypothesis of \eqref{eq:decision_bd1_inte_200} proves the desired result.
Again, thinning the MDPs in \eqref{eq:decision_bd1_inte_200} with a higher rate $2\mu_0+\mu_2$, we obatin
    \begin{align}
    \lefteqn{v(i+1,1,0,1) - v(i,2,0,0)}& \nonumber \\
    &\quad = \frac{h_2}{2\mu_0+\mu_2} + \frac{\mu_0}{2\mu_0+\mu_2}\Big[\min\{v(i+1,0,1,1), v(i+1,0,0,2)\} - \min\{v(i,1,1,0), v(i,1,0,1)\}\Big] \nonumber \\
    &\qquad + \frac{\mu_0}{2\mu_0+\mu_2} \Big[v(i+1,1,0,1) - \min\{v(i,1,1,0), v(i,1,0,1)\}\Big] \nonumber \\
    &\qquad + \frac{\mu_2}{2\mu_0+\mu_2} \Big[v(i,2,0,0) - v(i,2,0,0)\Big] \label{eq:decision_bd1_inte6}.
\end{align}
As we proved at $i=0$, a similar logic applies to find that the difference $\min\{v(i+1,0,1,1), v(i+1,0,0,2)\} - \min\{v(i,1,1,0), v(i,1,0,1)\}$ in the second term is lower bounded by $\frac{h_2}{\mu_2}$ by resorting to the established results \eqref{eq:decision_bd1_inte_110} and \eqref{eq:decision_bd1_inte_101}. 
\textcolor{blue}{Similar to} \eqref{eq:decision_bd1_inte3}, the difference $v(i+1,1,0,1) - \min\{v(i,1,1,0), v(i,1,0,1)\}$ in the third term in \eqref{eq:decision_bd1_inte4} is lower bounded by $\frac{h_2}{\mu_2}$ as well, using the verified results \eqref{eq:decision_bd1_inte1} and \eqref{eq:decision_bd1_200}. 
Applying these bounds back in \eqref{eq:decision_bd1_inte6} proves \eqref{eq:decision_bd1_inte_200}.
\end{proof}

\begin{proof}[Proof of Statement \ref{state4:general_rates}]
Notice that $(i,j,k,\ell) \in \X_{\widehat{D}}$ contains three cases $(1,1,0)$, $(1,0,1)$ and $(2,0,0)$ for $(j,k,\ell)$, and thus \eqref{eq:decision_bd2} in Statement \ref{state4:general_rates} is implied by the following 
\begin{align} 
    D(i,1,1,0) = v(i,0,2,0) - v(i,0,1,1) &\leq \frac{h_1}{\mu_1} - \frac{h_2}{\mu_2} +\frac{ih_0}{2}\left(\frac{1}{\mu_1} - \frac{1}{\mu_2}\right), \label{eq:decision_bd2_110} \\
    D(i,1,0,1) = v(i,0,1,1) - v(i,0,0,2) &\leq \frac{h_1}{\mu_1} - \frac{2h_2}{\mu_2} + \frac{ih_0}{2}\left(\frac{1}{\mu_1} - \frac{2}{\mu_2}\right), \label{eq:decision_bd2_101}\\
    D(i,2,0,0) = v(i,1,1,0) - v(i,1,0,1) &\leq \frac{h_1}{\mu_1} - \frac{h_2}{\mu_2} +\frac{(i+1)h_0}{2}\left(\frac{1}{\mu_1} - \frac{1}{\mu_2}\right).\label{eq:decision_bd2_200}
\end{align}
\textcolor{blue}{We prove inequalities \eqref{eq:decision_bd2_110}--\eqref{eq:decision_bd2_200} together by induction on $i$}. 
When $i = 0$, by recalling the (in)equalities of values that have been computed in the proof of Statement \ref{state1:general_rates}, the differences in \eqref{eq:decision_bd2_110} - \eqref{eq:decision_bd2_200} respectively become
\begin{align*}
    D(0,1,1,0) = v(0,0,2,0) - v(0,0,1,1) &= \frac{h_1}{\mu_1} - \frac{h_2}{\mu_2},\\
    D(0,1,0,1) = v(0,0,1,1) - v(0,0,0,2) &= \frac{h_1}{\mu_1} - \frac{2h_2}{\mu_2},\\
    D(0,2,0,0) = v(0, 1, 1, 0) - v(0, 1, 0, 1)
    &\leq \frac{h_1}{\mu_1} - \frac{h_2}{\mu_2} \leq \frac{h_1}{\mu_1} - \frac{h_2}{\mu_2} + \frac{h_0}{2}\left(\frac{1}{\mu_1} - \frac{1}{\mu_2}\right).
\end{align*}
% By the optimality equations and the result \eqref{eq:always_action2_i200} in Proposition \ref{prop:always_action2} in the main text, the difference in \eqref{eq:decision_bd2_inte} is
% \begin{align*}
%     \lefteqn{v(1,0,1,1) - v(0,2,0,0)}&\\
%     &\quad = \left(\frac{h_0+h_1+h_2}{\mu_1+\mu_2} + \frac{\mu_1}{\mu_1+\mu_2}v(0,1,0,1)+ \frac{\mu_2}{\mu_1+\mu_2}v(0,1,1,0)\right)\\
%     &\qquad -\left(\frac{2h_0}{2\mu_0} + \min\{v(0,1,1,0),v(0,1,0,1)\}\right)\\
%     &\quad = \left(\frac{h_0+h_1+h_2}{\mu_1+\mu_2} + \frac{\mu_1}{\mu_1+\mu_2}v(0,1,0,1)+ \frac{\mu_2}{\mu_1+\mu_2}v(0,1,1,0)\right)\\
%     &\qquad -\left(\frac{h_0}{\mu_0} + v(0,1,0,1)\right)\\
%     &\quad = \frac{\mu_0-(\mu_1+\mu_2)}{\mu_0(\mu_1+\mu_2)}h_0+ \frac{h_1+h_2}{\mu_1+\mu_2} + \frac{\mu_2}{\mu_1+\mu_2}\Big[v(0,1,1,0)-v(0,1,0,1)\Big]\\
%     &\quad \leq \frac{\mu_0-(\mu_1+\mu_2)}{\mu_0(\mu_1+\mu_2)}h_0+ \frac{h_1+h_2}{\mu_1+\mu_2} + \frac{\mu_2}{\mu_1+\mu_2}\left(\frac{h_1}{\mu_1} - \frac{h_2}{\mu_2} + \frac{h_0}{2}\left(\frac{1}{\mu_1} - \frac{1}{\mu_2}\right)\right)\\
%     &\quad = \frac{h_1}{\mu_1} +h_0\left(\frac{1}{2\mu_1} - \frac{1}{\mu_0}\right) \leq \frac{h_1}{\mu_1} +\frac{2h_0}{2}\left(\frac{1}{\mu_1} - \frac{1}{\mu_0}\right),
% \end{align*}
% where the inequality uses the base case of \eqref{eq:decision_bd2_200} at $i=0$.

Now for $i \geq 1$, to establish the results at $i$, suppose all the inequalities hold at $i-1$. Consider the decision state $(i,1,1,0)$ first. 
\textcolor{blue}{A simple algebraic manipulation of the optimality equation \eqref{eq:opt1} for states $(i, 0, 2, 0)$ and $(i, 0, 1, 1)$ yields}
\begin{align}
    \lefteqn{D(i,1,1,0) = v(i, 0, 2, 0) - v(i, 0, 1, 1)}& \nonumber \\
    &\quad = \left(\frac{i h_0+2h_1}{2\mu_1} + v(i-1, 1, 1, 0)\right) \nonumber \\
    &\qquad - \left(\frac{i h_0+h_1+ h_2}{\mu_1+\mu_2} + \frac{\mu_1}{\mu_1+\mu_2} v(i-1, 1, 0, 1) + \frac{\mu_2}{ \mu_1+\mu_2} v(i-1, 1, 1, 0)\right) \nonumber \\
    &\quad = \frac{i (\mu_2-\mu_1) h_0}{2\mu_1(\mu_1+\mu_2)} + \frac{\mu_2h_1 - \mu_1h_2}{\mu_1(\mu_1+\mu_2)} + \frac{\mu_1}{ \mu_1+\mu_2} \Big[v(i-1, 1, 1, 0) - v(i-1, 1, 0, 1)\Big]. \label{eq:i110-diff}
\end{align} 
By the inductive hypothesis of \eqref{eq:decision_bd2_200}, \textcolor{black}{equation \eqref{eq:i110-diff} admits the desired upper bound:}
\begin{align*}
    \lefteqn{D(i,1,1,0) = v(i, 0, 2, 0) - v(i, 0, 1, 1)}&\\
    &\quad \leq \frac{i (\mu_2-\mu_1) h_0}{2\mu_1(\mu_1+\mu_2)} + \frac{\mu_2}{\mu_1+\mu_2}\left(\frac{h_1}{\mu_1} - \frac{h_2}{\mu_2}\right)+ \frac{\mu_1}{ \mu_1+\mu_2} \left(\frac{h_1}{\mu_1} - \frac{h_2}{\mu_2} + \frac{ih_0}{2}\left(\frac{1}{\mu_1} - \frac{1}{\mu_2}\right)\right)\\
    &\quad = \frac{h_1}{\mu_1} - \frac{h_2}{\mu_2} + \frac{ih_0}{2}\left(\frac{1}{\mu_1} - \frac{1}{\mu_2}\right).
\end{align*}
Similarly, \textcolor{blue}{applying the optimality equation \eqref{eq:opt1} to the states $(i, 0, 1, 1)$ and $(i, 0, 0, 2)$ yields}
\begin{align}
    \lefteqn{D(i,1,0,1) = v(i,0,1,1) - v(i,0,0,2)}& \nonumber \\
    &\quad = \left(\frac{ih_0+h_1+h_2}{\mu_1+\mu_2} + \frac{\mu_1}{\mu_1+\mu_2}v(i-1,1,0,1)+\frac{\mu_2}{\mu_1+\mu_2}v(i-1,1,1,0)\right) \nonumber \\
    &\qquad - \left(\frac{ih_0+2h_2}{\mu_2} + v(i-1,1,0,1)\right) \nonumber \\
    &\quad = \frac{-i\mu_1h_0}{\mu_2(\mu_1+\mu_2)} + \frac{h_1+h_2}{\mu_1+\mu_2} - \frac{2h_2}{\mu_2}  + \frac{\mu_2}{\mu_1+\mu_2}\Big[v(i-1,1,1,0) - v(i-1,1,0,1)\Big]. \label{eq:i101-diff}
\end{align} 
\textcolor{black}{By the inductive assumption of \eqref{eq:decision_bd2_200}, equation \eqref{eq:i101-diff} yields}
\begin{align*}
    \lefteqn{D(i,1,0,1) = v(i,0,1,1) - v(i,0,0,2)}&\\
    &\quad \leq \frac{-i\mu_1h_0}{\mu_2(\mu_1+\mu_2)} + \frac{h_1+h_2}{\mu_1+\mu_2} - \frac{2h_2}{\mu_2} + \frac{\mu_2}{\mu_1+\mu_2}\left(\frac{h_1}{\mu_1} - \frac{h_2}{\mu_2} + \frac{ih_0}{2}\left(\frac{1}{\mu_1} -\frac{1}{\mu_2}\right)\right)\\
    &\quad = \frac{h_1}{\mu_1} - \frac{2h_2}{\mu_2} + \frac{ih_0}{2}\left(\frac{1}{\mu_1} -\frac{2}{\mu_2}\right).
\end{align*}
Now consider the decision state $(i,2,0,0)$. Recall that $\mu_1 \leq \mu_2$.
\begin{align}
    \lefteqn{D(i,2,0,0) = v(i, 1, 1, 0) - v(i, 1, 0, 1)}& \nonumber \\
    &\quad = \left(\frac{(i+1) h_0+h_1}{\mu_0+\mu_1} + \frac{\mu_0}{\mu_0+\mu_1}\min\{v(i,0,2,0),v(i,0,1,1)\} + \frac{\mu_1}{ \mu_0+\mu_1} v(i-1, 2, 0, 0)\right) \nonumber \\
    &\qquad - \left(\frac{(i+1) h_0+ h_2}{\mu_0+\mu_2} + \frac{\mu_0}{\mu_0+\mu_2}\min\{v(i,0,1,1),v(i,0,0,2)\} + \frac{\mu_2}{ \mu_0+\mu_2} v(i-1, 2, 0, 0)\right) \nonumber \\
    &\quad = \frac{(i+1) (\mu_2-\mu_1) h_0}{(\mu_0+\mu_2)( \mu_0+\mu_1)} + \frac{h_1}{\mu_0+\mu_1} - \frac{h_2}{\mu_0+\mu_2}  \nonumber \\
    &\qquad + \frac{\mu_0}{\mu_0+\mu_2} \Big[\min\{v(i,0,2,0),v(i,0,1,1)\}-\min\{v(i,0,1,1),v(i,0,0,2)\}\Big] \nonumber \\
    &\qquad +  \frac{\mu_0\cdot(\mu_2-\mu_1)}{(\mu_0+\mu_2)( \mu_0+\mu_1)}\Big[\min\{v(i,0,2,0),v(i,0,1,1)\}-v(i-1,2,0,0)\Big] \label{eq:i200-diff}
\end{align}
For the difference $\min\{v(i,0,2,0),v(i,0,1,1)\}-\min\{v(i,0,1,1),v(i,0,0,2)\}$ in the fourth term in \eqref{eq:i200-diff}, replace the first minimum with an upper bound $v(i, 0, 2, 0)$ if $v(i, 0, 1, 1) \leq v(i, 0, 0, 2)$ and with $v(i, 0, 1, 1)$ otherwise. Therefore we have
\begin{align}
    \lefteqn {\min\{v(i,0,2,0),v(i,0,1,1)\} - \min\{v(i,0,1,1),v(i,0,0,2)\}}& \nonumber \\
    &\quad \leq \max\{v(i,0,2,0) - v(i,0,1,1), v(i,0,1,1) - v(i,0,0,2)\} \nonumber \\
    &\quad \leq \max\left\{\frac{h_1}{\mu_1} - \frac{h_2}{\mu_2} +\frac{ih_0}{2}\left(\frac{1}{\mu_1} - \frac{1}{\mu_2}\right),\frac{h_1}{\mu_1} - \frac{2h_2}{\mu_2} +\frac{ih_0}{2}\left(\frac{1}{\mu_1} - \frac{2}{\mu_2}\right)\right\} \nonumber \\
    &\quad = \frac{h_1}{\mu_1} - \frac{h_2}{\mu_2} +\frac{ih_0}{2}\left(\frac{1}{\mu_1} - \frac{1}{\mu_2}\right), \label{eq:decision_bd2_inte2}
\end{align}where the last inequality applies \eqref{eq:decision_bd2_110} and \eqref{eq:decision_bd2_101} that have been established. 

\noindent Consider the difference $\min\{v(i,0,2,0),v(i,0,1,1)\}-v(i-1,2,0,0)$ in the fifth term in \eqref{eq:i200-diff}. \textcolor{blue}{Suppose $a^*(i-1,2,0,0)$ is an optimal action at state $(i-1,2,0,0)$. We obtain an upper bound on the minimum by choosing the value function associated with applying action $a^*(i-1,2,0,0)$ at state $(i,1,1,0)$. Specifically,}
\begin{enumerate} [label= \textbf{Case} \arabic*:, leftmargin=2.5\parindent]
    \item {\color{blue} If $v(i-1,1,1,0) < v(i-1,1,0,1)$, choosing an upper bound $v(i,0,2,0)$ in the minimum $\min\{v(i,0,2,0),v(i,0,1,1)\}$ yields
    \begin{small}
    \begin{align*}
        \lefteqn{\min\{v(i,0,2,0),v(i,0,1,1)\} - v(i-1,2,0,0)}& \nonumber \\
        &\quad \leq v(i,0,2,0) - v(i-1,2,0,0) \nonumber\\
        &\quad = \left( \frac{ih_0+2h_1}{2\mu_1}+ v(i-1,1,1,0)\right) - \left(\frac{(i+1)h_0}{2\mu_0} + \min\{v(i-1,1,1,0),v(i-1,1,0,1)\}\right) \nonumber \\
        &\quad = \frac{ih_0+2h_1}{2\mu_1} - \frac{(i+1)h_0}{2\mu_0} \nonumber \\
        &\quad \leq \frac{h_1}{\mu_1} + \frac{(i+1)h_0}{2}\left(\frac{1}{\mu_1} - \frac{1}{\mu_0}\right), 
    \end{align*}
    \end{small}where the first equality is obtained by applying the optimality equations \eqref{eq:opt-zero} to states $(i, 0, 2, 0)$, and \eqref{eq:opt-zero}--\eqref{eq:opt1} to states $(i-1, 2, 0, 0)$, while the second equality follows from the case assumption.}
    \item \textcolor{blue}{If $v(i-1,1,1,0) \geq v(i-1,1,0,1)$, choosing an upper bound $v(i,0,1,1)$ in the minimum $\min\{v(i,0,2,0),v(i,0,1,1)\}$ yields
    \begin{align}
        \lefteqn{\min\{v(i,0,2,0),v(i,0,1,1)\} - v(i-1,2,0,0)}& \nonumber \\
        &\quad \leq v(i,0,1,1) - v(i-1,2,0,0). \label{eq:decision_bd2_case1_inte1}
    \end{align}}
    \textcolor{blue}{Applying the optimality equations \eqref{eq:opt-zero} to states $(i, 0, 1, 1)$, and \eqref{eq:opt-zero}--\eqref{eq:opt1} to states $(i-1, 2, 0, 0)$, we obtain}
    \textcolor{blue}{
    \begin{align}
        \lefteqn{v(i,0,1,1) - v(i-1,2,0,0)}& \nonumber \\
        &\quad = \left( \frac{ih_0+h_1+h_2}{\mu_1+\mu_2}+ \frac{\mu_1}{\mu_1+\mu_2}v(i-1,1,0,1)+ \frac{\mu_2}{\mu_1+\mu_2}v(i-1,1,1,0)\right) \nonumber \\
        &\qquad - \left(\frac{(i+1)h_0}{2\mu_0} + \min\{v(i-1,1,1,0),v(i-1,1,0,1)\}\right) \nonumber \\
        &\quad = \left( \frac{ih_0+h_1+h_2}{\mu_1+\mu_2}+ \frac{\mu_1}{\mu_1+\mu_2}v(i-1,1,0,1)+ \frac{\mu_2}{\mu_1+\mu_2}v(i-1,1,1,0)\right) \nonumber \\
        &\qquad - \left(\frac{(i+1)h_0}{2\mu_0} + v(i-1,1,0,1)\right) \nonumber \\
        &\quad = \frac{ih_0}{\mu_1+\mu_2} - \frac{(i+1)h_0}{2\mu_0} + \frac{h_1+h_2}{\mu_1+\mu_2} \nonumber\\
        &\qquad + \frac{\mu_2}{\mu_1+\mu_2}\Big[v(i-1,1,1,0)-v(i-1,1,0,1)\Big], \label{eq:decision_bd2_case1_inte2}
    \end{align}where the second-to-last step follows from the case assumption.} \textcolor{blue}{Applying the inductive hypothesis of \eqref{eq:decision_bd2_200} in \eqref{eq:decision_bd2_case1_inte2} yields
    \begin{align}
        \lefteqn{v(i,0,1,1) - v(i-1,2,0,0)}& \nonumber \\
        &\quad \leq \frac{ih_0}{\mu_1+\mu_2} - \frac{(i+1)h_0}{2\mu_0} + \frac{h_1+h_2}{\mu_1+\mu_2} + \frac{\mu_2}{\mu_1+\mu_2}\left(\frac{h_1}{\mu_1} - \frac{h_2}{\mu_2} + \frac{ih_0}{2}\left(\frac{1}{\mu_1} - \frac{1}{\mu_2}\right)\right)\nonumber\\
        &\quad = \frac{h_1}{\mu_1} + \frac{ih_0}{2\mu_1} - \frac{(i+1)h_0}{2\mu_0} \nonumber \\
        &\quad \leq \frac{h_1}{\mu_1} + \frac{(i+1)h_0}{2}\left(\frac{1}{\mu_1} - \frac{1}{\mu_0}\right).\label{eq:decision_bd2_case1_inte3}
    \end{align}}
    {\color{blue}Substituting \eqref{eq:decision_bd2_case1_inte3} into \eqref{eq:decision_bd2_case1_inte1} yields}
    \begin{align}
        \textcolor{blue}{\lefteqn{\min\{v(i,0,2,0),v(i,0,1,1)\} - v(i-1,2,0,0)}}& \nonumber\\
        &\quad \textcolor{blue}{\leq \frac{h_1}{\mu_1} + \frac{(i+1)h_0}{2}\left(\frac{1}{\mu_1} - \frac{1}{\mu_0}\right).} \label{eq:decision_bd2_inte3}
    \end{align}
\end{enumerate}
\textcolor{blue}{Therefore, in either case, inequality \eqref{eq:decision_bd2_inte3} holds.}
Using the bounds in \eqref{eq:decision_bd2_inte2} and \eqref{eq:decision_bd2_inte3} back in \eqref{eq:i200-diff} yields
\begin{align*}
    \lefteqn{D(i,2,0,0) = v(i, 1, 1, 0) - v(i, 1, 0, 1)}&\\
    &\quad \leq \frac{(i+1) (\mu_2-\mu_1) h_0}{(\mu_0+\mu_2)( \mu_0+\mu_1)} + \frac{h_1}{\mu_0+\mu_1} - \frac{h_2}{\mu_0+\mu_2}\\
    &\qquad + \frac{\mu_0}{\mu_0+\mu_2} \left(\frac{h_1}{\mu_1} - \frac{h_2}{\mu_2} + \frac{ih_0}{2}\left(\frac{1}{\mu_1} - \frac{1}{\mu_2}\right)\right)\\
    &\qquad +  \frac{\mu_0\cdot(\mu_2-\mu_1)}{(\mu_0+\mu_2)( \mu_0+\mu_1)}\left(\frac{h_1}{\mu_1} + \frac{(i+1)h_0}{2}\left(\frac{1}{\mu_1} - \frac{1}{\mu_0}\right)\right)\\
    &\quad \leq \frac{(i+1) (\mu_2-\mu_1) h_0}{(\mu_0+\mu_2)( \mu_0+\mu_1)} + \frac{h_1}{\mu_0+\mu_1} - \frac{h_2}{\mu_0+\mu_2}\\
    &\qquad + \frac{\mu_0}{\mu_0+\mu_2} \left(\frac{h_1}{\mu_1} - \frac{h_2}{\mu_2} + \frac{(i+1)h_0}{2}\left(\frac{1}{\mu_1} - \frac{1}{\mu_2}\right)\right)\\
    &\qquad +  \frac{\mu_0\cdot(\mu_2-\mu_1)}{(\mu_0+\mu_2)( \mu_0+\mu_1)}\left(\frac{h_1}{\mu_1} + \frac{(i+1)h_0}{2}\left(\frac{1}{\mu_1} - \frac{1}{\mu_0}\right)\right)\\
    &\quad = \frac{h_1}{\mu_1} - \frac{h_2}{\mu_2} + \frac{(i+1)h_0}{2}\left(\frac{1}{\mu_1} - \frac{1}{\mu_2}\right),
\end{align*}
where the last inequality uses $\mu_2 \geq \mu_1$.

{\color{blue}
Finally, we note that \eqref{eq:decision_bd2_inte} follows from a similar case analysis to that used to obtain \eqref{eq:decision_bd2_inte3}, by replacing $i-1$ with $i$ and using the corresponding established result in place of the inductive hypothesis of \eqref{eq:decision_bd2_200}.}
\end{proof}

\subsection{Proof of Theorem \ref{thm:large_queue}}
\begin{proof} [Proof of Theorem \ref{thm:large_queue}]
    \textcolor{blue}{Statement \ref{state:empty_Station2_slower} is implied by Statement \ref{state:larger_mu1_empty_Station2} of Proposition \ref{prop:thre_larger_mu1} below and Statement \ref{state3:general_rates} of Lemma \ref{lemma:general_rates}. 
    Statement \ref{state:empty_Station2_equal} is an immediate application of Statement \ref{state3:general_rates} in Lemma \ref{lemma:general_rates} and Proposition \ref{prop:decision_lb}. 
    Statements \ref{state:empty_Station2_faster} is supported by Propositions \ref{prop:empty-station2-large-i} and \ref{prop:decision_lb}.
    Finally, Statement \ref{state:queue_in_Station2_slower} is a consequence of Proposition \ref{prop:action1_large_i}, while Statements \ref{state:queue_in_Station2_faster_slow_initial} and \ref{state:queue_in_Station2_faster_fast_initial} follow from Statements \ref{state:queue_Station2_action2_large_i} and \ref{state:queue_Station2_action1_large_i} in Proposition \ref{prop:queue_Station2_large_i}, respectively.}
\end{proof}
% Since the state variables $j$, $k$, and $\ell$ vary over a finite set, Statement \ref{state:main-one-faster} in Theorem \ref{thm:large_queue} is implied by the following proposition.
\begin{proposition}
    \label{prop:thre_larger_mu1}
    \textcolor{blue}{
    The following results hold for decision states $(i,j,k,\ell) \in \X_{\widehat{D}}$. 
    \begin{enumerate}
        \item \label{state:larger_mu1_queue_in_Station2}
        Suppose $\mu_1 \geq \mu_2$ and consider $\ell = 1$, i.e., the decision states of the form $(i,1,0,1)$. There exists a finite $i'_{1,0,1}$ such that $D(i,1,0,1) \leq 0$ for all $i \geq i'_{1,0,1}$.
        \item \label{state:larger_mu1_empty_Station2}
        Suppose $\mu_1 > \mu_2$ and consider $\ell = 0$, i.e., the decision states of the form $(i,1,1,0)$ or $(i,2,0,0)$. There exists a finite $i'_{j,k,\ell}$ such that $D(i,j,k,\ell) \leq 0$ for all $i \geq i'_{j,k,\ell}$.
    \end{enumerate}}
\end{proposition}

\begin{proof} [Proof of Proposition \ref{prop:thre_larger_mu1}]
\textcolor{blue}{For Statement \ref{state:larger_mu1_queue_in_Station2}, suppose $\mu_1 \geq \mu_2$ and} consider the decision states $(i,1,0,1)$. Without loss of generality assume $i \geq 1$ since we only focus on the behaviour when $i$ is large. 
\textcolor{blue}{Note that \eqref{eq:i101-diff} follows directly from the optimality equation \eqref{eq:opt1}, regardless of whether $\mu_1 \geq \mu_2$:}
\begin{align}
    \lefteqn{D(i,1,0,1) = v(i, 0, 1, 1) - v(i, 0, 0, 2)}& \nonumber \\
    & \quad = \frac{-i \mu_1 h_0}{\mu_2( \mu_1+\mu_2)} + \textcolor{blue}{\frac{h_1+h_2}{\mu_1+\mu_2} - \frac{2h_2}{\mu_2}} + \frac{ \mu_2}{ \mu_1+\mu_2} \Big[ v(i-1, 1, 1, 0) - v(i-1,1, 0, 1) \Big] \nonumber \\
    &\quad \leq \frac{-i \mu_1 h_0}{\mu_2( \mu_1+\mu_2)} + \textcolor{blue}{\frac{h_1+h_2}{\mu_1+\mu_2} - \frac{2h_2}{\mu_2}} + \frac{\mu_2}{ \mu_1+\mu_2} \left(\frac{h_1}{\mu_1} - \frac{h_2}{\mu_2}\right)\nonumber \\
    &\quad = \textcolor{blue}{\left(\frac{h_1}{\mu_1} - \frac{2h_2}{\mu_2}\right) + \frac{-i \mu_1 h_0}{\mu_2( \mu_1+\mu_2)}},\label{eq:i101-diff-final-bound2}
\end{align}
where the inequality holds by \eqref{eq:decision_bd1} since $\mu_1 \geq \mu_2$.\\
Defining $c_1 := \textcolor{blue}{\frac{h_1}{\mu_1} - \frac{2h_2}{\mu_2}}$ 
yields that 
\begin{align}
  D(i,1,0,1) = v(i, 0, 1, 1) - v(i, 0, 0, 2) & \leq  \frac{-i \mu_1 h_0}{\mu_2( \mu_1+\mu_2)} + c_1. \label{eq:bound-c1}
\end{align}
Since $c_1$ is a constant independent of $i$, defining $i'_{1,0,1} := c_1\frac{\mu_2(\mu_1+\mu_2)}{\mu_1h_0}$, implies for any $i \geq i'_{1,0,1}$ we have $v(i, 0, 1, 1) \leq v(i, 0, 0, 2)$ as desired.

\textcolor{blue}{We now suppose $\mu_1>\mu_2$ and prove Statement \ref{state:larger_mu1_empty_Station2} for states $(i,1,1,0)$ and $(i,2,0,0)$, where $i \geq 1$.
Since \eqref{eq:i110-diff} also holds regardless of whether $\mu_1 > \mu_2$,} the proof for $(i,1,1,0)$ is analogous to that for $(i,1,0,1)$ and is omitted for brevity. The comparable bound in this case is
\begin{align}
    D(i,1,1,0) = v(i, 0, 2, 0) - v(i, 0, 1, 1) & \leq  \frac{-i (\mu_1-\mu_2) h_0}{2\mu_1( \mu_1+\mu_2)} + c_2, \label{eq:bound-c2}
\end{align}
where $c_2 := \textcolor{blue}{\frac{h_1}{\mu_1} - \frac{h_2}{\mu_2}}$. \textcolor{blue}{Since $\mu_1 > \mu_2$ implies $\frac{(\mu_1-\mu_2) h_0}{2\mu_1( \mu_1+\mu_2)} > 0$, the result follows from letting $i'_{1,1,0} = c_2 \frac{2\mu_1( \mu_1+\mu_2)} {(\mu_1-\mu_2) h_0}$.}
% For states $(i,1,1,0)$,
% \begin{align*}
%    \lefteqn{v(i, 0, 2, 0) - v(i, 0, 1, 1)} \\
%        &\quad = \frac{i h_0+2h_1}{2\mu_1} - \frac{i h_0+h_1+h_2}{\mu_1+\mu_2}\\
%        & \qquad + v(i-1, 1, 1, 0) - \Big[ \frac{\mu_1}{\mu_1+\mu_2} v(i-1, 1, 0, 1) + \frac{\mu_2}{ \mu_1+\mu_2} v(i-1, 1, 1, 0)\Big]\\
%        & \quad =  \frac{-i (\mu_1-\mu_2) h_0}{2\mu_1( \mu_1+\mu_2)} + \frac{\mu_2h_1 - \mu_1h_2}{\mu_1( \mu_1+\mu_2)} + \frac{\mu_1}{ \mu_1+\mu_2} \Big[v(i-1, 1, 1, 0) - v(i-1, 1, 0, 1)\Big]\\
%        & \quad \leq \frac{-i (\mu_1-\mu_2) h_0}{2\mu_1( \mu_1+\mu_2)} + \frac{\mu_2h_1 - \mu_1h_2}{\mu_1( \mu_1+\mu_2)} +  \frac{\mu_1}{ \mu_1+\mu_2} \Big(\frac{h_1}{\mu_1} - \frac{h_2}{\mu_2}\Big),
% \end{align*}
% where the inequality holds by \eqref{eq:decision_bd1}. Let $c_2 := \frac{\mu_2h_1 - \mu_1h_2}{\mu_1( \mu_1+\mu_2)} +  \frac{\mu_1}{ \mu_1+\mu_2} \Big(\frac{h_1}{\mu_1} - \frac{h_2}{\mu_2}\Big)$ and note that it is a constant independent of $i$. Thus, 
% \begin{align*}
%          v(i, 0, 2, 0) - v(i, 0, 1, 1) & \leq  \frac{-i (\mu_1-\mu_2) h_0}{2\mu_1( \mu_1+\mu_2)} + c_2,
% \end{align*}
% Defining $i'_{1,1,0} := c_2\frac{2\mu_1(\mu_1+\mu_2)}{(\mu_1-\mu_2)h_0}$ yields the result. 
Now consider the decision states $(i,2,0,0)$, thinning the MDPs by a higher rate $\mu_0+\mu_1$ yields
\begin{align*}
    \lefteqn{D(i,2,0,0) = v(i, 1, 1, 0) - v(i, 1, 0, 1)}\\
        & \quad = \frac{h_1- h_2}{\mu_0+\mu_1}+ \frac{\mu_0}{\mu_0+\mu_1} \Big[\min\{v(i, 0, 2, 0),v(i, 0, 1, 1)\} - \min\{v(i, 0, 1, 1),v(i, 0, 0, 2)\}\Big]\\
        & \qquad + \frac{\mu_2}{\mu_0+\mu_1} \Big[v(i-1, 2, 0, 0) - v(i-1, 2, 0, 0)\Big] + \frac{\mu_1-\mu_2}{\mu_0+\mu_1} \Big[v(i-1, 2, 0, 0) - v(i, 1, 0, 1)\Big]\\
        & \quad \leq \frac{h_1- h_2}{\mu_0+\mu_1} + \frac{\mu_1-\mu_2}{\mu_0+\mu_1}\textcolor{blue}{\left(- \min\left\{\frac{h_1}{\mu_1},\frac{h_2}{\mu_2}\right\}\right)}\\
        & \qquad + \frac{\mu_0}{\mu_0+\mu_1} \Big[\min\{v(i, 0, 2, 0),v(i, 0, 1, 1)\} - \min\{v(i, 0, 1, 1),v(i, 0, 0, 2)\}\Big]
\end{align*}
where the inequality follows from \eqref{eq:onemore2} since $(i-1,2,0,0) - v(i,1,0,1) \leq \textcolor{blue}{- \min\left\{\frac{h_1}{\mu_1},\frac{h_2}{\mu_2}\right\}}$. There are several cases to consider. 
\begin{enumerate}[label= \textbf{Case} \arabic*:, leftmargin=3\parindent]
    \item If $v(i, 0, 1, 1) \leq v(i, 0, 0, 2)$, replacing the first minimum with an upper bound $v(i, 0, 2, 0)$ yields
        \begin{align*}
            \lefteqn{v(i,1,1,0) - v(i,1,0,1)}\\
            & \quad \leq \frac{\mu_0}{\mu_0+\mu_1} \Big[v(i, 0, 2, 0) - v(i, 0, 1, 1)\Big] + \textcolor{blue}{\frac{h_1- h_2}{\mu_0+\mu_1} + \frac{\mu_1-\mu_2}{\mu_0+\mu_1}\left(- \min\left\{\frac{h_1}{\mu_1},\frac{h_2}{\mu_2}\right\}\right)}\\
            & \quad \leq \frac{\mu_0}{\mu_0+\mu_1} \left(\frac{-i (\mu_1-\mu_2) h_0}{2\mu_1( \mu_1+\mu_2)} + c_2 \right) + \textcolor{blue}{\frac{h_1- h_2}{\mu_0+\mu_1} + \frac{\mu_1-\mu_2}{\mu_0+\mu_1}\left(- \min\left\{\frac{h_1}{\mu_1},\frac{h_2}{\mu_2}\right\}\right)}\\
            & \quad = \frac{-i\mu_0(\mu_1-\mu_2)h_0}{2\mu_1(\mu_0+\mu_1)(\mu_1+\mu_2)} + c_3,
        \end{align*}
    where the last inequality follows from \eqref{eq:bound-c2}. Here, $c_3$ is a constant independent of $i$, defined as $c_3 = \frac{\mu_0}{\mu_0+\mu_1}c_2 + \textcolor{blue}{\frac{h_1- h_2}{\mu_0+\mu_1} + \frac{\mu_1-\mu_2}{\mu_0+\mu_1}\left(- \min\left\{\frac{h_1}{\mu_1},\frac{h_2}{\mu_2}\right\}\right)}$.

    \item If $v(i, 0, 1, 1) \geq v(i, 0, 0, 2)$, replace the first minimum by $v(i, 0, 1, 1)$.
    \begin{align*}
        \lefteqn{v(i,1,1,0) - v(i,1,0,1)}\\
            & \quad \leq \frac{\mu_0}{\mu_0+\mu_1} \Big[v(i, 0, 1, 1) - v(i, 0, 0, 2)\Big] + \textcolor{blue}{\frac{h_1- h_2}{\mu_0+\mu_1} + \frac{\mu_1-\mu_2}{\mu_0+\mu_1}\left(- \min\left\{\frac{h_1}{\mu_1},\frac{h_2}{\mu_2}\right\}\right)}\\
            & \quad \leq \frac{\mu_0}{\mu_0+\mu_1} \left(\frac{-i \mu_1 h_0}{\mu_2( \mu_1+\mu_2)} + c_1\right) + \textcolor{blue}{\frac{h_1- h_2}{\mu_0+\mu_1} + \frac{\mu_1-\mu_2}{\mu_0+\mu_1}\left(- \min\left\{\frac{h_1}{\mu_1},\frac{h_2}{\mu_2}\right\}\right)}\\
            & \quad = \frac{-i\mu_0\mu_1h_0}{\mu_2(\mu_0+\mu_1)(\mu_1+\mu_2)} + c_4,
    \end{align*}
    where the last inequality holds by \eqref{eq:bound-c1}. Here, $c_4$ is another constant independent of $i$, defined as $c_4 = \frac{\mu_0}{\mu_0+\mu_1}c_1 + \textcolor{blue}{\frac{h_1- h_2}{\mu_0+\mu_1} + \frac{\mu_1-\mu_2}{\mu_0+\mu_1}\left(- \min\left\{\frac{h_1}{\mu_1},\frac{h_2}{\mu_2}\right\}\right)}$.
\end{enumerate}
The result holds for the decision states $(i,2,0,0)$ by defining $i'_{2,0,0}$ as follows
\begin{align*}
    i'_{2,0,0} = \max\left\{\frac{2\mu_1(\mu_0+\mu_1)(\mu_1+\mu_2)}{\mu_0(\mu_1-\mu_2)h_0}c_3,\frac{\mu_2(\mu_0+\mu_1)(\mu_1+\mu_2)}{\mu_0\mu_1h_0}c_4\right\}.
\end{align*}
\end{proof}

\begin{proposition} \label{prop:decision_lb}
    Suppose $\mu_2 \geq \mu_1$ \textcolor{blue}{and consider decision states $(i,j,k,\ell) \in \X_{\widehat{D}}$ such that $\ell =0$, i.e., states of the form $(i,1,1,0)$ and $(i,2,0,0)$. The following holds for all $i$:
    \begin{align*}
        D(i,j,k,\ell) \geq \min\left\{\frac{h_1}{\mu_1} - \frac{h_2}{\mu_2},0\right\}.
    \end{align*}
    This is equivalent to:}
    \begin{align} 
        \textcolor{blue}{v(i,0,2,0) - v(i,0,1,1) \geq \min\left\{\frac{h_1}{\mu_1} - \frac{h_2}{\mu_2},0\right\},} \label{eq:decision_lb_i110}\\
        \textcolor{blue}{v(i,1,1,0) - v(i,1,0,1) \geq \min\left\{\frac{h_1}{\mu_1} - \frac{h_2}{\mu_2},0\right\}.} \label{eq:decision_lb_i200}
    \end{align}
\end{proposition}
\begin{proof}[Proof of Proposition \ref{prop:decision_lb}]
We prove the result by induction on the number of jobs in the system $i$. When $i = 0$, the base case for \eqref{eq:decision_lb_i110} is
\begin{align*}
    v(0,0,2,0) - v(0,0,1,1) = \frac{2h_1}{\mu_1} - \left(\frac{h_1}{\mu_1} + \frac{h_2}{\mu_2}\right) = \frac{h_1}{\mu_1} - \frac{h_2}{\mu_2} \geq \textcolor{blue}{\min\left\{\frac{h_1}{\mu_1} - \frac{h_2}{\mu_2},0\right\}}.
\end{align*}
As for the base case $i=0$ for \eqref{eq:decision_lb_i200}, thinning the MDPs by a higher rate $\mu_0+\mu_2$ yields
\begin{align}
    & v(0,1,1,0) - v(0,1,0,1) \nonumber \\
    & \quad = \frac{h_1-h_2}{\mu_0+\mu_2} +\frac{\mu_0}{\mu_0+\mu_2} \Big[ \min\{v(0,0,2,0),v(0,0,1,1)\}- \min\{v(0,0,1,1),v(0,0,0,2)\} \Big] \nonumber \\
    & \quad \qquad+ \frac{\mu_1}{\mu_0+\mu_2}\Big[v(0,1,0,0) - v(0,1,0,0)\Big] + \frac{\mu_2-\mu_1}{\mu_0+\mu_2}\Big[v(0,1,1,0) - v(0,1,0,0)\Big]. \label{eq:decision_diff_0200}
\end{align}
\textcolor{blue}{Consider the difference $\min\{v(0,0,2,0),v(0,0,1,1)\}- \min\{v(0,0,1,1),v(0,0,0,2)\}$ in the second minimum in \eqref{eq:decision_diff_0200}.} 
Replace the second minimum with the upper bound  $v(0,0,1,1)$ to obtain
\textcolor{blue}{
    \begin{align}
    \lefteqn{\min\{v(0,0,2,0),v(0,0,1,1)\}- \min\{v(0,0,1,1),v(0,0,0,2)\}}&\nonumber\\
    & \quad \geq \min\{v(0,0,2,0),v(0,0,1,1)\}- v(0,0,1,1) \nonumber\\
    & \quad = \min\{v(0,0,2,0) -v(0,0,1,1), 0\} \nonumber\\
    &\quad \geq \min\left\{\min\left\{\frac{h_1}{\mu_1} - \frac{h_2}{\mu_2},0\right\}, 0\right\} \nonumber \\
    &\quad = \min\left\{\frac{h_1}{\mu_1} - \frac{h_2}{\mu_2},0\right\}, \label{eq:decision_diff_0200_inte1}
\end{align}where the last inequality follows from the previous base case for \eqref{eq:decision_lb_i110}. Applying \eqref{eq:decision_diff_0200_inte1} in \eqref{eq:decision_diff_0200} yields}
\begin{align*}
    & v(0,1,1,0) - v(0,1,0,1) \\
    & \quad \geq  \frac{h_1-h_2}{\mu_0+\mu_2} \textcolor{blue}{+ \frac{\mu_0}{\mu_0+\mu_2}\left(\min\left\{\frac{h_1}{\mu_1} - \frac{h_2}{\mu_2},0\right\}\right)} + \frac{\mu_2-\mu_1}{\mu_0+\mu_2}\Big[v(0,1,1,0) - v(0,1,0,0)\Big]\\
    & \quad \geq  \frac{h_1 - h_2}{\mu_0+\mu_2} \textcolor{blue}{+ \frac{\mu_0}{\mu_0+\mu_2}\left(\min\left\{\frac{h_1}{\mu_1} - \frac{h_2}{\mu_2},0\right\}\right)} + \left(\frac{\mu_2-\mu_1}{\mu_0+\mu_2}\right)\frac{h_1}{\mu_1}\\
    & \quad =\frac{\mu_2}{\mu_0+\mu_2}\left(\frac{h_1}{\mu_1}-\frac{h_2}{\mu_2}\right)\textcolor{blue}{+ \frac{\mu_0}{\mu_0+\mu_2}\left(\min\left\{\frac{h_1}{\mu_1} - \frac{h_2}{\mu_2},0\right\}\right)}\\
    &\quad \geq \textcolor{blue}{\min\left\{\frac{h_1}{\mu_1} - \frac{h_2}{\mu_2},0\right\}},
\end{align*}where the second inequality is by \eqref{eq:empty-queue-bound1} since $v(0,1,1,0) - v(0,1,0,0) \geq \frac{h_1}{\mu_1}$.

Now suppose both inequalities \eqref{eq:decision_lb_i110} and \eqref{eq:decision_lb_i200} are true for $i-1$, where $i \geq 1$, and prove them at $i$. Thinning the MDPs in \eqref{eq:decision_lb_i110} with a higher rate $\mu_1+\mu_2$ yields
\begin{align*}
    \lefteqn{v(i,0,2,0) - v(i,0,1,1)} \\
    & \quad = \frac{h_1-h_2}{\mu_1+\mu_2} +     \frac{\mu_1}{\mu_1+\mu_2}\Big[v(i-1,1,1,0) - v(i-1,1,0,1)\Big]\\
    & \qquad+ \frac{\mu_1}{\mu_1+\mu_2}\Big[v(i-1,1,1,0) - v(i-1,1,1,0)\Big]+ \frac{\mu_2-\mu_1}{\mu_1+\mu_2}\Big[v(i,0,2,0) - v(i-1,1,1,0)\Big]\\
    & \quad \geq  \frac{h_1 - h_2}{\mu_1+\mu_2} \textcolor{blue}{+\frac{\mu_1}{\mu_1+\mu_2}\left(\min\left\{\frac{h_1}{\mu_1} - \frac{h_2}{\mu_2},0\right\}\right)} + \frac{\mu_2-\mu_1}{\mu_1+\mu_2}\Big[v(i,0,2,0) - v(i-1,1,1,0)\Big]\\
    & \quad \geq  \frac{h_1 - h_2}{\mu_1+\mu_2} \textcolor{blue}{+\frac{\mu_1}{\mu_1+\mu_2}\left(\min\left\{\frac{h_1}{\mu_1} - \frac{h_2}{\mu_2},0\right\}\right)} + \frac{\mu_2-\mu_1}{\mu_1+\mu_2}\textcolor{blue}{\left(\min\left\{\frac{h_1}{\mu_1}, \frac{h_2}{\mu_2}\right\}\right)}\\
    & \quad =\textcolor{blue}{\frac{\mu_1}{\mu_1+\mu_2}}\left(\frac{h_1}{\mu_1}-\frac{h_2}{\mu_2}\right) \textcolor{blue}{+\frac{\mu_2}{\mu_1+\mu_2}\left(\min\left\{\frac{h_1}{\mu_1} - \frac{h_2}{\mu_2},0\right\}\right)}\\
    &\quad \geq \textcolor{blue}{\min\left\{\frac{h_1}{\mu_1} - \frac{h_2}{\mu_2},0\right\}},
\end{align*}
where the first inequality follows from applying the inductive hypothesis at $i-1$ for \eqref{eq:decision_lb_i200}, and the second inequality is by \eqref{eq:onemore1}. This proves \eqref{eq:decision_lb_i110}.  For \eqref{eq:decision_lb_i200}, thinning the MDPs with a higher rate $\mu_0+\mu_2$ yields
\begin{align}
    & v(i,1,1,0) - v(i,1,0,1) \nonumber \\
    & \quad = \frac{h_1-h_2}{\mu_0+\mu_2} +\frac{\mu_0}{\mu_0+\mu_2} \Big[ \min\{v(i,0,2,0),v(i,0,1,1)\}- \min\{v(i,0,1,1),v(i,0,0,2)\} \Big] \nonumber \\
    & \quad \qquad+ \frac{\mu_1}{\mu_0+\mu_2}\Big[v(i-1,2,0,0) - v(i-1,2,0,0)\Big] \nonumber \\
    & \quad \qquad+ \frac{\mu_2-\mu_1}{\mu_0+\mu_2}\Big[v(i,1,1,0) - v(i-1,2,0,0)\Big]. \label{eq:decision_diff_i200}
\end{align}
\textcolor{blue}{By an argument similar to that used to obtain \eqref{eq:decision_diff_0200_inte1}, the following holds:
\begin{align}
    & \min\{v(i,0,2,0),v(i,0,1,1)\}- \min\{v(i,0,1,1),v(i,0,0,2)\} \geq \min\left\{\frac{h_1}{\mu_1} - \frac{h_2}{\mu_2},0\right\}. \label{eq:decision_lb_i200_inte1}
\end{align}
Subsitituting \eqref{eq:decision_lb_i200_inte1} in \eqref{eq:decision_diff_i200} yields}
\begin{align*}
    \lefteqn{v(i,1,1,0) - v(i,1,0,1)}  \\
    & \quad \geq  \frac{h_1-h_2}{\mu_0+\mu_2} \textcolor{blue}{+ \frac{\mu_0}{\mu_0+\mu_2}\left(\min\left\{\frac{h_1}{\mu_1} - \frac{h_2}{\mu_2},0\right\}\right)} + \frac{\mu_2-\mu_1}{\mu_0+\mu_2}\Big[v(i,1,1,0) - v(i-1,2,0,0)\Big]\\
    & \quad \geq  \frac{h_1 - h_2}{\mu_0+\mu_2} \textcolor{blue}{+ \frac{\mu_0}{\mu_0+\mu_2}\left(\min\left\{\frac{h_1}{\mu_1} - \frac{h_2}{\mu_2},0\right\}\right)}  + \frac{\mu_2-\mu_1}{\mu_0+\mu_2}\textcolor{blue}{\left(\min\left\{\frac{h_1}{\mu_1}, \frac{h_2}{\mu_2}\right\}\right)}\\
    & \quad =\frac{\mu_1}{\mu_0+\mu_2}\left(\frac{h_1}{\mu_1}-\frac{h_2}{\mu_2}\right) \textcolor{blue}{+ \frac{\mu_0+\mu_2-\mu_1}{\mu_0+\mu_2}\left(\min\left\{\frac{h_1}{\mu_1} - \frac{h_2}{\mu_2},0\right\}\right)}\\
    &\quad \geq \textcolor{blue}{\min\left\{\frac{h_1}{\mu_1} - \frac{h_2}{\mu_2},0\right\}},
\end{align*}
where the second inequality follows from \eqref{eq:onemore1}, thereby confirming \eqref{eq:decision_lb_i200}.
\end{proof}
\begin{proposition} \label{prop:empty-station2-large-i}
    \textcolor{blue}{
    Suppose $\mu_1 < \mu_2$ and consider decision states $(i,j,k,\ell) \in \X_{\widehat{D}}$ such that $\ell=0$, i.e., states of the form $(i,1,1,0)$ and $(i,2,0,0)$. There exists a finite $i'_{j,k,\ell}$ such that $D(i,j,k,\ell) \geq 0$ for all $i \geq i'_{j,k,\ell}$.}
\end{proposition}

\begin{proof} [\textcolor{blue}{Proof of Proposition \ref{prop:empty-station2-large-i}.}]
    \color{blue}{Without loss of generality, we assume $i\geq1$, since we only focus on the sign of $D$ for sufficiently large $i$. For the decision states $(i,1,1,0)$ with $i \geq 1$, equation \eqref{eq:i110-diff} holds without assumptions:
    \begin{align}
    \lefteqn{D(i,1,1,0) = v(i, 0, 2, 0) - v(i, 0, 1, 1)}& \nonumber \\
    &\quad = \frac{i (\mu_2-\mu_1) h_0}{2\mu_1(\mu_1+\mu_2)} + \frac{\mu_2h_1 - \mu_1h_2}{\mu_1(\mu_1+\mu_2)} + \frac{\mu_1}{ \mu_1+\mu_2} \Big[v(i-1, 1, 1, 0) - v(i-1, 1, 0, 1)\Big] \nonumber \\
    &\quad \geq \frac{i (\mu_2-\mu_1) h_0}{2\mu_1(\mu_1+\mu_2)} + \frac{\mu_2h_1 - \mu_1h_2}{\mu_1(\mu_1+\mu_2)} + \frac{\mu_1}{ \mu_1+\mu_2}\left(\min\left\{\frac{h_1}{\mu_1} - \frac{h_2}{\mu_2},0\right\}\right) \nonumber \\
    &\quad = \frac{i (\mu_2-\mu_1) h_0}{2\mu_1(\mu_1+\mu_2)} + \tilde{c}_1, \label{eq2:i110-diff}
\end{align}where the inequality follows from \eqref{eq:decision_lb_i200} in Proposition \ref{prop:decision_lb}. Here, $\tilde{c}_1$, independent of $i$, is defined as $\tilde{c}_1 = \frac{\mu_2h_1 - \mu_1h_2}{\mu_1(\mu_1+\mu_2)} + \frac{\mu_1}{ \mu_1+\mu_2}\left(\min\left\{\frac{h_1}{\mu_1} - \frac{h_2}{\mu_2},0\right\}\right)$. Since $\mu_1 < \mu_2$, the right-hand side of \eqref{eq2:i110-diff} becomes non-negative for all sufficiently large $i$, and the result follows.}

\textcolor{blue}{Consider now the decision states $(i,2,0,0)$ where $i\geq 1$.
Since $\mu_2 > \mu_1$, equation \eqref{eq:i200-diff} continues to hold for $i \geq 1$. That is,
\begin{align}
    \lefteqn{D(i,2,0,0) = v(i, 1, 1, 0) - v(i, 1, 0, 1)} \nonumber\\
    & \quad =\frac{(i+1) (\mu_2-\mu_1) h_0}{(\mu_0+\mu_2)( \mu_0+\mu_1)} + \frac{h_1}{\mu_0+\mu_1} - \frac{h_2}{\mu_0+\mu_2} \nonumber\\
    & \qquad + \ \frac{\mu_0}{\mu_0+\mu_2} \Big[\min\{v(i,0,2,0),v(i,0,1,1)\}-\min\{v(i,0,1,1),v(i,0,0,2)\}\Big] \nonumber\\
    & \qquad +   \frac{\mu_0\cdot(\mu_2-\mu_1)}{(\mu_0+\mu_2)( \mu_0+\mu_1)}\Big[\min\{v(i,0,2,0),v(i,0,1,1)\}-v(i-1,2,0,0)\}\Big]. \label{eq3:i200-diff}
\end{align}}%
Using arguments similar to those used to obtain \eqref{eq:decision_diff_0200_inte1} or \eqref{eq:decision_lb_i200_inte1}, both of which rely on \eqref{eq:decision_lb_i110} in Proposition \ref{prop:decision_lb}, the difference $\min\{v(i,0,2,0),v(i,0,1,1)\}-\min\{v(i,0,1,1),v(i,0,0,2)\}$ in the fourth term of \eqref{eq3:i200-diff} admits the following lower bound: 
\begin{align}
    & \min\{v(i,0,2,0),v(i,0,1,1)\}- \min\{v(i,0,1,1),v(i,0,0,2)\} \geq \min\left\{\frac{h_1}{\mu_1} - \frac{h_2}{\mu_2},0\right\}. \label{eq3:i200-diff-fourth-term}
\end{align}
% \textcolor{blue}{
% For the difference $\min\{v(i,0,2,0),v(i,0,1,1)\}-\min\{v(i,0,1,1),v(i,0,0,2)\}$ in the fourth term in \eqref{eq3:i200-diff}, replacing the second minimum with an upper bound $v(i, 0, 1, 1)$ gives
% \begin{align}
%     \lefteqn {\min\{v(i,0,2,0),v(i,0,1,1)\} - \min\{v(i,0,1,1),v(i,0,0,2)\}}& \nonumber \\
%     &\quad \geq \min\{v(i,0,2,0), v(i,0,1,1)\} - v(i,0,1,1) \nonumber \\
%     &\quad = \min\{v(i,0,2,0) - v(i,0,1,1), 0\}  \nonumber \\
%     &\quad \geq \min\left\{\min\left\{\frac{h_1}{\mu_1} - \frac{h_2}{\mu_2},0\right\}, 0\right\}  \nonumber \\
%     &\quad = \min\left\{\frac{h_1}{\mu_1} - \frac{h_2}{\mu_2},0\right\}, \label{eq3:i200-diff-fourth-term}
% \end{align}where the last inequality follows from \eqref{eq:decision_lb_i110} in Proposition \ref{prop:decision_lb} (recall $\mu_2 > \mu_1$).}
\textcolor{blue}{
Consider the difference $\min\{v(i,0,2,0),v(i,0,1,1)\}-v(i-1,2,0,0)\}$ in the fifth term in \eqref{eq3:i200-diff}.}
\begin{enumerate} [label= \textbf{Case} \arabic*:, leftmargin=2.5\parindent]
    \item \textcolor{blue}{Suppose $v(i,0,2,0) \geq v(i,0,1,1)$. The following holds:} 
    \textcolor{blue}{
    \begin{small}
    \begin{align}
        \lefteqn{\min\{v(i,0,2,0),v(i,0,1,1)\} - v(i-1,2,0,0)}& \nonumber \\
        &\quad =v(i,0,1,1) - v(i-1,2,0,0) \nonumber \\
        &\quad = \left( \frac{ih_0+h_1+h_2}{\mu_1+\mu_2}+ \frac{\mu_1}{\mu_1+\mu_2}v(i-1,1,0,1)+ \frac{\mu_2}{\mu_1+\mu_2}v(i-1,1,1,0)\right) \nonumber \\
        &\qquad - \left(\frac{(i+1)h_0}{2\mu_0} + \min\{v(i-1,1,1,0),v(i-1,1,0,1)\}\right) \nonumber \\
        &\quad \geq   \frac{ih_0+h_1+h_2}{\mu_1+\mu_2} - \frac{(i+1)h_0}{2\mu_0}+ \frac{\mu_2}{\mu_1+\mu_2}\Big[v(i-1,1,1,0)-v(i-1,1,0,1)\Big], \label{eq3:i200-diff-case1}
    \end{align}
    \end{small}where the second-to last step holds by applying the optimality equations \eqref{eq:opt-zero} to states $(i, 0, 1, 1)$ and \eqref{eq:opt-zero}--\eqref{eq:opt1} to states $(i-1, 2, 0, 0)$, while the last step follows from choosing an upper bound $v(i-1,1,0,1)$ in the minimum.} 
    \textcolor{blue}{Applying \eqref{eq:decision_lb_i200} in Proposition \ref{prop:decision_lb} into \eqref{eq3:i200-diff-case1} yields
    \begin{align}
        \lefteqn{\min\{v(i,0,2,0),v(i,0,1,1)\} - v(i-1,2,0,0)}& \nonumber \\
        & \quad \geq \frac{ih_0+h_1+h_2}{\mu_1+\mu_2} - \frac{(i+1)h_0}{2\mu_0} + \frac{\mu_2}{\mu_1+\mu_2}\left(\min\left\{\frac{h_1}{\mu_1} - \frac{h_2}{\mu_2},0\right\}\right) \nonumber \\
        &\quad = \left(\frac{1}{\mu_1+\mu_2}-\frac{1}{2\mu_0}\right)(i+1)h_0 + \tilde{c}_2 \nonumber \\
        &\quad \geq \left(\frac{1}{\mu_2}-\frac{1}{\mu_0}\right)\frac{(i+1)h_0}{2} + \tilde{c}_2, \label{eq:tilde-c2}
    \end{align}where the last step follows from $\frac{1}{\mu_1+\mu_2}\geq \frac{1}{2\mu_2}$ since $\mu_1 < \mu_2$. Here, $\tilde{c}_2$ is a constant independent of $i$, defined as $\tilde{c}_2 = \frac{h_1+h_2-h_0}{\mu_1+\mu_2}+ \frac{\mu_2}{\mu_1+\mu_2}\left(\min\left\{\frac{h_1}{\mu_1} - \frac{h_2}{\mu_2},0\right\}\right)$.}
    \item \textcolor{blue}{If, instead, $v(i,0,2,0) < v(i,0,1,1)$, then we have
    \begin{small}
    \begin{align*}
        \lefteqn{\min\{v(i,0,2,0),v(i,0,1,1)\} - v(i-1,2,0,0)}& \nonumber \\
        &\quad = v(i,0,2,0) - v(i-1,2,0,0) \\
        &\quad = \left( \frac{ih_0+2h_1}{2\mu_1}+ v(i-1,1,1,0)\right) - \left(\frac{(i+1)h_0}{2\mu_0} + \min\{v(i-1,1,1,0),v(i-1,1,0,1)\}\right) \nonumber \\
        &\quad \geq \frac{ih_0+2h_1}{2\mu_1} - \frac{(i+1)h_0}{2\mu_0},
    \end{align*}
    \end{small}where the second-to last step holds by applying the optimality equations \eqref{eq:opt-zero} to states $(i, 0, 2, 0)$ and \eqref{eq:opt-zero}--\eqref{eq:opt1} to states $(i-1, 2, 0, 0)$, while the final inequality comes from replacing the minimum with an upper bound $v(i-1,1,1,0)$.} 
    \textcolor{blue}{Using $\mu_1 < \mu_2$ again yields
    \begin{align}
        \lefteqn{\min\{v(i,0,2,0),v(i,0,1,1)\} - v(i-1,2,0,0)}& \nonumber \\
        &\quad \geq \frac{ih_0+2h_1}{2\mu_2} - \frac{(i+1)h_0}{2\mu_0} \nonumber\\
        &\quad = \left(\frac{1}{\mu_2}-\frac{1}{\mu_0}\right)\frac{(i+1)h_0}{2} + \tilde{c}_2', \label{eq:tilde-c2-prime}
    \end{align}where $\tilde{c}_2' = \frac{2h_1-h_0}{2\mu_1}$ is another constant independent of $i$.}
\end{enumerate}
\textcolor{blue}{
Comparing the lower bounds in \eqref{eq:tilde-c2} and \eqref{eq:tilde-c2-prime}, it follows that 
\begin{align}
    \min\{v(i,0,2,0),v(i,0,1,1)\} - v(i-1,2,0,0) \geq \left(\frac{1}{\mu_2}-\frac{1}{\mu_0}\right)\frac{(i+1)h_0}{2} + \tilde{c}_3, \label{eq:tilde-c3}
\end{align}where $\tilde{c}_3= \min\{\tilde{c}_2, \tilde{c}_2'\}$, also independent of $i$.}
\textcolor{blue}{Substituting the lower bounds in \eqref{eq3:i200-diff-fourth-term} and \eqref{eq:tilde-c3} into \eqref{eq3:i200-diff} yields
\begin{align}
    \lefteqn{D(i,2,0,0) = v(i, 1, 1, 0) - v(i, 1, 0, 1)}&\nonumber\\
    & \quad \geq \frac{(i+1) (\mu_2-\mu_1) h_0}{(\mu_0+\mu_2)( \mu_0+\mu_1)} + \frac{h_1}{\mu_0+\mu_1} - \frac{h_2}{\mu_0+\mu_2} \textcolor{blue}{+ \frac{\mu_0}{\mu_0+\mu_2} \left(\min\left\{\frac{h_1}{\mu_1} - \frac{h_2}{\mu_2},0\right\}\right)}\nonumber\\
    & \qquad + \frac{\mu_0\cdot(\mu_2-\mu_1)}{(\mu_0+\mu_2)( \mu_0+\mu_1)} \left(\textcolor{blue}{\left(\frac{1}{\mu_2}-\frac{1}{\mu_0}\right)\frac{(i+1)h_0}{2} + \tilde{c}_3}\right)\nonumber\\
    % &\quad = \frac{\mu_0\cdot(\mu_2-\mu_1)}{(\mu_0+\mu_2)( \mu_0+\mu_1)}\left(\frac{1}{\mu_2} + \frac{1}{\mu_0}\right)\frac{(i+1)h_0}{2} + \tilde{c}_4\nonumber\\
    & \quad =\frac{\mu_2-\mu_1}{\mu_2(\mu_0+\mu_1)}\frac{(i+1)h_0}{2} + \tilde{c}_4, \label{eq:tilde-c4}
\end{align}
where $\tilde{c}_4 = \frac{\mu_0\cdot(\mu_2-\mu_1)}{(\mu_0+\mu_2)( \mu_0+\mu_1)}\tilde{c}_3 + \frac{h_1}{\mu_0+\mu_1} - \frac{h_2}{\mu_0+\mu_2} \textcolor{blue}{+ \frac{\mu_0}{\mu_0+\mu_2} \left(\min\left\{\frac{h_1}{\mu_1} - \frac{h_2}{\mu_2},0\right\}\right)}$ is another constant independent of $i$.} 
\textcolor{blue}{The result follows since the right-hand side of \eqref{eq:tilde-c4} is non-negative for all sufficiently large $i$ when $\mu_1<\mu_2$.}
\end{proof}
% \endgroup
\begin{proposition}
    \label{prop:action1_large_i}
    Suppose \textcolor{blue}{$\frac{1}{\mu_1} \leq \frac{2}{\mu_2}$ and consider the decision states $(i,1,0,1) \in \X_{\widehat{D}}$. There exists a finite $i'$ such that $D(i,1,0,1) \leq 0$ for all $i \geq i'$.} If, in addition, $\frac{h_1}{\mu_1} \leq \frac{2h_2}{\mu_2}$ holds, then \textcolor{blue}{$D(i,1,0,1) \leq 0$} for all $i$.
\end{proposition}

\begin{proof}[Proof of Proposition \ref{prop:action1_large_i}] 
    \textcolor{blue}{We prove the result by considering two cases under $\frac{1}{\mu_1} \leq \frac{2}{\mu_2}$: $\mu_1 < \mu_2$ and $\mu_1 \geq \mu_2$.} 
\begin{enumerate}[label= \textbf{Case} \arabic*:, leftmargin=2.5\parindent]
    \item \textcolor{blue}{Consider first $\mu_1 < \mu_2$. 
    Notice that the second result is a direct application of \eqref{eq:decision_bd2} in Statement \ref{state4:general_rates} of Lemma \ref{lemma:general_rates}. In particular, $\frac{1}{\mu_2} < \frac{1}{\mu_1} \leq \frac{2}{\mu_2}$, together with $\frac{h_1}{\mu_1} \leq \frac{2h_2}{\mu_2}$, implies that $D(i,1,0,1) \leq 0$ for  \emph{all} $i$. It remains to prove the first result under} $\frac{1}{\mu_2} < \frac{1}{\mu_1} \leq \frac{2}{\mu_2}$ (and that $\frac{h_1}{\mu_1} > \frac{2h_2}{\mu_2}$). Assume without loss of generality $i \geq 2$. The proof proceeds by considering two subcases. 
\begin{enumerate}[label= \textbf{Subcase} (\alph*):, leftmargin=2\parindent]
    \item If $\frac{1}{\mu_1} < \frac{2}{\mu_2}$, the result holds by \eqref{eq:decision_bd2} in Statement \ref{state4:general_rates} in Lemma \ref{lemma:general_rates}. 
    \item If $\frac{1}{\mu_1} = \frac{2}{\mu_2}$, consider $v(i,0,1,1) - v(i,0,0,2)$. \textcolor{blue}{Recall that \eqref{eq:i101-diff} holds for $i \geq 1$ without assumptions on the service rates:}
    \begin{align}
        \lefteqn{D(i,1,0,1) = v(i,0,1,1) - v(i,0,0,2)} \nonumber\\
        & \quad = \ \frac{-i\mu_1h_0}{\mu_2(\mu_1+\mu_2)} + \textcolor{blue}{\frac{h_1+h_2}{\mu_1+\mu_2} - \frac{2h_2}{\mu_2}} \nonumber\\
        & \qquad +\frac{\mu_2}{\mu_1+\mu_2}\Big[v(i-1,1,1,0) - v(i-1,1,0,1)\Big]. \label{eq2:i101-diff}
    \end{align}
    Consider the generic difference $v(i,1,1,0) - v(i,1,0,1)$ (where we have replaced $i-1$ with $i$ in \eqref{eq2:i101-diff} and assume $i\geq 1$) \textcolor{blue}{and notice that \eqref{eq:i200-diff} continues to hold when $\mu_2 > \mu_1$:}
    \begin{small}
        \begin{align}
            \lefteqn{D(i,2,0,0) = v(i, 1, 1, 0) - v(i, 1, 0, 1)} \nonumber\\
                & \quad =\frac{(i+1) (\mu_2-\mu_1) h_0}{(\mu_0+\mu_2)( \mu_0+\mu_1)} + \frac{h_1}{\mu_0+\mu_1} - \frac{h_2}{\mu_0+\mu_2} \label{eq2:i200-diff}\\
                & \qquad + \ \frac{\mu_0}{\mu_0+\mu_2} \Big[\min\{v(i,0,2,0),v(i,0,1,1)\}-\min\{v(i,0,1,1),v(i,0,0,2)\}\Big] \nonumber\\
                & \qquad +   \frac{\mu_0\cdot(\mu_2-\mu_1)}{(\mu_0+\mu_2)( \mu_0+\mu_1)}\Big[\min\{v(i,0,2,0),v(i,0,1,1)\}-v(i-1,2,0,0)\Big]. \nonumber\end{align}  
    \end{small}%
    Consider the difference $\textcolor{blue}{\min\{v(i,0,2,0),v(i,0,1,1)\}}-\min\{v(i,0,1,1),v(i,0,0,2)\textcolor{blue}{\}}$ in the fourth term of \eqref{eq2:i200-diff}. Since \textcolor{blue}{$\mu_1 < \mu_2$ and $\frac{h_1}{\mu_1} > \frac{2h_2}{\mu_2} > \frac{h_2}{\mu_2}$, under which} inequality \eqref{eq:decision_lb_i110} in Proposition \ref{prop:decision_lb} implies that $v(i,0,2,0) \geq v(i,0,1,1)$, yielding
        \begin{align}
            \lefteqn{\textcolor{blue}{\min\{v(i,0,2,0),v(i,0,1,1)\}-\min\{v(i,0,1,1),v(i,0,0,2)\}}}& \nonumber\\
                &\quad = v(i,0,1,1) - \min\{v(i,0,1,1),v(i,0,0,2)\} \nonumber \\
                &\quad =\max\{0, v(i,0,1,1) - v(i,0,0,2)\} \nonumber\\
                & \quad \leq \max\left\{0,\frac{h_1}{\mu_1} - \frac{2h_2}{\mu_2} +\frac{ih_0}{2}\left(\frac{1}{\mu_1} - \frac{2}{\mu_2}\right)\right\}\nonumber\\
                & \quad = \textcolor{blue}{\frac{h_1}{\mu_1} - \frac{2h_2}{\mu_2}}, \label{eq2:i200-diff-third-term}
        \end{align}
        where the inequality holds by  \eqref{eq:decision_bd2} and the subsequent equality holds by recalling $\frac{1}{\mu_1} = \frac{2}{\mu_2}$ \textcolor{blue}{and $\frac{h_1}{\mu_1} > \frac{2h_2}{\mu_2}$.}

        \textcolor{blue}{Next, we consider the difference $\min\{v(i,0,2,0),v(i,0,1,1)\}-v(i-1,2,0,0)$ in the fifth term in \eqref{eq2:i200-diff}. 
        Since $\mu_2 = 2\mu_1 > \mu_1$, inequality \eqref{eq:decision_bd2_inte} in Statement \ref{state4:general_rates} of Lemma \ref{lemma:general_rates} holds. Replacing $i$ with $i-1$, this implies that the corresponding inequality also holds at level $i-1$:}
        \begin{align}
            v(i,0,1,1)-v(i-1,2,0,0) & \leq  \frac{h_1}{\mu_1} + \frac{(i+1)h_0}{2}\left(\frac{1}{\mu_1} - \frac{1}{\mu_0}\right).
            \label{eq2:i200-diff-fourth-term}
        \end{align}
        Using the bounds in \eqref{eq2:i200-diff-third-term} and \eqref{eq2:i200-diff-fourth-term} in \eqref{eq2:i200-diff} yields
        \begin{small}
        \begin{align}
            \lefteqn{D(i,2,0,0) = v(i, 1, 1, 0) - v(i, 1, 0, 1)} \nonumber\\
                &\quad \leq \frac{(i+1) (\mu_2-\mu_1) h_0}{(\mu_0+\mu_2)( \mu_0+\mu_1)} + \frac{h_1}{\mu_0+\mu_1} - \frac{h_2}{\mu_0+\mu_2} + \frac{\mu_0}{\mu_0+\mu_2} \left(\textcolor{blue}{\frac{h_1}{\mu_1} - \frac{2h_2}{\mu_2}}\right) \nonumber\\
                & \qquad +  \frac{\mu_0\cdot(\mu_2-\mu_1)}{(\mu_0+\mu_2)( \mu_0+\mu_1)}\left(\frac{h_1}{\mu_1} + \frac{(i+1)h_0}{2}\left(\frac{1}{\mu_1} - \frac{1}{\mu_0}\right)\right) \nonumber\\
                & \quad = \left(\frac{h_1}{\mu_1} - \frac{h_2}{\mu_2} + \frac{(i+1)h_0}{2}\left(\frac{1}{\mu_1} - \frac{1}{\mu_2}\right)\right) - \frac{\mu_0}{\mu_0+\mu_2}\left(\textcolor{blue}{\frac{h_2}{\mu_2}}+\frac{(i+1)h_0}{2\mu_2}\right), \label{eq2:i200-diff-last-term-bound}
        \end{align}
        \end{small}where the last step follows from $\frac{1}{\mu_1} = \frac{2}{\mu_2}$.
    Using the bound in \eqref{eq2:i200-diff-last-term-bound} (replacing $i$ with $i-1$) in \eqref{eq2:i101-diff} and again recalling $\frac{1}{\mu_1} = \frac{2}{\mu_2}$ yield 
    \begin{align}
        \lefteqn{D(i,1,0,1) = v(i,0,1,1) - v(i,0,0,2)} \nonumber \\
        & \quad \leq \frac{-i\mu_1h_0}{\mu_2(\mu_1+\mu_2)} + \textcolor{blue}{\frac{h_1+h_2}{\mu_1+\mu_2} - \frac{2h_2}{\mu_2}} \nonumber\\
        & \qquad + \frac{\mu_2}{\mu_1+\mu_2}\left(\frac{h_1}{\mu_1} - \frac{h_2}{\mu_2} + \frac{ih_0}{2}\left(\frac{1}{\mu_1} -\frac{1}{\mu_2}\right)- \frac{\mu_0}{\mu_0+\mu_2}\left(\textcolor{blue}{\frac{h_2}{\mu_2}}+\frac{ih_0}{2\mu_2}\right)\right) \nonumber \\
        & \quad = \frac{h_1}{\mu_1} - \frac{2h_2}{\mu_2} - \frac{\mu_2}{\mu_1+\mu_2}\frac{\mu_0}{\mu_0+\mu_2}\left(\textcolor{blue}{\frac{h_2}{\mu_2}}+\frac{ih_0}{2\mu_2}\right). \label{eq:i101-diff-final-bound}
    \end{align}
    This bound on the right-hand side of \eqref{eq:i101-diff-final-bound} is negative for $i$ large enough. The result follows.
    \end{enumerate}
    \item \textcolor{blue}{When $\mu_1 \geq \mu_2$, the first result follows directly from Statement \ref{state:larger_mu1_queue_in_Station2} of Proposition \ref{prop:thre_larger_mu1}. It therefore remains to establish the second result under the additional condition $\frac{h_1}{\mu_1} \leq \frac{2h_2}{\mu_2}$. 
    The derivations in the proof of Proposition \ref{prop:thre_larger_mu1} leading to \eqref{eq:i101-diff-final-bound2} continue to hold for $i \geq 1$ (recall $\mu_1 \geq \mu_2$):
    \begin{align}
        % \lefteqn{v(i,0,1,1) - v(i,0,0,2)} \nonumber\\
        % & \quad \leq 
        % \frac{-i\mu_1h_0}{\mu_2(\mu_1+\mu_2)} + \frac{h_1+h_2}{\mu_1+\mu_2} - \frac{2h_2}{\mu_2} +\frac{\mu_2}{\mu_1+\mu_2}\left(\frac{h_1}{\mu_1} - \frac{h_2}{\mu_2}\right) \nonumber \\
        % & \quad =
        D(i,1,0,1) = v(i,0,1,1) - v(i,0,0,2) \leq \left(\frac{h_1}{\mu_1} - \frac{2h_2}{\mu_2}\right) + \frac{-i\mu_1h_0}{\mu_2(\mu_1+\mu_2)}. \label{eq2:i101-diff-final-bound2}
    \end{align}
    When $i=0$, we have $v(0,0,1,1) - v(0,0,0,2) = \frac{h_1}{\mu_1} - \frac{2h_2}{\mu_2}$, which also satisfies \eqref{eq2:i101-diff-final-bound2}.
    It then follows that, if $\frac{h_1}{\mu_1} \leq \frac{2h_2}{\mu_2}$, the right-hand side of \eqref{eq2:i101-diff-final-bound2} is non-positive for all $i$.}
\end{enumerate}
    \textcolor{blue}{The proof completes by combining the cases $\mu_1 < \mu_2$ and $\mu_1\geq \mu_2$ under the condition $\frac{1}{\mu_1} \leq  \frac{2}{\mu_2}$.}
\end{proof}

\begin{proposition}
    \label{prop:queue_Station2_large_i}
    Suppose $\frac{1}{\mu_1} > \frac{2}{\mu_2}$. Consider the decision states $(i,1,0,1)$ where $i \geq 0$. 
    \begin{enumerate}
        \item \label{state:queue_Station2_action2_large_i}
        There exists a $\mu_0'$ small enough such that for any fixed $\mu_0 < \mu_0'$, there is an $i'$ sufficiently large \textcolor{blue}{with $D(i,1,0,1) \geq 0$} for all $i \geq i''$.
        \item \label{state:queue_Station2_action1_large_i}
        {\color{blue}There exists a $\mu_0'$ large enough such that for any fixed $\mu_0 > \mu_0'$, there is an $i'$ sufficiently large with $D(i,1,0,1) \leq 0$ for all $i \geq i'$.}
    \end{enumerate}
\end{proposition}

\begin{proof} [Proof of Statement \ref{state:queue_Station2_action2_large_i} in Proposition \ref{prop:queue_Station2_large_i}.]
The goal is to find $\mu_0'$ small enough so that for any 
$\mu_0 <\mu_0'$, the difference $v(i,0,1,1) - v(i,0,0,2) \geq 0$ for all $i$ sufficiently large. Since we are seeking $i$ large, without loss of generality, assume $i \geq 2$ so that there is enough work to keep the servers busy if they are assigned at Station 0. 
{\color{blue}The derivation in the proof of Statement \ref{state4:general_rates} of Lemma \ref{lemma:general_rates} leading to \eqref{eq:i101-diff} holds without assumptions on the service rates.} That is (recall we assumed $i \geq 2.$ If $i =0$ the transitions would be different),
\begin{align}
    \lefteqn{D(i,1,0,1) = v(i,0,1,1) - v(i,0,0,2)} \nonumber\\
    &\quad = \ \frac{-i\mu_1h_0}{\mu_2(\mu_1+\mu_2)} + \textcolor{blue}{\frac{h_1+h_2}{\mu_1+\mu_2} - \frac{2h_2}{\mu_2}} +\frac{\mu_2}{\mu_1+\mu_2}\Big(v(i-1,1,1,0) - v(i-1,1,0,1)\Big) \label{eq3:i101-diff}.
\end{align}
Consider the difference in fourth term in \eqref{eq3:i101-diff}, $v(i-1,1,1,0) - v(i-1,1,0,1)$. For ease of notation, replace $i-1$ with $i$ such that $i \geq 1$. 
\textcolor{blue}{Notice that the arithmetic from \eqref{eq:tilde-c4} in Proposition \ref{prop:empty-station2-large-i} continues to hold when $\mu_1 < \mu_2$. That is,
\begin{align*}
    D(i,2,0,0) = v(i, 1, 1, 0) - v(i, 1, 0, 1) \geq \frac{\mu_2-\mu_1}{\mu_2(\mu_0+\mu_1)}\frac{(i+1)h_0}{2} + \tilde{c}_4, 
\end{align*}
where $\tilde{c}_4 = \frac{\mu_0\cdot(\mu_2-\mu_1)}{(\mu_0+\mu_2)( \mu_0+\mu_1)}\tilde{c}_3 + \frac{h_1}{\mu_0+\mu_1} - \frac{h_2}{\mu_0+\mu_2} \textcolor{blue}{+ \frac{\mu_0}{\mu_0+\mu_2} \left(\min\left\{\frac{h_1}{\mu_1} - \frac{h_2}{\mu_2},0\right\}\right)}$.} 
Replacing $i$ in the inequality above with $i-1$ and a final substitution into \eqref{eq3:i101-diff} yields
\begin{align}
    \lefteqn{D(i,1,0,1) = v(i,0,1,1) - v(i,0,0,2)} \nonumber\\
    & \qquad \geq  \ \frac{-i\mu_1h_0}{\mu_2(\mu_1+\mu_2)} + \textcolor{blue}{\frac{h_1+h_2}{\mu_1+\mu_2} - \frac{2h_2}{\mu_2}} + \frac{\mu_2}{\mu_1+\mu_2}\left(\frac{\mu_2-\mu_1}{\mu_2(\mu_0+\mu_1)}\frac{ih_0}{2} + \tilde{c}_4\right) \nonumber\\
    & \qquad = \ \frac{ih_0}{\mu_1+\mu_2} \left(\frac{\mu_2-\mu_1}{2(\mu_0+\mu_1)} - \frac{\mu_1}{\mu_2}\right) + \tilde{c}_5, \label{eq:tilde-c5}
\end{align}
where $\tilde{c}_5 = \frac{\mu_2}{\mu_1+\mu_2}\tilde{c}_4 + \textcolor{blue}{\frac{h_1+h_2}{\mu_1+\mu_2} - \frac{2h_2}{\mu_2}}$ is a constant independent of $i$ as well.

Let $x:=\frac{\mu_2}{\mu_1}$ and $y:=\frac{\mu_0}{\mu_1}$ so that 
\begin{align*}
    & \frac{\mu_2-\mu_1}{2(\mu_0+\mu_1)} - \frac{\mu_1}{\mu_2} = \frac{x-1}{2(y+1)} - \frac{1}{x} = \frac{x^2-x-2-2y}{2x(y+1)}.
\end{align*}
When $\frac{1}{\mu_1}>\frac{2}{\mu_2}$, i.e., $\mu_2 > 2\mu_1$ we have $x>2$. Thus, $x^2-x-2 = (x-2)(x+1) > 0$. If $y < \frac{x^2-x-2}{2}$, or equivalently $\frac{\mu_0}{\mu_1} < \frac{1}{2}\big(\frac{\mu_2}{\mu_1}-2\big)\big(\frac{\mu_2}{\mu_1}+1\big)$, we have $\frac{\mu_2-\mu_1}{2(\mu_0+\mu_1)} - \frac{\mu_1}{\mu_2} > 0$. This can be done by choosing $\mu_0$ small enough, for example, by letting $\mu_0'=\frac{\mu_1}{2}\big(\frac{\mu_2}{\mu_1}-2\big)\big(\frac{\mu_2}{\mu_1}+1\big)$. The result follows from choosing $i$ large enough in \eqref{eq:tilde-c5}.
\end{proof}

{\color{blue}
\begin{proof} [Proof of Statement \ref{state:queue_Station2_action1_large_i} in Proposition \ref{prop:queue_Station2_large_i}.]
    We begin by defining two functions of $\mu_0$ below:
    \begin{align}
        g(\mu_0)&:= \frac{h_1}{\mu_1} - \frac{2h_2}{\mu_2} - \frac{\mu_2}{\mu_1+\mu_2}\frac{\mu_0}{\mu_0+\mu_2} \left(\frac{h_1}{\mu_1} - \frac{h_2}{\mu_2}\right), \label{eq:g-mu0}\\
        f(\mu_0) &:= \frac{1}{\mu_1} - \frac{2}{\mu_2} - \frac{\mu_2}{\mu_1+\mu_2}\frac{\mu_0}{\mu_0+\mu_2} \left(\frac{1}{\mu_1} - \frac{1}{\mu_2}\right). \label{eq:f-mu0}
    \end{align}
    Note that $f(\mu_0)$ is a continuous and (strictly) monotone decreasing function of $\mu_0 > 0$ with the following properties:
    \begin{align*}
        &\lim_{\mu_0 \to 0} f(\mu_0) = \frac{1}{\mu_1} - \frac{2}{\mu_2} > 0,\\
        &\lim_{\mu_0 \to \infty}f(\mu_0) = \frac{1}{\mu_1} - \frac{2}{\mu_2} - \frac{\mu_2}{\mu_1+\mu_2} \left(\frac{1}{\mu_1} - \frac{1}{\mu_2}\right) = 2\left(\frac{1}{\mu_1+\mu_2} - \frac{1}{\mu_2}\right) <0.
    \end{align*}
    It follows that there exists a unique $\mu_0'$ such that $f(\mu_0') = 0$. We now fix any $\mu_0 > \mu_0'$, which implies $f(\mu_0) <0$. Consequently, the following quantity $i_0$ is well-defined:
    \begin{align}
        i_0 := \min\left\{i\in Z_+\Bigg|i \geq 2, \frac{ih_0}{2}\cdot f(\mu_0) + g(\mu_0) < 0\right\}.
    \end{align} 
    For ease of exposition, we denote $c:=\frac{i_0h_0}{2}\cdot f(\mu_0) + g(\mu_0) < 0$.
    We, without loss of generality, consider $i\geq 2$. Recall that inequality \eqref{eq:i101-diff} follows directly from simple algebraic manipulation of the optimality equation \eqref{eq:opt1}:
    \begin{align}
        \lefteqn{D(i,1,0,1) = v(i,0,1,1) - v(i,0,0,2)} \nonumber\\
        &\quad = \ \frac{-i\mu_1h_0}{\mu_2(\mu_1+\mu_2)} + \frac{h_1+h_2}{\mu_1+\mu_2} - \frac{2h_2}{\mu_2} +\frac{\mu_2}{\mu_1+\mu_2}\Big(v(i-1,1,1,0) - v(i-1,1,0,1)\Big) \label{eq4:i101-diff}.
    \end{align}
    Again, for the difference $v(i-1,1,1,0) - v(i-1,1,0,1)$ appearing in the fourth term in \eqref{eq4:i101-diff}, inequality \eqref{eq:i200-diff} remains valid when $\mu_1 < \mu_2$ (with $i$ replaced by $i-1 \geq 1$):
    \begin{small}
    \begin{align}
        \lefteqn{D(i-1,2,0,0) = v(i-1, 1, 1, 0) - v(i-1, 1, 0, 1)}& \nonumber \\
        &\quad = \frac{i (\mu_2-\mu_1) h_0}{(\mu_0+\mu_2)( \mu_0+\mu_1)} + \frac{h_1}{\mu_0+\mu_1} - \frac{h_2}{\mu_0+\mu_2}  \nonumber \\
        &\qquad + \frac{\mu_0}{\mu_0+\mu_2} \Big[\min\{v(i-1,0,2,0),v(i-1,0,1,1)\}-\min\{v(i-1,0,1,1),v(i-1,0,0,2)\}\Big] \nonumber \\
        &\qquad +  \frac{\mu_0\cdot(\mu_2-\mu_1)}{(\mu_0+\mu_2)( \mu_0+\mu_1)}\Big[\min\{v(i-1,0,2,0),v(i-1,0,1,1)\}-v(i-2,2,0,0)\Big] \label{eq4:i200-diff}
    \end{align}
    \end{small}%
For the difference of two minima in the fourth term in \eqref{eq4:i200-diff}, choosing an upper bound $v(i-1,0,1,1)$ in the first minimum yields
\begin{align}
    \lefteqn {\min\{v(i-1,0,2,0),v(i-1,0,1,1)\}-\min\{v(i-1,0,1,1),v(i-1,0,0,2)\}}& \nonumber \\
    &\quad \leq v(i-1,0,1,1) - \min\{v(i-1,0,1,1),v(i-1,0,0,2)\} \nonumber \\
    &\quad = \max\{0, v(i-1,0,1,1) - v(i-1,0,0,2)\} \nonumber \\
    &\quad = \max\{0, D(i-1,1,0,1)\}. \label{eq4:i200-diff-third-term}
\end{align}
Recall that $\mu_2 > 2\mu_1 > \mu_1$, under which \eqref{eq:decision_bd2_inte} in Statement \ref{state4:general_rates} of Lemma \ref{lemma:general_rates} holds. Substituting \eqref{eq4:i200-diff-third-term} and \eqref{eq:decision_bd2_inte} (evaluated at $i-2$) into \eqref{eq4:i200-diff}, we obtain
\begin{align}
    \lefteqn{D(i-1,2,0,0) = v(i-1, 1, 1, 0) - v(i-1, 1, 0, 1)}& \nonumber \\
    &\quad \leq \frac{i (\mu_2-\mu_1) h_0}{(\mu_0+\mu_2)( \mu_0+\mu_1)} + \frac{h_1}{\mu_0+\mu_1} - \frac{h_2}{\mu_0+\mu_2} + \frac{\mu_0}{\mu_0+\mu_2} \cdot\max\{0, D(i-1,1,0,1)\} \nonumber \\
    &\qquad +  \frac{\mu_0\cdot(\mu_2-\mu_1)}{(\mu_0+\mu_2)( \mu_0+\mu_1)}\left(\frac{h_1}{\mu_1} + \frac{ih_0}{2}\left(\frac{1}{\mu_1} - \frac{1}{\mu_0}\right)\right) \nonumber \\
    &\quad = \frac{h_1}{\mu_1} - \frac{h_2}{\mu_2} + \frac{ih_0}{2}\left(\frac{1}{\mu_1} - \frac{1}{\mu_2}\right) \nonumber \\
    &\qquad +  \frac{\mu_0}{\mu_0+\mu_2} \left(\max\{0, D(i-1,1,0,1)\} - \left(\frac{h_1}{\mu_1} - \frac{h_2}{\mu_2}\right) - \frac{ih_0}{2}\left(\frac{1}{\mu_1} - \frac{1}{\mu_2}\right)\right). \label{eq4:i200-diff-ub}
\end{align}
Applying \eqref{eq4:i200-diff-ub} in \eqref{eq4:i101-diff}, we obtain that for all $i \geq 2$,
\begin{align*}
    \lefteqn{D(i,1,0,1) = v(i,0,1,1) - v(i,0,0,2)} \nonumber\\
    &\quad \leq \ \frac{-i\mu_1h_0}{\mu_2(\mu_1+\mu_2)} + \frac{h_1+h_2}{\mu_1+\mu_2} - \frac{2h_2}{\mu_2} +\frac{\mu_2}{\mu_1+\mu_2}\left(\frac{h_1}{\mu_1} - \frac{h_2}{\mu_2} + \frac{ih_0}{2}\left(\frac{1}{\mu_1} - \frac{1}{\mu_2}\right)\right) \\
    &\qquad + \frac{\mu_2}{\mu_1+\mu_2} \left(\frac{\mu_0}{\mu_0+\mu_2} \left(\max\{0, D(i-1,1,0,1)\} - \left(\frac{h_1}{\mu_1} - \frac{h_2}{\mu_2}\right) - \frac{ih_0}{2}\left(\frac{1}{\mu_1} - \frac{1}{\mu_2}\right)\right)\right)\\
    &\quad = \frac{h_1}{\mu_1} - \frac{2h_2}{\mu_2} + \frac{ih_0}{2}\left(\frac{1}{\mu_1} - \frac{2}{\mu_2}\right) \\
    &\qquad + \frac{\mu_2}{\mu_1+\mu_2}\frac{\mu_0}{\mu_0+\mu_2} \left(\max\{0, D(i-1,1,0,1)\} - \left(\frac{h_1}{\mu_1} - \frac{h_2}{\mu_2}\right) - \frac{ih_0}{2}\left(\frac{1}{\mu_1} - \frac{1}{\mu_2}\right)\right)\\
    &\quad = \frac{\mu_2}{\mu_1+\mu_2}\frac{\mu_0}{\mu_0+\mu_2} \max\{0, D(i-1,1,0,1)\} + \frac{ih_0}{2}\cdot f(\mu_0) + g(\mu_0),
\end{align*}where the last step refers to the definitions of $f(\mu_0)$ and $g(\mu_0)$ given in \eqref{eq:f-mu0} and \eqref{eq:g-mu0}, respectively.

Now consider $i \geq i_0$. Since $f(\mu_0)<0$, it follows that $\frac{ih_0}{2}\cdot f(\mu_0) + g(\mu_0) \leq \frac{i_0h_0}{2}\cdot f(\mu_0) + g(\mu_0) = c < 0$. Therefore, for all $i \geq i_0$, we have 
\begin{align}
    D(i,1,0,1) &\leq \frac{\mu_2}{\mu_1+\mu_2}\frac{\mu_0}{\mu_0+\mu_2} \max\{0, D(i-1,1,0,1)\} + c \nonumber \\
    &\leq \max\{0, D(i-1,1,0,1)\} + c, \label{eq4:i200-diff-ub-final}
\end{align}where the last step follows since $\max\{0, D(i-1,1,0,1)\} \geq 0$.
Suppose towards a contradiction that $D(i,1,0,1) > 0$ for all $i\geq i_0$. Then for any $i \geq i_0+1$, inequality \eqref{eq4:i200-diff-ub-final} yields
\begin{align*}
        D(i,1,0,1) \leq D(i-1,1,0,1) + c,
\end{align*}
leading to 
\begin{align*}
        D(i,1,0,1) & \leq D(i_0,1,0,1) + c (i-i_0).
\end{align*}
This yields a contradiction, as the right-hand side becomes negative for sufficiently large $i$ (recall that $c<0$). Therefore, there exists an $i' \geq i_0$ such that  $D(i',1,0,1) \leq 0$. It then follows from \eqref{eq4:i200-diff-ub-final} that $D(i,1,0,1) = v(i,0,1,1) - v(i,0,0,2) \leq c < 0$ for all $i\geq i'$, which completes the proof.
\end{proof}}

\subsection{Preliminaries for Theorem \ref{thm:threshold}}
\begin{lemma}\label{lemma:equal_rates}
    Suppose $\mu_1 = \mu_2 := \mu$, the following inequalities hold.
    \begin{equation}
        v(0, 0, 1, 0) - v(0, 0, 0, 1) - v(0, 1, 1, 0) + v(0, 1, 0, 1) \geq 0. \label{eq:equal_rate_bdry}
    \end{equation}
    Suppose \textcolor{blue}{further that $h_1 \geq h_2$}. For all $i \geq 0$, we have
    \begin{align}
    v(i, 1, 1, 0) - v(i, 1, 0, 1)
        - v(i, 0, 1, 1) + v(i, 0, 0, 2) - \frac{2h_2+ih_0}{\mu} \leq 0. \label{eq:equal_rates_inte}
\end{align}
\end{lemma}
\begin{proof} [Proof of Lemma \ref{lemma:equal_rates}] Consider first \eqref{eq:equal_rate_bdry}. The (in)equalities of values that have been calculated in proving Statement \ref{state1:general_rates} yield
\begin{align*}
    \lefteqn{v(0, 0, 1, 0) - v(0, 0, 0, 1) - v(0, 1, 1, 0) + v(0, 1, 0, 1)}&\\
    &\quad \geq \frac{h_1}{\mu} - \frac{h_2}{\mu} - \left(\frac{h_0}{\mu_0} + \frac{h_1}{\mu} + \min\left\{\frac{h_1}{\mu}, \frac{h_2}{\mu}\right\}\right) + \left(\frac{h_0}{\mu_0} + \frac{h_2}{\mu} + \min\left\{\frac{h_1}{\mu}, \frac{h_2}{\mu}\right\}\right)\\
    &\quad = 0.
\end{align*}
Proof by induction on $i$ for \eqref{eq:equal_rates_inte}. When $i = 0$, recall the (in)equalities of the values computed in the proof of Statement \ref{state1:general_rates} in Lemma \ref{lemma:general_rates}, we have
    \begin{small}
    \begin{align*}
        \lefteqn{v(0, 1, 1, 0) - v(0, 1, 0, 1) - v(0, 0, 1, 1) + v(0, 0, 0, 2) - \frac{2h_2}{\mu}}&\\ 
        &\quad \leq \left(\frac{h_0}{\mu_0} + \frac{h_1}{\mu} + \min\left\{\frac{h_1}{\mu}, \frac{h_2}{\mu}\right\}\right) - \left(\frac{h_0}{\mu_0} + \frac{h_2}{\mu} + \min\left\{\frac{h_1}{\mu}, \frac{h_2}{\mu}\right\}\right) - \left(\frac{h_1}{\mu} + \frac{h_2}{\mu}\right) + \frac{3h_2}{\mu} - \frac{2h_2}{\mu}\\
        &\quad \leq -\frac{h_2}{\mu} \leq 0.
    \end{align*}
    \end{small}Suppose the inequality holds at $i-1$, where $i \geq 1$, consider it at $i$. 
    \textcolor{blue}{We begin by presenting the expansions of the differences $D(i,1,1,0)$, $D(i,2,0,0)$ and $D(i,1,0,1)$ for later references. For decision states $(i,1,1,0)$, evaluating the optimality equation \eqref{eq:opt1} at states $(i,0,2,0)$ and $(i,0,1,1)$ yields
    \begin{small}
    \begin{align}
        \lefteqn{D(i,1,1,0) = v(i, 0, 2, 0) - v(i, 0, 1, 1)} \nonumber \\
        & \quad = \left(\frac{ih_0+2h_1}{2\mu} + v(i-1, 1, 1, 0)\right) - \left(\frac{ih_0+h_1+h_2}{2\mu} + \frac{1}{2} v(i-1, 1, 0, 1) + \frac{1}{2} v(i-1, 1, 1, 0)\right) \nonumber\\
        & \quad = \frac{h_1 - h_2}{2\mu} + \frac{1}{2}\Big[v(i-1, 1, 1, 0) - v(i-1, 1, 0, 1)\Big]. \label{eq:i110-diff-equal}
    \end{align}
    \end{small} For decision states $(i,2,0,0)$,}
    \textcolor{blue}{
    \begin{small}
    \begin{align}
        \lefteqn{D(i,2,0,0) = v(i, 1, 1, 0) - v(i, 1, 0, 1)}& \nonumber \\
        &\quad = \left(\frac{(i+1)h_0+h_1}{\mu_0+\mu} + \frac{\mu_0}{\mu_0+\mu} \min\{v(i, 0, 2, 0),v(i, 0, 1, 1)\} + \frac{\mu}{\mu_0+\mu} v(i-1, 2, 0, 0)\right) \nonumber \\
        &\qquad - \left(\frac{(i+1)h_0+h_2}{\mu_0+\mu} + \frac{\mu_0}{\mu_0+\mu} \min\{v(i, 0, 1, 1),v(i, 0, 0, 2)\} + \frac{\mu}{\mu_0+\mu} v(i-1, 2, 0, 0)\right)\nonumber \\
        & \quad = \frac{h_1-h_2}{\mu_0+\mu} + \frac{\mu_0}{\mu_0+\mu} \Big[ \min\{v(i, 0, 2, 0),v(i, 0, 1, 1)\} - \min\{v(i, 0, 1, 1),v(i, 0, 0, 2)\} \Big]. \label{eq:i200-diff-equal} 
    \end{align}
    \end{small}Notice that when $i=0$, although the transitions differs (with $(i-1,2,0,0)$ replaced by $(0,1,0,0)$), equation \eqref{eq:i200-diff-equal} still holds.} 
    \textcolor{blue}{Similarly, the difference $D(i,1,0,1)$, where $i\geq 1$, can be expressed as follows:
    \begin{small}
    \begin{align}
        \lefteqn{D(i,1,0,1) = v(i, 0, 1, 1) - v(i, 0, 0, 2)}& \nonumber \\
        &\quad = \left(\frac{ih_0+h_1+h_2}{2\mu} + \frac{1}{2} v(i-1, 1, 1, 0) + \frac{1}{2} v(i-1, 1, 0, 1)\right) - \left(\frac{ih_0+2h_2}{\mu} + v(i-1, 1, 0, 1)\right) \nonumber \\
        & \quad = \frac{h_1 - 3h_2 - ih_0}{2\mu} + \frac{1}{2} \Big[ v(i-1, 1, 1, 0) - v(i-1, 1, 0, 1)\Big]. \label{eq:i101-diff-equal}
    \end{align}
    \end{small}%
    Substituting \eqref{eq:i200-diff-equal} and \eqref{eq:i101-diff-equal} into the left-hand side of \eqref{eq:equal_rates_inte} yields}
    \begin{align}
        \lefteqn{v(i, 1, 1, 0) - v(i, 1, 0, 1) - v(i, 0, 1, 1) + v(i, 0, 0, 2) - \frac{2h_2+ih_0}{\mu}}& \nonumber \\
        &\quad = \frac{h_1-h_2}{\mu_0+\mu} - \frac{h_1 - 3h_2 - ih_0}{2\mu} - \frac{2h_2+ih_0}{\mu} \nonumber \\
        &\qquad + \frac{\mu_0}{\mu_0+\mu} \Big[ \min\{v(i, 0, 2, 0),v(i, 0, 1, 1)\} - \min\{v(i, 0, 1, 1),v(i, 0, 0, 2)\} \Big] \nonumber \\
        &\qquad - \frac{1}{2} \Big[ v(i-1, 1, 1, 0) - v(i-1, 1, 0, 1)\Big]. \label{eq:eq:equal_rates_inte1}
    \end{align}
There are two cases to discuss based on whether $v(i, 0, 1, 1) \leq v(i, 0, 0, 2)$.
\begin{enumerate}[label= \textbf{Case} \arabic*:, leftmargin=2.5\parindent]
    \item If $v(i, 0, 1, 1) \leq v(i, 0, 0, 2)$, replacing the first minimum with an upper bound $(i,0,2,0)$ yields
    \begin{small}
    \begin{align*}
        \lefteqn{v(i, 1, 1, 0) - v(i, 1, 0, 1) - v(i, 0, 1, 1) + v(i, 0, 0, 2) - \frac{2h_2+ih_0}{\mu}}& \\
        &\quad \leq  \frac{h_1-h_2}{\mu_0+\mu} - \frac{h_1 - 3h_2 - ih_0}{2\mu}  - \frac{2h_2+ih_0}{\mu}\\
        &\qquad + \frac{\mu_0}{\mu_0+\mu} \Big[  v(i, 0, 2, 0) - v(i, 0, 1, 1)\Big] - \frac{1}{2} \Big[v(i-1, 1, 1, 0) - v(i-1, 1, 0, 1)\Big]\\
        &\quad = - \frac{2h_2+ih_0}{\mu} + \textcolor{black}{\frac{2h_2+ih_0}{2\mu}} + \frac{\mu_0}{2(\mu_0+\mu)} \Big[v(i-1, 1, 1, 0) - v(i-1, 1, 0, 1) \\
        &\qquad - \min\{v(i-1, 0, 2, 0),v(i-1, 0, 1, 1)\}  + \min\{v(i-1, 0, 1, 1),v(i-1, 0, 0, 2)\}\Big],
    \end{align*}
    \end{small}\textcolor{blue}{where the last step follows from substituting \eqref{eq:i110-diff-equal} and \eqref{eq:i200-diff-equal}.}
    \item If $v(i, 0, 1, 1) > v(i, 0, 0, 2)$, choosing an upper bound $(i,0,1,1)$ in the first minimum in \eqref{eq:eq:equal_rates_inte1} yields
    \begin{small}
    \begin{align*}
        \lefteqn{v(i, 1, 1, 0) - v(i, 1, 0, 1) - v(i, 0, 1, 1) + v(i, 0, 0, 2) - \frac{2h_2+ih_0}{\mu}}& \\
        &\quad \leq  \frac{h_1-h_2}{\mu_0+\mu} - \frac{h_1 - 3h_2 - ih_0}{2\mu}  - \frac{2h_2+ih_0}{\mu} \\
        &\qquad + \frac{\mu_0}{\mu_0+\mu} \Big[ v(i, 0, 1, 1) - v(i, 0, 0, 2)\Big] - \frac{1}{2}\Big[v(i-1, 1, 1, 0) - v(i-1, 1, 0, 1)\Big]\\
        &\quad = - \frac{2h_2+ih_0}{\mu} + \textcolor{black}{\frac{2h_2+ih_0}{2(\mu_0+\mu)}} + \frac{\mu_0}{2(\mu_0+\mu)} \Big[v(i-1, 1, 1, 0) - v(i-1, 1, 0, 1) \\
        & \qquad - \min\{v(i-1, 0, 2, 0),v(i-1, 0, 1, 1)\} + \min\{v(i-1, 0, 1, 1),v(i-1, 0, 0, 2)\}\Big],
    \end{align*}
    \end{small}\textcolor{blue}{where the final step is obtained by applying \eqref{eq:i101-diff-equal} and \eqref{eq:i200-diff-equal}.}
\end{enumerate}
\textcolor{blue}{In either case, the following holds:}
\begin{align*}
    \lefteqn{v(i, 1, 1, 0) - v(i, 1, 0, 1) - v(i, 0, 1, 1) + v(i, 0, 0, 2) - \frac{2h_2+ih_0}{\mu}}& \\
    &\quad \leq - \frac{2h_2+ih_0}{\mu} + \textcolor{black}{\frac{2h_2+ih_0}{2\mu}} + \frac{\mu_0}{2(\mu_0+\mu)} \Big[v(i-1, 1, 1, 0) - v(i-1, 1, 0, 1) \\
    & \qquad - \min\{v(i-1, 0, 2, 0),v(i-1, 0, 1, 1)\} + \min\{v(i-1, 0, 1, 1),v(i-1, 0, 0, 2)\}\Big].
\end{align*}
\textcolor{blue}{Since $h_1 \geq h_2$ and $\mu_1 = \mu_2 = \mu$, inequality \eqref{eq:decision_lb_i110} in Proposition \ref{prop:decision_lb} implies that $v(i-1, 0, 2, 0) \geq v(i-1, 0, 1, 1)$ in the first minimum.} 
Replacing the second minimum with an upper bound $v(i-1, 0, 0, 2)$ yields
\begin{align*}
    \lefteqn{v(i, 1, 1, 0) - v(i, 1, 0, 1) - v(i, 0, 1, 1) + v(i, 0, 0, 2) - \frac{2h_2+ih_0}{\mu}}& \\
    &\quad \leq - \frac{2h_2+ih_0}{\mu} + \textcolor{black}{\frac{2h_2+ih_0}{2\mu}} \\
    &\qquad + \frac{\mu_0}{2(\mu_0+\mu)} \Big[v(i-1, 1, 1, 0) - v(i-1, 1, 0, 1) - v(i-1, 0, 1, 1)+ v(i-1, 0, 0, 2)\Big]\\
    &\quad \leq - \frac{2h_2+ih_0}{\mu} + \textcolor{black}{\frac{2h_2+ih_0}{2\mu}} + \frac{\mu_0}{2(\mu_0+\mu)} \left(\frac{2h_2+(i-1)h_0}{\mu}\right)\\
    &\quad \leq 0,
\end{align*}
where the second-to-last inequality follows from the inductive hypothesis.
\end{proof}

\subsection{Proof of Theorem \ref{thm:threshold}}\label{sec:equal-service-rates}
In this section, we prove Theorem \ref{thm:threshold}, which characterizes the threshold structure of \textcolor{blue}{an optimal policy} under the additional assumption that $\mu_1 = \mu_2 := \mu$.
% so that Assumption \ref{assm:cost-to-serve} becomes just $h_1 \geq h_2$. 
Collectively, in the proofs we refer to these threshold results as the ``monotonicity'' results. 

% \begin{proof} [Proof of Statements \ref{thm:equal_rate1} and \ref{thm:equal_rate2} in Theorem \ref{thm:threshold}]
%     The first result is implied by Proposition \ref{prop:decision_lb}. The second follows similarly as the proof for the state $(i,1,0,1)$ in Proposition \ref{prop:thre_larger_mu1}. We have omitted the proof for brevity.  
% \end{proof}

\begin{proof} [Proof of Statement \ref{state:threshold_i} in Theorem \ref{thm:threshold} (Monotonicity in $i$)]
To show the policy is monotone in $i$, we need to show for all $i\geq0$,
\begin{align*}
    v(i, j-1, k+1, \ell) \leq v(i, j-1, k, \ell+1)
        \Rightarrow v(i+1, j-1, k+1, \ell) \leq v(i+1, j-1, k, \ell+1).
\end{align*}
This is implied by the following inequalities:
\begin{align}
    \lefteqn{D(i,1,1,0) - D(i+1,1, 1, 0) } \nonumber \\
        & \quad = v(i, 0,2, 0) - v(i, 0, 1, 1) - \big[v(i+1, 0,2, 0) - v(i+1, 0, 1, 1)\big] \geq 0, \label{eq1:mono_i}
\end{align}
\begin{align}
    \lefteqn{D(i, 1,0, 1) - v(i+1, 1, 0, 1) } \nonumber \\
        & \qquad = v(i, 0,1, 1) - v(i, 0, 0, 2) - \big[v(i+1, 0,1, 1) - v(i+1, 0, 0, 2)\big] \geq 0, \label{eq2:mono_i}
\end{align}
\begin{align}
    \lefteqn{D(i, 2, 0, 0) - v(i+1, 2, 0, 0) }  \nonumber\\
        & \qquad = v(i, 1, 1, 0) - v(i, 1, 0, 1) - \big[v(i+1, 1,1, 0) - v(i+1, 1, 0, 1)\big] \geq 0, \label{eq3:mono_i}
\end{align}
for all $i \geq 0$. Proof by induction on $i$.
Suppose $i = 0$. For \eqref{eq1:mono_i}, rearranging optimality equations \eqref{eq:opt-zero}--\eqref{eq:opt1} yields
\begin{align*}
    \lefteqn{v(0, 0, 2, 0) - v(0, 0, 1, 1) - v(1, 0, 2, 0) + v(1, 0, 1, 1)}\\
     & \qquad = \frac{2h_1}{2\mu} - \frac{h_1+h_2}{2\mu}- \frac{h_0+2h_1}{2\mu} + \frac{h_0+h_1+h_2}{2\mu}\\
     & \qquad \quad + \frac{1}{2} \Big[ v(0, 0, 1, 0) - v(0, 0, 0, 1) - v(0, 1, 1, 0) + v(0, 1, 0, 1)\Big]\\
     & \qquad \quad + \frac{1}{2}\Big[ v(0, 0, 1, 0) - v(0, 0, 1, 0) - v(0, 1, 1, 0) + v(0, 1, 1, 0)\Big]\\
     & \qquad = \frac{1}{2}\Big[ v(0, 0, 1, 0) - v(0, 0, 0, 1) - v(0, 1, 1, 0) + v(0, 1, 0, 1)\Big]\\
     & \qquad \geq  0,
\end{align*} 
where the last inequality holds by \eqref{eq:equal_rate_bdry} in Lemma \ref{lemma:equal_rates}. Inequality \eqref{eq2:mono_i} follows similarly and is omitted for brevity.

%Consider \eqref{eq2:mono_i}.
%\begin{align*}
%    \lefteqn{v(0, 0, 1, 1) - v(0, 0, 0, 2) - v(1, 0, 1, 1) + v(1, 0, 0, 2)}\\
%    & \qquad = \frac{h_1+h_2}{2\mu} - \frac{2h_2}{\mu} - \frac{h_0+h_1+h_2}{2\mu} + \frac{h_0+2h_2}{\mu} \\
%    & \qquad \quad + \Big[\frac{1}{2} v(0, 0, 0, 1) + \frac{1}{2} v(0, 0, 1, 0)\Big] - v(0, 0, 0, 1) \\
%    & \qquad \quad - \Big[\frac{1}{2} v(0, 1, 0, 1) + \frac{1}{2} v(0, 1, 1, 0)\Big] + v(0, 1, 0, 1) \\
%    & \qquad \quad =  \frac{h_0}{\mu} - \frac{h_0}{2\mu} + \frac{1}{2} \Big[v(0, 0, 1, 0) - v(0, 0, 0, 1) - v(0, 1, 1, 0) + v(0, 1, 0, 1)\Big]\\
%    & \qquad \quad \geq 0,
%\end{align*}
%where the inequality follows from \eqref{eq:boundary}.

Consider \eqref{eq3:mono_i}.
\begin{align*}
    \lefteqn{v(0, 1, 1, 0) - v(0, 1, 0, 1) - v(1, 1, 1, 0) + v(1, 1, 0, 1)}\\
    & \qquad = \frac{h_0+h_1}{\mu_0+\mu} - \frac{h_0+h_2}{\mu_0 + \mu} - \frac{2h_0+h_1}{\mu_0 + \mu} + \frac{2h_0+h_2}{\mu_0 + \mu}\\
    & \qquad \quad+ \frac{\mu}{\mu_0 + \mu} \Big[ v(0, 1, 0, 0) - v(0, 1, 0, 0) - v(0, 2, 0, 0) + v(0, 2, 0, 0)\Big]\\
    & \qquad \quad + \frac{\mu_0}{\mu_0 + \mu}  \Big[ \min \{v(0, 0, 2, 0),v(0, 0, 1, 1)\} - \min \{v(0, 0, 1, 1),v(0, 0, 0, 2)\}\\
    & \qquad \quad - \min \{v(1, 0, 2, 0),v(1, 0, 1, 1)\} + \min \{v(1, 0, 1, 1),v(1, 0, 0, 2)\}\Big]\\
    & \qquad = \frac{\mu_0}{\mu_0 + \mu}  \Big[ \min \{v(0, 0, 2, 0),v(0, 0, 1, 1)\} - \min \{v(0, 0, 1, 1),v(0, 0, 0, 2)\}\\
    & \qquad \quad - \min \{v(1, 0, 2, 0),v(1, 0, 1, 1)\} + \min \{v(1, 0, 1, 1),v(1, 0, 0, 2)\}\Big].
\end{align*}
There are several cases to consider, 
\begin{enumerate}[label= \textbf{Case} \arabic*:, leftmargin=3\parindent]
    \item Suppose $v(0, 0, 2, 0) \leq v(0, 0, 1, 1)$. Consider the subcases
    \begin{enumerate}[label= \textbf{Subcase} (\alph*):, leftmargin=3\parindent]
        \item If $v(1, 0, 1, 1) \leq v(1, 0, 0, 2)$, choosing $v(0,0,1,1)$ and $v(1,0,2,0)$ respectively in the second and third minima yields
        \begin{small}
        \begin{align*}
            \lefteqn{v(0, 1, 1, 0) - v(0, 1, 0, 1) - v(1, 1, 1, 0) + v(1, 1, 0, 1)} \\
                & \qquad \geq \frac{\mu_0}{\mu_0 + \mu} \Big[ v(0, 0, 2, 0) - v(0, 0, 1, 1) - v(1, 0, 2, 0) + v(1, 0, 1, 1) \Big] \geq 0,
        \end{align*} 
        \end{small}by \eqref{eq1:mono_i} at $i=0$.
        \item If $v(1, 0, 1, 1) \geq v(1, 0, 0, 2)$, choosing $v(0,0,0,2)$ and $v(1,0,2,0)$ respectively in the second and third minima yields
        \begin{small}
        \begin{align*}
            \lefteqn{v(0, 1, 1, 0) - v(0, 1, 0, 1) - v(1, 1, 1, 0) + v(1, 1, 0, 1)}\\
                & \qquad \geq \frac{\mu_0}{\mu_0 + \mu} \Big[ v(0, 0, 2, 0) - v(0, 0, 0, 2) - v(1, 0, 2, 0) + v(1, 0, 0, 2) \Big]\geq 0,
        \end{align*} 
        \end{small}where the inequality holds by adding \eqref{eq1:mono_i} and \eqref{eq2:mono_i} for $i=0$ and canceling terms.
    \end{enumerate}
    \item Suppose $v(0, 0, 2, 0) \geq v(0, 0, 1, 1)$. Consider the subcases
    \begin{enumerate}[label= \textbf{Subcase} (\alph*):, leftmargin=3\parindent]
        \item If $v(1, 0, 1, 1) \leq v(1, 0, 0, 2)$, by choosing $v(0,0,1,1)$ and $v(1,0,1,1)$ in the second and third minima, respectively,
        \begin{small}
        \begin{align*}
            \lefteqn{v(0, 1, 1, 0) - v(0, 1, 0, 1) - v(1, 1, 1, 0) + v(1, 1, 0, 1)} \\
                & \qquad \geq \frac{\mu_0}{\mu_0 + \mu} \Big[ v(0, 0, 1, 1) - v(0, 0, 1, 1) - v(1, 0, 1, 1) + v(1, 0, 1, 1) \Big] = 0.
        \end{align*}
        \end{small}
         \item If $v(1, 0, 1, 1) \geq v(1, 0, 0, 2)$, by choosing $v(0,0,0,2)$ and $v(1,0,1,1)$ in the second and third minima, respectively,
         \begin{small}
        \begin{align*}
            \lefteqn{v(0, 1, 1, 0) - v(0, 1, 0, 1) - v(1, 1, 1, 0) + v(1, 1, 0, 1)}\\
                & \qquad \geq \frac{\mu_0}{\mu_0 + \mu} \Big[ v(0, 0, 1, 1) - v(0, 0, 0, 2) - v(1, 0, 1, 1) + v(1, 0, 0, 2) \Big] \geq 0,
        \end{align*} 
        \end{small}where the inequality holds by \eqref{eq2:mono_i} at $i=0$.
    \end{enumerate}
\end{enumerate}

Suppose the inequalities \eqref{eq1:mono_i} - \eqref{eq3:mono_i} hold at $i-1$, then starting from \eqref{eq1:mono_i} we have (for $i \geq 1$),
\begin{align*}
    \lefteqn{v(i, 0, 2, 0) - v(i, 0, 1, 1) - v(i+1, 0, 2, 0) + v(i+1, 0, 1, 1)}\\
    & \qquad =  \frac{ih_0+2h_1}{2\mu} - \frac{ih_0+h_1+h_2}{2\mu} - \frac{(i+1)h_0+2h_1}{2\mu} + \frac{(i+1)h_0+h_1+h_2}{2\mu}\\
    & \qquad \quad+ \frac{1}{2} \Big[ v(i-1, 1, 1, 0) - v(i-1, 1, 0, 1) - v(i, 1, 1, 0) + v(i, 1, 0, 1)\Big]\\
    & \qquad \quad + \frac{1}{2} \Big[ v(i-1, 1, 1, 0) - v(i-1, 1, 1, 0) - v(i, 1, 1, 0) + v(i, 1, 1, 0)\Big]\\
    & \qquad \geq  0,
\end{align*}
where the inequality follows from the inductive assumption of \eqref{eq3:mono_i}. For \eqref{eq2:mono_i} at $i$, note
\begin{align*}
    \lefteqn{v(i, 0, 1, 1) - v(i, 0, 0, 2) - v(i+1, 0, 1, 1) + v(i+1, 0, 0, 2)}\\
        & \qquad = \frac{ih_0+h_1+h_2}{2\mu} - \frac{ih_0+2h_2}{\mu} - \frac{(i+1)h_0+h_1+h_2}{2\mu} + \frac{(i+1)h_0+2h_2}{\mu} \\
        & \qquad \quad + \Big[\frac{1}{2} v(i-1, 1, 0, 1) + \frac{1}{2} v(i-1, 1, 1, 0)\Big] - v(i-1, 1, 0, 1) \\
        & \qquad \quad - \Big[\frac{1}{2} v(i, 1, 0, 1) + \frac{1}{2} v(i, 1, 1, 0)\Big] + v(i, 1, 0, 1) \\
        & \qquad = \frac{h_0}{\mu} - \frac{h_0}{2\mu} + \frac{1}{2} \Big[v(i-1, 1, 1, 0) - v(i-1, 1, 0, 1) - v(i, 1, 1, 0) + v(i, 1, 0, 1)\Big]\\
        & \qquad \geq  0,
\end{align*}
where the inequality follows from the inductive assumption of \eqref{eq3:mono_i}. For \eqref{eq3:mono_i} at $i$,
\begin{align*}
    \lefteqn{v(i, 1, 1, 0) - v(i, 1, 0, 1) - v(i+1, 1, 1, 0) + v(i+1, 1, 0, 1)}\\
        & \qquad = \frac{(i+1)h_0+h_1}{\mu_0+\mu} - \frac{(i+1)h_0+h_2}{\mu_0+\mu} - \frac{(i+2)h_0+h_1}{\mu_0+\mu} + \frac{(i+2)h_0+h_2}{\mu_0+\mu}\\
        & \qquad \quad + \frac{\mu}{\mu_0+\mu} \Big[ v(i-1, 2, 0, 0) - v(i-1, 2, 0, 0) - v(i, 2, 0, 0) + v(i, 2, 0, 0)\Big]\\
        & \qquad \quad + \frac{\mu_0}{\mu_0+\mu} \Big[ \min \{v(i, 0, 2, 0),v(i, 0, 1, 1)\} - \min \{v(i, 0, 1, 1),v_n(i, 0, 0, 2)\}\\
        & \qquad \quad  - \min \{v(i+1, 0, 2, 0),v(i+1, 0, 1, 1)\} + \min \{v(i+1, 0, 1, 1),v(i+1, 0, 0, 2)\}\Big]\\
        & \qquad = \frac{\mu_0}{\mu_0+\mu} \Big[ \min \{v(i, 0, 2, 0),v(i, 0, 1, 1)\} - \min \{v(i, 0, 1, 1),v(i, 0, 0, 2)\}\\
        & \qquad \quad - \min \{v(i+1, 0, 2, 0),v(i+1, 0, 1, 1)\} + \min \{v(i+1, 0, 1, 1),v(i+1, 0, 0, 2)\}\Big] \geq 0,
\end{align*}
where the discussion of the minima is similar to what we proved when $i = 0$ and is omitted here for brevity.
\end{proof}
\begin{proof} [Proof of Statement \ref{state:threshold12} of Monotonicity in Theorem \ref{thm:threshold}]
Such \textcolor{blue}{a }threshold exists if for all $i \geq 0$,
\begin{align}
    \lefteqn{D(i, 1, 1, 0) - D(i, 1, 0, 1) }\nonumber \\
    &\quad = v(i, 0, 2, 0) - v(i, 0, 1, 1)  - v(i, 0, 1, 1) + v(i, 0, 0, 2) \nonumber \\
    &\quad =  v(i, 0, 2, 0) - 2v(i, 0, 1, 1)  + v(i, 0, 0, 2) \geq 0. \label{eq:mono12}
\end{align}
This result holds by noting that for all $i \geq 1$,
\begin{align*}
    \lefteqn{v(i, 0, 2, 0) - 2v(i, 0, 1, 1) + v(i, 0, 0, 2)}\\
        & \qquad = \frac{ih_0+2h_1}{2\mu} - 2\frac{ih_0+h_1+h_2}{2\mu} + \frac{ih_0+2h_2}{\mu}\\
        & \qquad \quad + v(i-1, 1, 1, 0) - 2 \Big[ \frac{1}{2} v(i-1, 1, 1, 0)+  \frac{1}{2} v(i-1, 1, 0, 1) \Big] + v(i-1, 1, 0, 1)\\
        & \qquad = \frac{ih_0+2h_1}{2\mu} - \frac{ih_0+h_1+h_2}{\mu} + \frac{ih_0+2h_2}{\mu}\\
        & \qquad = \frac{ih_0+2h_2}{2\mu} \geq 0.
\end{align*}
When $i=0$, the result follows similarly using different transitions.
\end{proof}
\begin{proof} [Proof of Statement \ref{state:threshold02} of Monotonicity in Theorem \ref{thm:threshold}]
Such \textcolor{blue}{a} threshold exists if for all $i \geq 0$,
\begin{align}
    \lefteqn{D(i,2,0,0) - D(i,1,0,1)}& \nonumber \\
    & \quad = v(i, 1, 1, 0) - v(i, 1, 0, 1)  - v(i, 0, 1, 1) + v(i, 0, 0, 2) \geq 0. \label{eq:mono02}
\end{align}
\textcolor{blue}{The proof proceeds by considering two cases: $h_1 \leq h_2$ and $h_1 > h_2$.}
\begin{enumerate} [label= \textbf{Case} \arabic*:, leftmargin=2.5\parindent]
    \item \textcolor{blue}{Suppose $h_1 \leq h_2$. Since $\mu_1 = \mu_2 = \mu$, Statement \ref{eq:decision_bd1} in Lemma \ref{lemma:general_rates} together with Proposition \ref{prop:decision_lb} implies that the following holds for all $i \geq 0$:
    \begin{align}
        v(i, 1, 1, 0) - v(i, 1, 0, 1) = \frac{h_1-h_2}{\mu}. \label{eq:i200-equal}
    \end{align}}
    \textcolor{blue}{Another application of Statement \ref{eq:decision_bd1} in Lemma \ref{lemma:general_rates} yields
    \begin{align}
        v(i, 0, 1, 1) - v(i, 0, 0, 2) \leq \frac{h_1-h_2}{\mu}. \label{eq:i101-ub-equal}
    \end{align}}
    \textcolor{blue}{The relations \eqref{eq:i200-equal} and \eqref{eq:i101-ub-equal} together yield the desired result \eqref{eq:mono02}.}
    \item \textcolor{blue}{We now suppose that $h_1 > h_2$ and proceed by induction on $i$. 
    When $i = 0$, a little algebra in the optimality equations \eqref{eq:opt-zero} (recall Statement \ref{state1:general_rates} in Lemma \ref{lemma:general_rates})} yields
    \begin{align}
        \lefteqn{v(0, 1, 1, 0) - v(0, 1, 0, 1) - v(0, 0, 1, 1) + v(0, 0, 0, 2)} \nonumber \\
        & \quad = \left(\frac{h_0}{\mu_0} + \frac{h_1}{\mu} + \textcolor{blue}{\min\left\{\frac{h_1}{\mu}, \frac{h_2}{\mu}\right\}}\right) \nonumber \\
        &\qquad - \left(\frac{h_0}{\mu_0}+\frac{h_2}{\mu} + \frac{\mu_0}{\mu_0+\mu}\min\left\{\frac{h_1}{\mu},\frac{2h_2}{\mu}\right\} + \textcolor{blue}{\frac{\mu_2}{{\mu_0+\mu_2}}\min\left\{\frac{h_1}{\mu}, \frac{h_2}{\mu}\right\}}\right) \nonumber \\
        &\qquad - \left(\frac{h_1}{\mu} + \frac{h_2}{\mu}\right) + \frac{3h_2}{\mu} \nonumber \\
        & \quad \textcolor{blue}{=\frac{h_2}{\mu} +  \frac{\mu_0}{\mu_0+\mu} \left(\min\left\{\frac{h_1}{\mu}, \frac{h_2}{\mu}\right\}-\min\left\{\frac{h_1}{\mu}, \frac{2h_2}{\mu}\right\}\right).} \label{eq:mono02-base}
    \end{align}\textcolor{blue}{For the difference in the second term in \eqref{eq:mono02-base}, choosing $\frac{h_1}{\mu}$ in the second minimum if $\frac{h_1}{\mu}\leq \frac{h_2}{\mu}$ and $\frac{2h_2}{\mu}$ otherwise, we obtain}
    \textcolor{blue}{
    \begin{align}
        &\min\left\{\frac{h_1}{\mu}, \frac{h_2}{\mu}\right\}-\min\left\{\frac{h_1}{\mu}, \frac{2h_2}{\mu}\right\} \geq \min\left\{\frac{h_1}{\mu}-\frac{h_1}{\mu}, \frac{h_2}{\mu}-\frac{2h_2}{\mu}\right\}  = -\frac{h_2}{\mu}. \label{eq:mono02-base-second-term}
    \end{align}Substituting \eqref{eq:mono02-base-second-term} into \eqref{eq:mono02-base} yields the non-negativity of \eqref{eq:mono02-base}, thereby establishing the base case of \eqref{eq:mono02}.}
    We now prove \eqref{eq:mono02} for a generic $i$ by assuming it holds at $i-1$, where $i \geq 1$. 
    \textcolor{blue}{Note that \eqref{eq:i200-diff-equal} and \eqref{eq:i101-diff-equal} hold by straightforward algebraic computation of optimality equations \eqref{eq:opt-zero}--\eqref{eq:opt1}, without additional assumptions. Substituting them into the left-hand side of \eqref{eq:mono02} yields}
\begin{align*}
    \lefteqn{v(i, 1, 1, 0) - v(i, 1, 0, 1) - v(i, 0, 1, 1) + v(i, 0, 0, 2)}\\
        & \qquad = \frac{h_1-h_2}{\mu_0+\mu} - \frac{h_1 - 3h_2 - ih_0}{2\mu} \\
        & \qquad \quad + \frac{\mu_0}{\mu_0+\mu} \Big[ \min\{v(i, 0, 2, 0),v(i, 0, 1, 1)\} - \min\{v(i, 0, 1, 1),v(i, 0, 0, 2)\} \Big]\\
        & \qquad \quad - \frac{1}{2} \Big[ v(i-1, 1, 1, 0) - v(i-1, 1, 0, 1)\Big].
\end{align*}
\textcolor{black}{Recall that $h_1 > h_2$ and $\mu_1 = \mu_2 = \mu$, implying $v(i, 0, 2, 0) - v(i, 0, 1, 1) \geq 0$ by Proposition \ref{prop:decision_lb} (see \eqref{eq:decision_lb_i110}). Choosing an upper bound $ v(i, 0, 0, 2)$ in the second minimum $\min\{v(i, 0, 1, 1),v(i, 0, 0, 2)\}$ yields}
\begin{small}
\begin{align}
    \lefteqn{v(i, 1, 1, 0) - v(i, 1, 0, 1) - v(i, 0, 1, 1) + v(i, 0, 0, 2)}  \nonumber\\
        & \qquad \geq \frac{h_1-h_2}{\mu_0+\mu} - \frac{h_1 - 3h_2 - ih_0}{2\mu} \nonumber \\
        & \qquad \quad  + \frac{\mu_0}{\mu_0+\mu} \Big[ v(i, 0, 1, 1) - v(i, 0, 0, 2) \Big] - \frac{1}{2} \Big[ v(i-1, 1, 1, 0) - v(i-1, 1, 0, 1)\Big] \nonumber \\
        & \qquad = \frac{h_1-h_2}{\mu_0+\mu} - \frac{h_1 - 3h_2 - ih_0}{2\mu} + \frac{\mu_0}{\mu_0+\mu}\left(\frac{h_1 - 3h_2 - ih_0}{2\mu}\right) - \frac{1}{2}\left(\frac{h_1-h_2}{\mu_0+\mu}\right)  \nonumber \\
        & \qquad \quad + \frac{\mu_0}{2(\mu_0+\mu)} \Big[ v(i-1, 1, 1, 0) - v(i-1, 1, 0, 1)  \nonumber \\
        & \qquad \quad - \min\{v(i-1, 0, 2, 0),v(i-1, 0, 1, 1)\}+ \min\{v(i-1, 0, 1, 1),v(i-1, 0, 0, 2)\}\Big]  \nonumber\\
        & \qquad = \frac{\mu_0}{2(\mu_0+\mu)} \Big[\frac{2h_2+ih_0}{\mu_0} + v(i-1, 1, 1, 0) - v(i-1, 1, 0, 1)  \nonumber \\
        & \qquad \quad - \min\{v(i-1, 0, 2, 0),v(i-1, 0, 1, 1)\} + \min\{v(i-1, 0, 1, 1),v(i-1, 0, 0, 2)\}\Big],  \nonumber 
\end{align}
\end{small}\textcolor{blue}{where the first equality follows from substituting \eqref{eq:i101-diff-equal} and \eqref{eq:i200-diff-equal}, with $i$ replaced by $i-1$ in \eqref{eq:i200-diff-equal}.} 
Noting the non-negativity of $\frac{2h_2+ih_0}{\mu_0}$ and replacing the first minimum with $v(i-1,0,1,1)$, we obtain
\begin{align}
    \lefteqn{v(i, 1, 1, 0) - v(i, 1, 0, 1) - v(i, 0, 1, 1) + v(i, 0, 0, 2)}  \nonumber\\
    & \qquad \geq \frac{\mu_0}{2(\mu_0+\mu)} \Big[ v(i-1, 1, 1, 0) - v(i-1, 1, 0, 1) \nonumber \\
    & \qquad \quad - v(i-1, 0, 1, 1) + \min\{v(i-1, 0, 1, 1),v(i-1, 0, 0, 2)\}\Big], \label{eq:mono02-final-diff}
\end{align}
There are now several subcases to consider for the last minimum in \eqref{eq:mono02-final-diff}.
\begin{enumerate}[label= \textbf{Subcase} (\alph*):, leftmargin=2.5\parindent]
    \item If $v(i-1, 0, 1, 1) \leq v(i-1, 0, 0, 2)$ in \eqref{eq:mono02-final-diff}, we have 
    \begin{align*}
        \lefteqn{v(i, 1, 1, 0) - v(i, 1, 0, 1) - v(i, 0, 1, 1) + v(i, 0, 0, 2)}\\
            & \qquad \geq \frac{\mu_0}{2(\mu_0+\mu)} \Big[v(i-1, 1, 1, 0) - v(i-1, 1, 0, 1)\Big] \geq 0,
    \end{align*}
    by \eqref{eq:decision_lb_i200} of Proposition \ref{prop:decision_lb} \textcolor{black}{(since $h_1>h_2$ and $\mu_1 = \mu_2 = \mu$)}.
    \item If $v(i-1, 0, 1, 1) > v(i-1, 0, 0, 2)$ in \eqref{eq:mono02-final-diff},
    \begin{align*}
        \lefteqn{v(i, 1, 1, 0) - v(i, 1, 0, 1) - v(i, 0, 1, 1) + v(i, 0, 0, 2)}\\
        & \qquad \geq \frac{\mu_0}{2(\mu_0+\mu)} \Big[v(i-1, 1, 1, 0) - v(i-1, 1, 0, 1) \\
        & \qquad \quad - v(i-1, 0, 1, 1)+ v(i-1, 0, 0, 2)\Big]\geq 0, 
    \end{align*}
    where the inequality follows from the inductive hypothesis.
\end{enumerate}
\end{enumerate}
\end{proof}
\begin{proof} [Proof of Statement \ref{state:threshold01} of Monotonicity in Theorem \ref{thm:threshold}]
Such a threshold exists if for all $i \geq 0$,
\begin{align}
    \lefteqn{D(i,1,1,0) - D(i,2,0,0)}&\nonumber\\
    & \quad = v(i, 0, 2, 0) - v(i, 0, 1, 1)  - v(i, 1, 1, 0) + v(i, 1, 0, 1) \geq 0. \label{eq:mono01}
\end{align}
    \textcolor{blue}{We proceed by considering two cases: $h_1 \leq h_2$ and $h_1 > h_2$. Suppose $h_1 \leq h_2$ holds. Then Statement \ref{eq:decision_bd1} in Lemma \ref{lemma:general_rates}, together with Proposition \ref{prop:decision_lb}, implies both \eqref{eq:i200-equal} and the following equation for all $i \geq 0$:
    \begin{align}
        v(i, 0, 2, 0) - v(i, 0, 1, 1) = \frac{h_1-h_2}{\mu}. \label{eq:i110-equal}
    \end{align}
    The result follows from substituting \eqref{eq:i110-equal} and \eqref{eq:i200-equal} into the left-hand side of \eqref{eq:mono01}.} 
    \textcolor{blue}{It remains to consider the case $h_1 > h_2$, for which we prove by induction on $i$.} If $i = 0$, recall the arithmetic in proving Statement \ref{state1:general_rates} in Lemma \ref{lemma:general_rates},
    \begin{align*}
        \lefteqn{v(0, 0, 2, 0) - v(0, 0, 1, 1) - v(0, 1, 1, 0) + v(0, 1, 0, 1)}\\
        & \quad =  \frac{2h_1}{\mu} - \left(\frac{h_1}{\mu} + \frac{h_2}{\mu}\right) - \left(\frac{h_0}{\mu_0} + \frac{h_1}{\mu} + \textcolor{blue}{\min\left\{\frac{h_1}{\mu}, \frac{h_2}{\mu}\right\}}\right) \\
        & \qquad + \left(\frac{h_0}{\mu_0}+\frac{h_2}{\mu} + \frac{\mu_0}{\mu_0+\mu}\min\left\{\frac{h_1}{\mu},\frac{2h_2}{\mu}\right\} + \textcolor{blue}{\frac{\mu_2}{{\mu_0+\mu_2}}\min\left\{\frac{h_1}{\mu}, \frac{h_2}{\mu}\right\}}\right)\\
        & \quad = \frac{\mu_0}{\mu_0+\mu}\left( \min\left\{\frac{h_1}{\mu},\frac{2h_2}{\mu}\right\}- \textcolor{blue}{\min\left\{\frac{h_1}{\mu}, \frac{h_2}{\mu}\right\}}\right)\\
        &\quad\geq0,
    \end{align*}
    where the last inequality \textcolor{blue}{holds because
    \begin{align*}
        \min\left\{\frac{h_1}{\mu}, \frac{2h_2}{\mu}\right\} - \min\left\{\frac{h_1}{\mu}, \frac{h_2}{\mu}\right\} \geq \min\left\{\frac{h_1}{\mu}-\frac{h_1}{\mu}, \frac{2h_2}{\mu}-\frac{h_2}{\mu}\right\} = 0.
    \end{align*}
    Suppose \eqref{eq:mono01} holds for $i-1$ and consider $i\geq 1$.}
    \textcolor{blue}{Notice that \eqref{eq:i110-diff-equal} and \eqref{eq:i200-diff-equal} continue to hold without additional assumptions. Applying \eqref{eq:i110-diff-equal} and \eqref{eq:i200-diff-equal} in the left-hand side of \eqref{eq:mono01} yields}
    \begin{align}
        \lefteqn{v(i, 0, 2, 0) - v(i, 0, 1, 1) - v(i, 1, 1, 0) + v(i, 1, 0, 1)} \nonumber \\
        & \quad = \frac{h_1 - h_2}{2\mu}  - \frac{h_1-h_2}{\mu_0+\mu} + \frac{1}{2}\Big[v(i-1, 1, 1, 0) - v(i-1, 1, 0, 1)\Big] \nonumber \\
        & \qquad - \frac{\mu_0}{\mu_0+\mu} \Big[ \min\{v(i, 0, 2, 0),v(i, 0, 1, 1)\} - \min\{v(i, 0, 1, 1),v(i, 0, 0, 2)\}\Big]. \label{eq:mono02-case2}
    \end{align}
    \textcolor{blue}{Since $h_1 > h_2$ and $\mu_1 = \mu_2 = \mu$, inequality \eqref{eq:decision_lb_i110} in Proposition \ref{prop:decision_lb} applies, impying that $v(i-1, 0, 2, 0)\geq v(i-1, 0, 1, 1)$. Hence \eqref{eq:mono02-case2} equals the following}
    \begin{align}
        \lefteqn{v(i, 0, 2, 0) - v(i, 0, 1, 1) - v(i, 1, 1, 0) + v(i, 1, 0, 1)}  \nonumber \\      
        & \qquad \textcolor{blue}{=} \frac{h_1 - h_2}{2\mu}  - \frac{h_1-h_2}{\mu_0+\mu} + \frac{1}{2}\Big[v(i-1, 1, 1, 0) - v(i-1, 1, 0, 1)\Big] \nonumber \\
        & \qquad \quad - \frac{\mu_0}{\mu_0+\mu} \Big[v(i, 0, 1, 1) - \min\{v(i, 0, 1, 1),v(i, 0, 0, 2)\}\Big]. \label{eq:mono02-case2-2}
    \end{align}
    Due to the minimum in \eqref{eq:mono02-case2-2}, there are several cases to consider.
    \begin{enumerate}[label= \textbf{Case} \arabic*:, leftmargin=3\parindent]
        \item If $v(i, 0, 1, 1) < v(i, 0, 0, 2)$, then \eqref{eq:mono02-case2-2} becomes
        \begin{align}
            \lefteqn{v(i, 0, 2, 0) - v(i, 0, 1, 1) - v(i, 1, 1, 0) + v(i, 1, 0, 1)} \nonumber \\
            & \quad \textcolor{blue}{=} \frac{h_1 - h_2}{2\mu}  - \frac{h_1-h_2}{\mu_0+\mu} + \frac{1}{2}\Big[v(i-1, 1, 1, 0) - v(i-1, 1, 0, 1)\Big] \nonumber \\
            & \quad = \frac{h_1 - h_2}{2\mu}  - \frac{h_1-h_2}{\mu_0+\mu} + \frac{1}{2} \bigg(\frac{h_1-h_2}{\mu_0+\mu} \nonumber \\
            & \qquad + \frac{\mu_0}{\mu_0+\mu} \Big[ \min\{v(i-1, 0, 2, 0),v(i-1, 0, 1, 1)\} \nonumber \\
            & \qquad - \min\{v(i-1, 0, 1, 1),v(i-1, 0, 0, 2)\}\Big]\bigg), \label{eq:mono01-case1}
        \end{align}
        \textcolor{blue}{where the last step substitutes the expression in \eqref{eq:i200-diff-equal} with $i$ replaced by $i-1$.} 
        Since $h_1 > h_2$ and $\mu_1 = \mu_2 = \mu$, inequality \eqref{eq:decision_lb_i110} in Proposition \ref{prop:decision_lb} applies, impying that $v(i-1, 0, 2, 0)\geq v(i-1, 0, 1, 1)$ in the first minimum in \eqref{eq:mono01-case1}:
        \begin{align*}
            \lefteqn{v(i, 0, 2, 0) - v(i, 0, 1, 1) - v(i, 1, 1, 0) + v(i, 1, 0, 1)}&\\
            & \quad \textcolor{blue}{=} \frac{h_1 - h_2}{2\mu} - \frac{1}{2} \frac{h_1-h_2}{\mu_0+\mu}\\
            & \qquad + \frac{\mu_0}{2(\mu_0+\mu)} \Big[ v(i-1, 0, 1, 1)- \min\{v(i-1, 0, 1, 1),v(i-1, 0, 0, 2)\}\Big] \\
            & \quad \geq  0,
        \end{align*}
        \textcolor{blue}{where the last equality follows from taking $v(i-1,0,1,1)$ as an upper bound in the minimum.}
        \item Suppose $v(i, 0, 1, 1) \geq v(i, 0, 0, 2)$ in \eqref{eq:mono02-case2-2}:
        \begin{align*}
            \lefteqn{v(i, 0, 2, 0) - v(i, 0, 1, 1) - v(i, 1, 1, 0) + v(i, 1, 0, 1)}&\\
            & \quad \textcolor{blue}{=} \frac{h_1 - h_2}{2\mu}  - \frac{h_1-h_2}{\mu_0+\mu}
            +\frac{1}{2}\Big[v(i-1, 1, 1, 0) - v(i-1, 1, 0, 1)\Big] \\
            & \qquad - \frac{\mu_0}{\mu_0+\mu} \Big[v(i, 0, 1, 1) - v(i, 0, 0, 2)\Big].
        \end{align*}
        \textcolor{blue}{Observe that \eqref{eq:i101-diff-equal} also holds without additional assumptions. Applying \eqref{eq:i200-diff-equal} (evaluated at $i-1$) and \eqref{eq:i101-diff-equal}, we obtain}
        \begin{align*}
            \lefteqn{v(i, 0, 2, 0) - v(i, 0, 1, 1) - v(i, 1, 1, 0) + v(i, 1, 0, 1)}&\\
            & \quad \textcolor{blue}{=} \frac{h_1 - h_2}{2\mu}  - \frac{h_1-h_2}{\mu_0+\mu} + \frac{1}{2} \left(\frac{h_1-h_2}{\mu_0+\mu}\right) - \frac{\mu_0}{\mu_0+\mu}\left(\frac{h_1 - 3h_2 - ih_0}{2\mu}\right)\\
            & \qquad + \frac{\mu_0}{2(\mu_0+\mu)} \Big[ \min\{v(i-1, 0, 2, 0),v(i-1, 0, 1, 1)\} \\
            & \qquad - \min\{v(i-1, 0, 1, 1),v(i-1, 0, 0, 2)\} - v(i-1, 1, 1, 0) + v(i-1, 1, 0, 1)\Big]\\
            & \quad = \frac{\mu_0}{2(\mu_0+\mu)} \Big[ \frac{2h_2+ih_0}{\mu} + \min\{v(i-1, 0, 2, 0),v(i-1, 0, 1, 1)\} \\
            & \qquad - \min\{v(i-1, 0, 1, 1),v(i-1, 0, 0, 2)\} - v(i-1, 1, 1, 0) + v(i-1, 1, 0, 1)\Big].
        \end{align*}
        \textcolor{black}{Another application of \eqref{eq:decision_lb_i110} in Proposition \ref{prop:decision_lb}, under $h_1 > h_2$ and $\mu_1=\mu_2=\mu$, implies $v(i-1, 0, 2, 0) \geq v(i-1, 0, 1, 1)$. Replacing the second minimum with an upper bound $v(i-1, 0, 0, 2)$ yields}
        \begin{align*}
            \lefteqn{v(i, 0, 2, 0) - v(i, 0, 1, 1) - v(i, 1, 1, 0) + v(i, 1, 0, 1)}&\\
            & \quad \geq  \frac{\mu_0}{2(\mu_0+\mu)} \Big[ \frac{2h_2+ih_0}{\mu} + v(i-1, 0, 1, 1) - v(i-1, 0, 0, 2)\\
            & \qquad - v(i-1, 1, 1, 0) + v(i-1, 1, 0, 1)\Big]\\
            &\quad \geq 0,
        \end{align*}
        where the final step follows from \eqref{eq:equal_rates_inte} in Lemma \ref{lemma:equal_rates}.
    \end{enumerate}
This completes the proof of \eqref{eq:mono01}, and the result follows.
\end{proof}

\bibliographystyle{unsrtnat}
\bibliography{policy_ref}
\end{document}